\documentclass[]{article}
\pdfoutput=1
\usepackage[utf8]{inputenc}
\usepackage[english]{babel}
\usepackage{fullpage}
\usepackage{amsmath,amssymb,amsfonts}
\usepackage{graphicx,epstopdf,epsfig}
\usepackage{caption,subcaption}
\usepackage{algorithm,algorithmic}
\usepackage{mathtools}
\mathtoolsset{showonlyrefs=true}
\usepackage{hyperref}
\usepackage{booktabs,multirow,lscape}

\usepackage[affil-it]{authblk}

\usepackage{tikz}

\newcommand{\circlearound}[1]{%
  \begin{tikzpicture}[baseline=(X.base), inner sep=0, outer sep=0]
    \node[draw,circle] (X)  {$#1$};
  \end{tikzpicture}%
}

\usepackage{placeins}

\author{Marie Kub\'{i}nov\'{a}%
  \thanks{Faculty of Mathematics and Physics, Charles University, Prague, Czech Republic. Electronic address: \texttt{kubinova@karlin.mff.cuni.cz}. Research supported in part by the grant SVV-2017-260455.}}
\affil{}
\date{}
\author{James G. Nagy%
  \thanks{Department of Mathematics and Computer Science, Emory University, GA, USA. Electronic address: \texttt{nagy@mathcs.emory.edu}. Research supported in part by US National Science Foundation under grant no.~DMS-1522760.}}
\affil{}

\title{Robust regression for mixed Poisson-Gaussian model}

\begin{document}
\maketitle

\begin{abstract}This paper focuses on efficient computational approaches to compute approximate
solutions of a linear inverse problem that is contaminated with mixed Poisson--Gaussian noise, and 
when there are additional outliers in the measured data. 
The Poisson--Gaussian noise leads to a weighted minimization problem, with solution-dependent
weights.  To address outliers, the standard least squares fit-to-data metric is replaced by the
Talwar robust regression function.  Convexity, regularization parameter selection schemes, and
incorporation of non-negative constraints are investigated.  A projected Newton algorithm is used
to solve the resulting constrained optimization problem, and a preconditioner
is proposed to accelerate conjugate gradient Hessian solves.  Numerical experiments
on problems from image deblurring 
illustrate the effectiveness of the methods.
\end{abstract}



\section{Introduction}\label{sec:intro}
In this paper we consider efficient computational approaches to compute approximate
solutions of a linear inverse problem,
\begin{equation}
b = Ax_\text{true} + \eta, \quad \quad A\in\mathbb{R}^{m\times n},\label{eq:inverse_problem}
\end{equation}
where $A$ is a known matrix, 
vector $b$ represents known acquired data, $\eta$ represents noise, and vector $x_\text{true}$ 
represents the unknown quantity that needs to be approximated.
We are particularly interested in imaging applications where $x_\text{true} \geq 0$ and $Ax_\text{true} \geq 0$.
Although this basic problem has been studied extensively 
(see, for example, \cite{Engl2000Regularization,Hansen2010Discrete,Mueller2012Linear,Vogel2002Computational}
and the references therein), the noise is typically assumed to come from a single source (or to be
represented by a single statistical distribution) and the data to contain no outliers.
In this paper we focus on a practical situation that arises in many imaging applications, and for which relatively
little work has been done, namely when the noise is comprised of a mixture of Poisson and
Gaussian components {\em and}
when there are outliers in the measured data.  While some research has been done
on the two topics separately (i.e., mixed Poisson--Gaussian noise models {\em or} 
outliers in measured data), to our knowledge no work has been done when the measured
data contains both issues. In the following, we review some of the approaches used to handle each of the issues.
\subsubsection*{Poisson--Gaussian noise}
A Poisson--Gaussian statistical model for the measured data takes the form
\begin{equation}
b_i = n_\text{obj}(i) + g(i), \ i = 1,\ldots,m, \quad \ n_\text{obj}(i) \sim \text{Pois}([Ax_\text{true}]_i), \ g(i) \sim \mathcal{N}(0,\sigma^2), \label{eq:noise}
\end{equation}
where $b_i$ is the $i$th component of the vector $b$ and $[Ax_\text{true}]_i$ the $i$th component of the true noise-free 
data $Ax_\text{true}$. We assume that the two random variables $n_\text{obj}(i)$ and $g(i)$ are independent. 
This mixed noise model arises in many important applications, such as when using charged coupled device (CCD) arrays,
x-ray detectors, and infrared 
photometers \cite{Bardsley2003nonnegatively,Gravel2004method,Luisier2011Image,Makitalo2013Optimal,Snyder1993Image}.
The Poisson part (sometimes referred to as shot noise) can arise
from the accumulation of photons over a detector, and the Gaussian part usually is due to {\em read-out} noise from
a detector, which can be generated by thermal fluctuations in interconnected electronics.

Since the log-likelihood function for the mixed Poisson--Gaussian model (\ref{eq:noise}) has an infinite series representation \cite{Snyder1993Image}, we assume a simplified model, where both random variables have the same type of distribution. 
There are two main approaches one can take to generate a simplified model.
The first approach is to add $\sigma^2$ to each component of the vector $b$, and from (\ref{eq:noise}) it then follows that
\begin{equation}
\mathbb{E}(b_i + \sigma^2) = [Ax_\text{true}]_i + \sigma^2 \quad \text{and} \quad \text{var}(b_i + \sigma^2) = [Ax_\text{true}]_i + \sigma^2. \label{eq:poiss_approx_noise}
\end{equation}
For large $\sigma$, the Gaussian random variable $g(i) + \sigma^2$ is well-approximated by a Poisson random variable with the Poisson parameter $\sigma^2$, and therefore $b_i + \sigma^2$ is also well approximated by a Poisson random variable with the Poisson parameter $[Ax_\text{true}]_i + \sigma^2$. The data fidelity function corresponding to the negative Poisson log-likelihood then has the form
\begin{equation}
\sum_{i=1}^m ([Ax_\text{true}]_i + \sigma^2) - (b_i +\sigma^2)\,\log([Ax_\text{true}]_i + \sigma^2);\label{eq:poisson}
\end{equation}
see also \cite{Snyder1993Image}. An alternative approach is to approximate the true negative log-likelihood by a weighted least-squares function, where the weights correspond to the measured data, i.e.,
\begin{equation}
 \sum_{i=1}^m \frac{1}{2}\left(\frac{[Ax]_i - b_i}{\sqrt{b_i + \sigma^2}}\right)^2;\label{eq:wls_basic}
\end{equation}
see \cite[Sec. 1.3]{Hansen2013Least}. 
A more accurate approximation can be achieved by replacing the measured data by the computed data 
(which depends on $x$), i.e., replace the fidelity function (\ref{eq:wls_basic}) by
\begin{equation}
 \sum_{i=1}^m \frac{1}{2}\left(\frac{[Ax]_i - b_i}{\sqrt{[Ax]_i + \sigma^2}}\right)^2;\label{eq:wls}
\end{equation}
see \cite{Stagliano2011Analysis} for more details. 
Additional additive Poisson noise (e.g., background emission) can be incorporated into the model in a straightforward way. 
\subsubsection*{Outliers}
For  data corrupted solely with Gaussian noise, i.e., 
\begin{equation}
b_i = [Ax]_i + g(i), \ i = 1,\ldots,m, \quad g(i) \sim \mathcal{N}(0,\sigma^2), \label{eq:gausssian_noise}
\end{equation} employing the negative log-likelihood leads to the standard least-squares functional
\begin{equation}
 \sum_{i=1}^m \frac{1}{2}\left([Ax]_i - b_i\right)^2.\label{eq:ls}
\end{equation} 
It is well known however that a computed solution based on least squares is not robust if outliers occur, meaning that even a small number of components with gross errors can cause a severe deterioration of our estimate. Robustness of the least squares fidelity function can be achieved by replacing the loss function $\frac{1}{2}t^2$ used in \eqref{eq:ls} by a function $\rho(t)$  as
\begin{equation}
 \sum_{i=1}^m \rho\left([Ax]_i - b_i\right),\label{eq:robust}
\end{equation}where the function $\rho$ is less stringent towards the gross errors and satisfies the following conditions:
\begin{enumerate}
\item $\rho(t) \geq 0$;
\item $\rho(t) = 0 \Leftrightarrow t =0$;
\item $\rho(-t) = \rho(t)$;
\item $\rho(t^\prime) \geq \rho(t)$, for $t^\prime \geq t \geq 0$;
\end{enumerate}
see also \cite[Sec. 1.5]{Hansen2013Least}. A list of eight most commonly used loss functions $\rho$ can be found in  \cite{Coleman1980system} or in MATLAB under \texttt{robustfit}; some of them are discussed in Section~\ref{sec:conv_anal}. Each of these functions also depends on a parameter $\beta$ (see Section~\ref{sec:beta}) defining the trade-off between the robustness and efficiency. Note that if we use this robust regression approach, in order to reduce the influence of possible outliers, we always sacrifice some efficiency at the model.

\bigskip

In this paper, we focus on combining these two approaches to suppress the influence of outliers for data with mixed noise \eqref{eq:noise}. Our work has been motivated by O'Leary \cite{OLeary1990Robust}, and more recent work by Calef \cite{Calef2013Iteratively}.  The initial ideas of our work were first outlined in the 
conference paper \cite{Kubinova2015Iteratively}.

The paper is organized as follows. In Section~\ref{sec:robust} we introduce a data-fidelity function suitable for data corrupted both with mixed Poisson--Gaussian noise and outliers. In Section~\ref{sec:reg_param} we propose a regularization parameter choice method for the regularization of the resulting inverse problem, and in Section~\ref{sec:optim} we focus on the optimization algorithm and the solution of the linear subproblems. Section~\ref{sec:num_exp} demonstrates the performance of the resulting method on image deblurring problems with various types of outliers.

Throughout the paper, $D$ (or $D$ with an accent) denotes a general real diagonal matrix, $e_i$ denotes the $i$th column of the identity matrix of a suitable size.

\section{Data-fidelity function}\label{sec:robust}
In Section \ref{sec:intro}, we reviewed fidelity functions \eqref{eq:poisson}, \eqref{eq:wls_basic}, and \eqref{eq:wls}, commonly used for problems with mixed Poisson--Gaussian noise and also robust loss  functions used to handle problems with Gaussian noise and outliers \eqref{eq:robust}. Since we need to deal with both issues simultaneously here, we propose combining both approaches. More specifically, combining a robust loss function with the weighted least squares problem \eqref{eq:wls}, so that the data fidelity function becomes
\begin{equation}
 J(x) = \sum_{i=1}^m \rho\left(\frac{[Ax]_i - b_i}{\sqrt{[Ax]_i + \sigma^2}}\right).
\label{eq:robust_wls}
\end{equation}
In the remainder of this section, we investigate the properties of the proposed data-fidelity function \eqref{eq:robust_wls} and the choice of the robustness parameter $\beta$, which is defined in the next subsection.
\subsection{Choice of the loss function -- convexity analysis} \label{sec:conv_anal}
For ordinary least squares, functions known under names Huber, logistic, Fair, and Talwar, shown in Figure~\ref{fig:loss_functions}, lead to an interval-wise convex data fidelity function, see \cite{OLeary1990Robust}, i.e., positive-semidefinite Hessian, which is favorable for Newton-type minimization algorithms. This however does not 
always hold in our
case where the weighted least squares formulation \eqref{eq:robust_wls} has solution-dependent weights. 

To see this, let us begin by denoting the components of the residual as $r_i \equiv [Ax]_i - b_i$ and the 
solution-dependent weights as $w_i \equiv \frac{1}{\sqrt{[Ax]_i + \sigma^2}}$. Then the gradient and the Hessian of \eqref{eq:robust_wls} can be rewritten as

\begin{align}
\text{grad}_J(x) &= A^Tz, & z_i &= \left(w^\prime_ir_i + w_i\right)\rho^{\prime}(w_ir_i);\\
\text{Hess}_J(x) &= A^TDA, & \quad D_{ii} &= (w^{\prime\prime}_ir_i + 2w_i^{\prime})\rho^{\prime}(w_ir_i) + (w^{\prime}_i r_i + w_i)^2\rho^{\prime\prime}(w_ir_i). \label{eq:hess_and_grad}
\end{align}
We investigate the entries $D_{ii}$ in order to examine the positive semi-definiteness of the Hessian $\text{Hess}_J(x)$. Recall that $A^TDA$ is positive semi-definite, if $D_{ii}\geq 0$.

Assuming $\rho^{\prime\prime}\geq 0$, the signs of the diagonal entries $D_{ii}$ in \eqref{eq:hess_and_grad} are
\begin{equation}
\text{sign}(D_{ii}) = \left(\circlearound{+}\cdot\text{sign}(r_i) + 2\circlearound{-}\right)\cdot\text{sign}(\rho^{\prime}(w_ir_i)) + \circlearound{+}\circlearound{+}, \label{eq:dii_sign}
\end{equation}
where we have replaced some of the quantities in the
expression for $D_{ii}$ shown in equation \eqref{eq:hess_and_grad} with the symbol $\circlearound{-}$ when the
value it replaces is always a negative number and with $\circlearound{+}$ when the value it replaces is always nonnegative. We will now investigate all possible cases with respect to $\text{sign}(\rho^{\prime}(w_ir_i))$:
\begin{itemize}
\item \textit{Case 1: $\rho^\prime(w_ir_i)< 0$} \\
$\rho^\prime(w_ir_i)< 0$  yields $r_i< 0$, and therefore $D_{ii} > 0$.
\item \textit{Case 2: $\rho^\prime(w_ir_i)= 0$} \\
$\rho^{\prime}(w_ir_i) = 0 $ yields $D_{ii} = 0$.
\item \textit{Case 3: $\rho^\prime(w_ir_i)> 0$} \\
 Substituting for $w_i$ and $r_i$ in \eqref{eq:hess_and_grad}, we obtain
\begin{align*}\label{eq:dii}
 D_{ii} &= \left(\frac{3}{4}([Ax]_i-b_i)([Ax]_i+\sigma^2)^{-5/2} - ([Ax]_i+\sigma^2)^{-3/2})\right)\rho^{\prime}(w_ir_i)\\ 
\nonumber & \quad + \left(-\frac{1}{2}([Ax]_i-b_i)([Ax]_i+\sigma^2)^{-3/2} + ([Ax]_i+\sigma^2)^{-1/2}\right)^2\rho^{\prime\prime}(w_ir_i).
\end{align*} 
For $[Ax]_i\gg b_i+\sigma^2$, to achieve $D_{ii}\geq 0$, 
\[
\sqrt{[Ax]_i}\cdot\rho^{\prime\prime}(\sqrt{[Ax]_i})\gtrsim \rho^{\prime}(\sqrt{[Ax]_i}),\\
\]
must hold. This corresponds to 
\[
\rho^{\prime}(t) \gtrsim t \quad \text{yielding} \quad  \rho(t) \gtrsim t^2/2,\\
\]
i.e., for large $[Ax]_i$, the loss function $\rho$ has to grow at least quadratically.
\end{itemize}

Concluding, for large $t$, the loss function $\rho(t)$ has to be either constant or grow at least quadratically, which is in contradiction with the idea of robust regression. Therefore, considering the functions from \cite{Coleman1980system}, the only loss function $\rho$ for which the data fidelity function (\ref{eq:robust_wls}) has positive semidefinite Hessian, is the function Talwar:
\begin{equation}
\rho(t) = \left\{\begin{array}{ll}
t^2/2, & |t|\leq\beta,\\
\beta^2/2, & |t| > \beta.
\end{array} \right.\label{eq:talwar}
\end{equation} 
\subsection{Selection of the robustness parameter}\label{sec:beta}
\begin{figure}
\centering
\begin{subfigure}[b]{.4\textwidth}
\centering
\includegraphics[width = .8\textwidth]{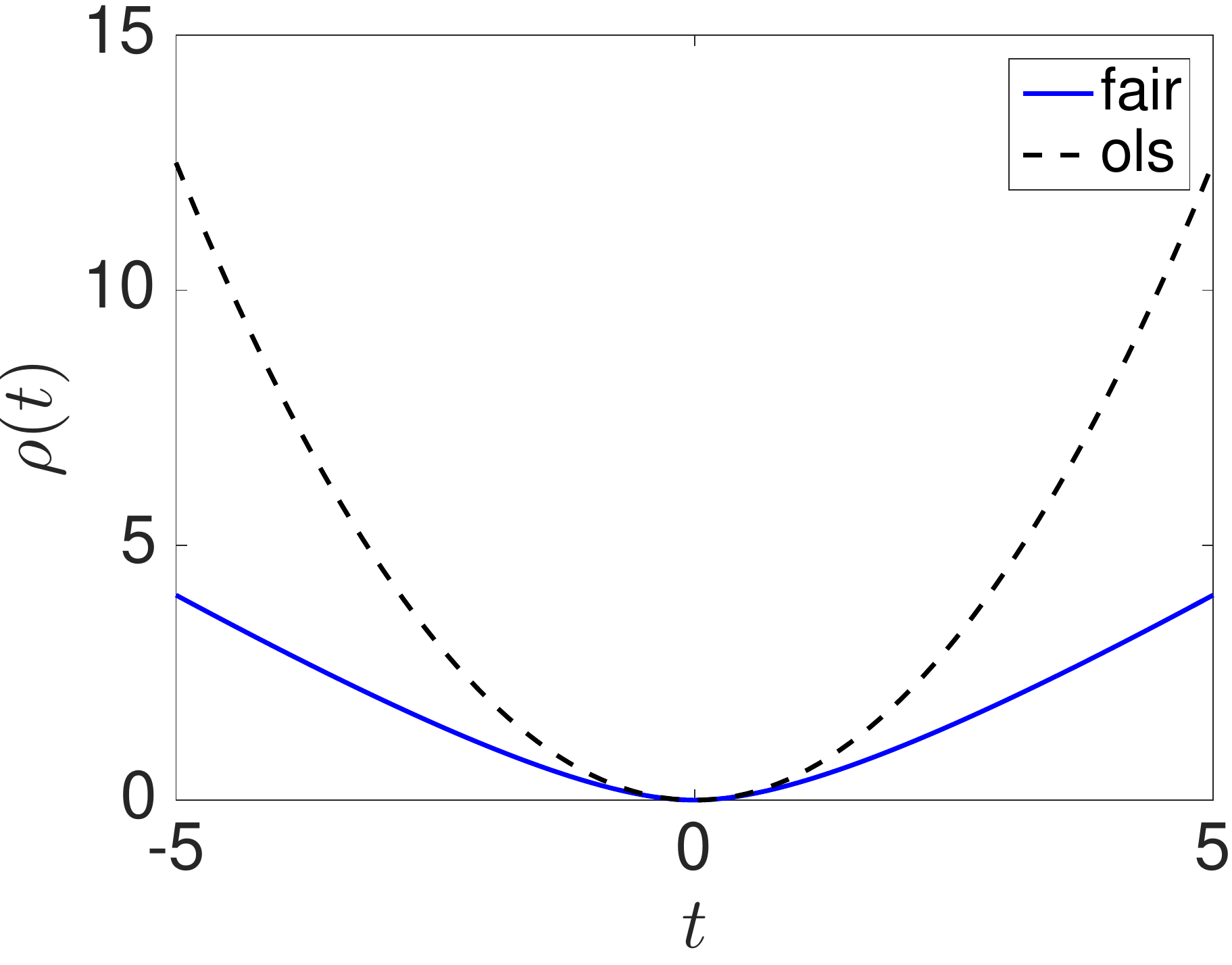} 
\caption{Fair}
\end{subfigure}
\begin{subfigure}[b]{.4\textwidth}
\centering
\includegraphics[width = .8\textwidth]{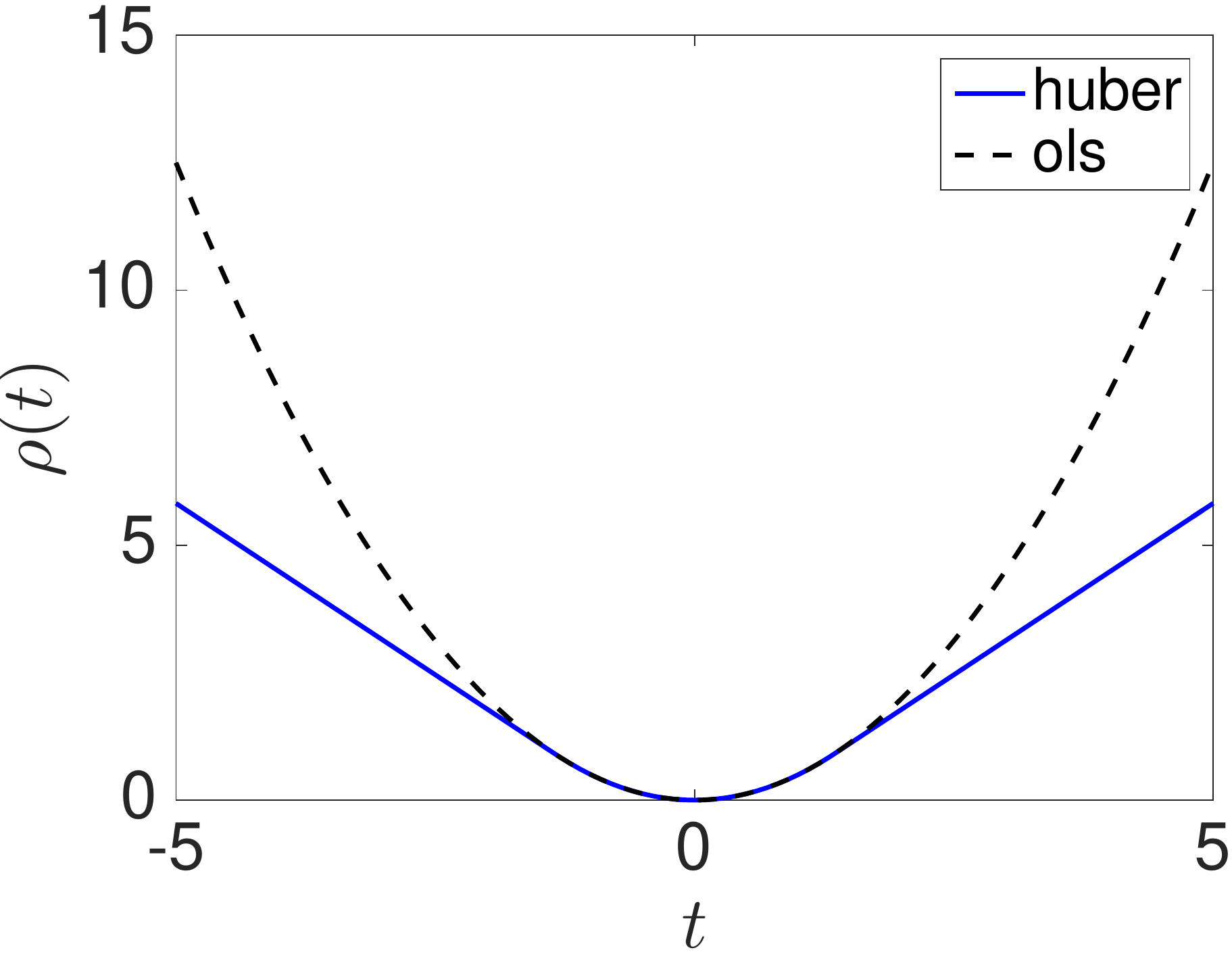} 
\caption{Huber}
\end{subfigure}

\begin{subfigure}[b]{.4\textwidth}
\centering
\includegraphics[width = .8\textwidth]{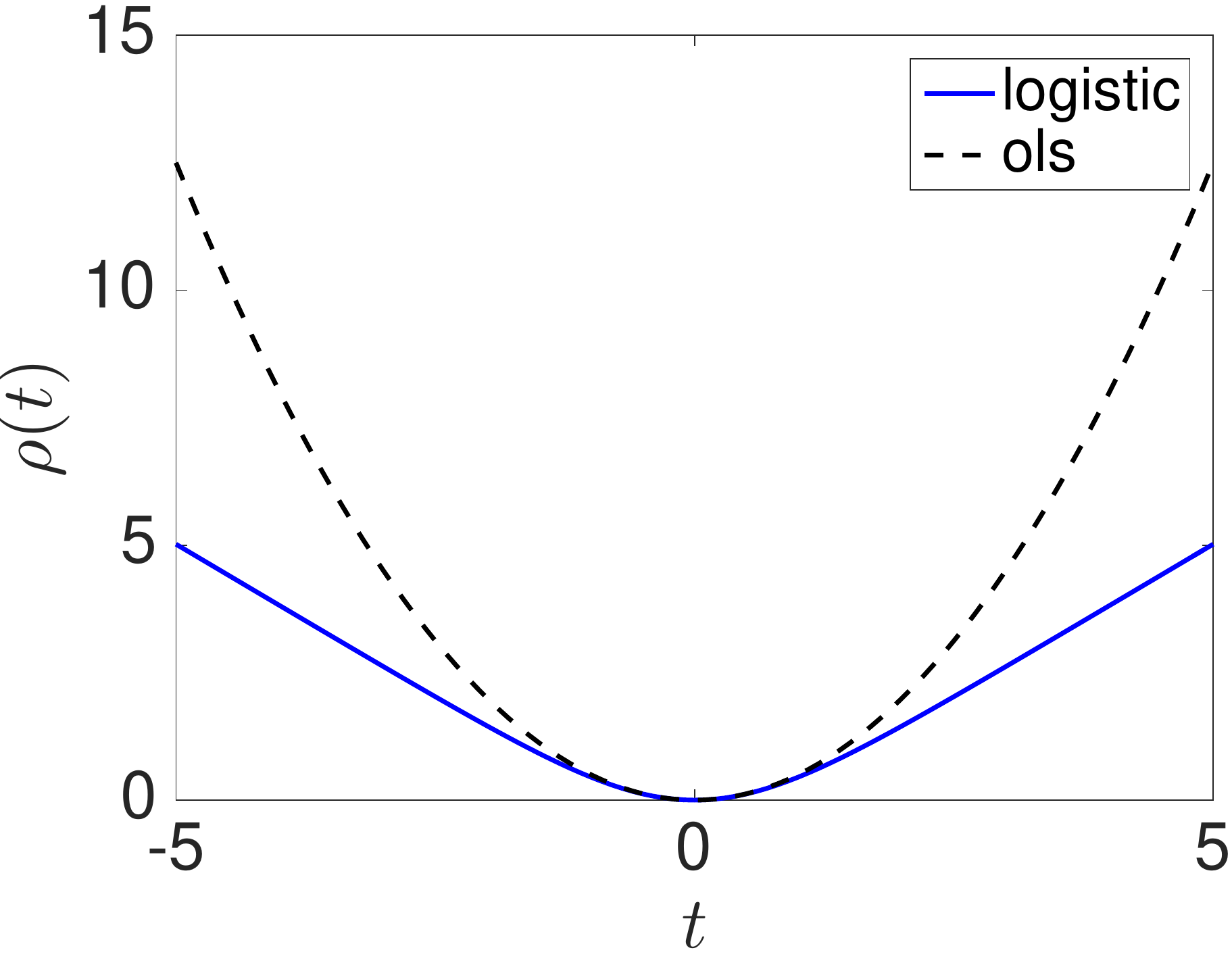} 
\caption{logistic}
\end{subfigure}
\begin{subfigure}[b]{.4\textwidth}
\centering
\includegraphics[width = .8\textwidth]{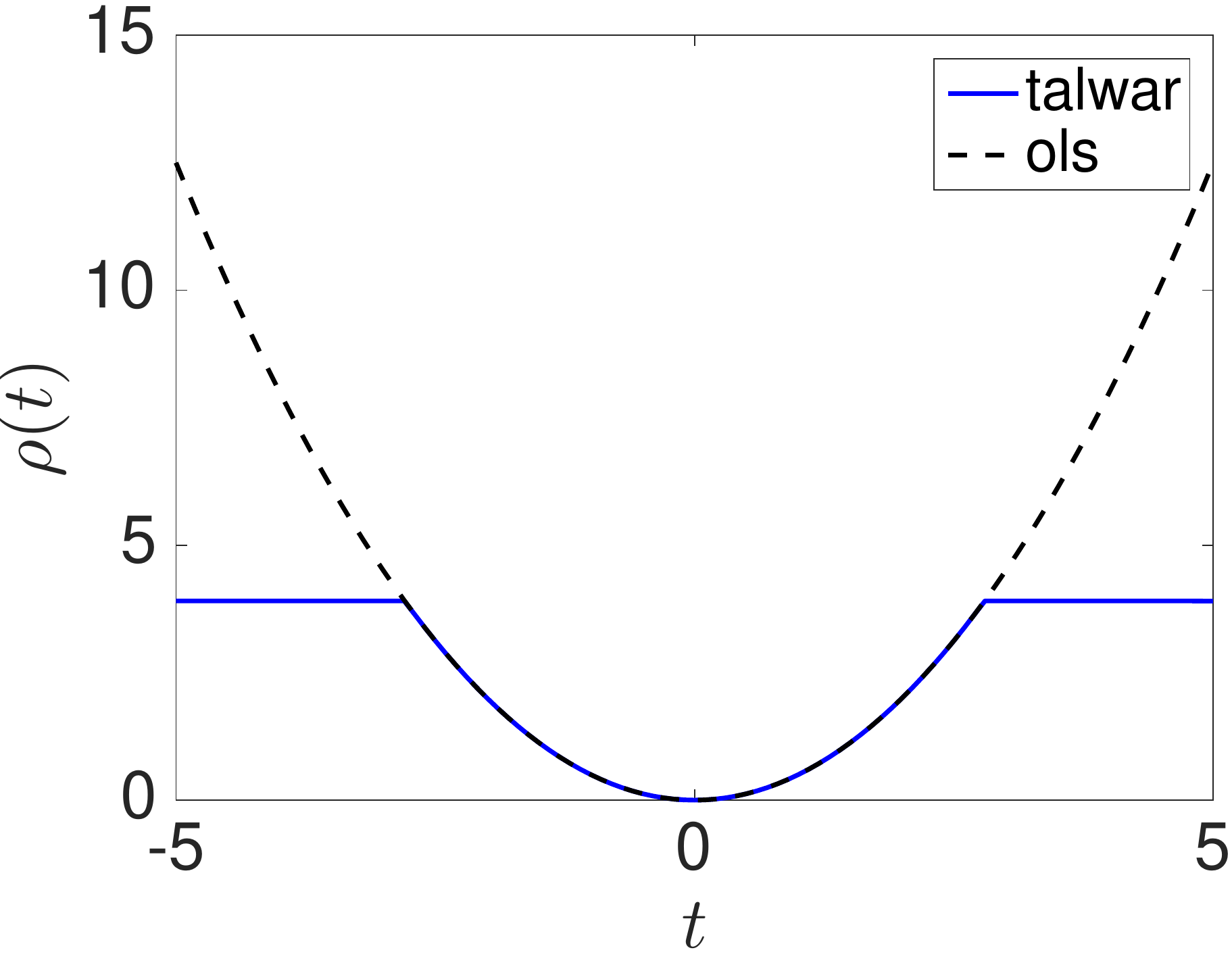} 
\caption{Talwar}\label{fig:rho_talwar}
\end{subfigure}

\caption{Loss functions Fair, Huber, logistic and Talwar for the tuning parameter $\beta$ corresponding to 95\% efficiency (solid line)  together with the standard loss function $t^2/2$ (dashed line). 
}\label{fig:loss_functions}
\end{figure}

Parameters $\beta$ for $95\%$ asymptotic efficiency with respect to the standard loss function $\frac{1}{2}t^2$ when the disturbances come from the unit normal distribution can again be found in \cite{Coleman1980system}. For Talwar, the 95\% efficiency tuning parameter is 
\begin{equation}
\beta_{95} = 2.795.\label{eq:beta_opt}
\end{equation} 
Note that in our specific case, the random variable inside the function $\rho$ in (\ref{eq:robust_wls}) is already rescaled to have unit variance and therefore approximately unit normal distribution. We may therefore apply the parameter $\beta_{95}$ without any further rescaling based on estimated variance, which is usually required in case of ordinary least squares with unknown variance of noise. Function Talwar with  $\beta = \beta_{95}$ is shown in Figure~\ref{fig:rho_talwar}. 

\subsection{Non-negativity constraints}
In many applications, such as imaging, the reconstruction will benefit from taking into account the prior information about the component-wise non-negativity of the true solution $x_\text{true}$. Here, however, imposing non-negativity is not just a question of visual appeal, it also guarantees the two estimates (\ref{eq:poisson}) and (\ref{eq:wls}) of the negative log-likelihood will provide similar results; see \cite{Stagliano2011Analysis}. Therefore, the component-wise non-negativity constraint is an integral part of the resulting optimization problem. However, employment of the non-negativity constraint results in the need of more sophisticated optimization tools. The use of one of the possible algorithms is discussed in Section~\ref{sec:optim}.

\section{Regularization and selection of the regularization parameter}\label{sec:reg_param}
As a consequence of noise and ill-posedness of the inverse problem (\ref{eq:inverse_problem}), some form of regularization needs to be employed in order to achieve a reasonable approximation of the true solution $x_\text{true}$. For computational convenience, we use Tikhonov regularization with a quadratic penalization term, i.e., we minimize the functional of the form
\begin{equation}
J_\lambda(x) \equiv \sum_{i=1}^m \rho\left(\frac{[Ax]_i - b_i}{\sqrt{[Ax]_i + \sigma^2}}\right) + \frac{\lambda}{2}\|Lx\|^2, \qquad x\geq 0.\label{eq:functional}
\end{equation}
We assume that a good regularization parameter $\lambda$ with respect to $L$ is used so that the penalty term is reasonably close to the prior and the residual therefore is close to noise. In case of robust regression, it is particularly important not to over-regularize, since this would lead to large residuals and too many components of the data $b$ would be considered outliers.  Methods for
choosing $\lambda$ are discussed in this section.

\subsection{Morozov's discrepancy principle}
Since the residual components are scaled, and for data without outliers we have the expected value
 \begin{equation}
E\left\{\frac{1}{n}\sum_{i=1}^{n}\frac{([Ax]_i - b_i)^2}{[Ax]_i + \sigma^2}\right\} = 1,\label{eq:discrepancy_estimate}
\end{equation}
an obvious choice would be to use the Morozov's discrepancy principle \cite{Morozov1966solution,Vogel2002Computational}. However, as reported in \cite{Stagliano2011Analysis}, even without outliers, the discrepancy principle based on \eqref{eq:discrepancy_estimate} tends to provide unsatisfactory reconstructions for problems with large signal-to-noise ratio. Therefore we will not consider this approach further.

\subsection{Generalized cross validation}\label{sec:gcv}
The generalized cross validation method \cite{Golub1979Generalized}\cite[chap. 7]{Vogel2002Computational} is a method derived from the standard leave-one-out cross validation. To apply this method for linear Tikhonov regularization, one selects the regularization parameter $\lambda$ such that it minimizes the GCV functional
\begin{equation}
\text{GCV}(\lambda) = \frac{n\|r_\lambda\|^2}{(\text{trace}(I-A_\lambda))^2},\label{eq:gcv_fun}
\end{equation}
where $r_\lambda = Ax_\lambda - b = (A_\lambda - I)b$ is the residual, $n$ is its length, and the influence matrix $A_\lambda$ takes the form $A_\lambda = A(A^TA + \lambda L^TL)^{-1}A^T$. Here, due to the non-negativity constraints and the weights, the residual and the influence matrix have a more complicated form. An approximation of the influence matrix for problems with mixed noise, but without outliers, has been proposed in \cite{Bardsley2009Regularization}. There the numerator of the GCV functional
takes the form  $n\|Wr_\lambda\|^2$ and the approximate influence matrix
\begin{equation}
A_\lambda = WA(D_\lambda(A^TW^2A + \lambda L^TL)D_\lambda)^{\dagger}D_\lambda A^TW,\label{eq:gcv_fun_new} 
\end{equation}
where $W$ and $D_\lambda$ are diagonal matrices:
\begin{align}
W_{ii} &=
\frac{1}{\sqrt{[Ax_\lambda]_i + \sigma^2}};
\label{eq:W_ii_old}\\
[D_\lambda]_{ii} &= \left\{
\begin{array}{cc}
1 & [x_\lambda]_i > 0, \\ 
0 &\text{otherwise,}
\end{array} \right. \nonumber
\end{align}
and \,${}^\dagger$ denotes the Moore-Penrose pseudoinverse. Matrix $D_\lambda$ only handles the non-negativity constraints, 
and therefore can be adopted directly. The matrix $W$ needs a special adjustment, due to the change of the loss function to Talwar. The aim is to construct a matrix $W$ satisfying
\[\|Wr_\lambda\|^2 = \sum_{i=1}^m \rho\left(\frac{[Ax_\lambda]_i - b_i}{\sqrt{[Ax_\lambda]_i + \sigma^2}}\right).\]
Substituting for $\rho$ from the definition of the function Talwar \eqref{eq:talwar}, we redefine the scaling matrix as
\begin{align}
W_{ii} &\equiv \left\{
\begin{array}{cc}
\frac{1}{\sqrt{[Ax_\lambda]_i + \sigma^2}} & \quad \left\vert\frac{[Ax_\lambda]_i - b_i}{\sqrt{[Ax_\lambda]_i + \sigma^2}}\right\vert\leq \beta, \\ 
\frac{\beta}{[Ax_\lambda]_i - b_i} &\text{otherwise};
\end{array} \right.\label{eq:W_ii}
\end{align}

In order to make the evaluation of \eqref{eq:gcv_fun_new} feasible for large-scale problems, we approximate the trace of a matrix using the random trace estimation \cite{Hutchinson1990stochastic,Vogel2002Computational} as $\text{trace}(M)\approx v^TMv$, where the entries of $v$ take values $\pm 1$ with equal probability. Applying the random trace estimation to \eqref{eq:gcv_fun_new}, we obtain
\[
(\text{trace}(I-A_\lambda))^2 \approx (v^Tv - v^TA_\lambda v)^2.
\]
Finally, $A_\lambda v$ is approximated by $WAy$, with $y$ obtained applying truncated conjugate gradient iteration to
\begin{equation}
(D_\lambda(A^TW^2A + \lambda L^TL)D_\lambda)y = D_\lambda A^TWv.\label{eq:gcv_lin_syst}
\end{equation}

\section{Minimization problem}\label{sec:optim}

In this section we discuss numerical methods to compute a minimum of \eqref{eq:functional}.
We consider incorporating a non-negative constraint and solution of linear subproblems, including
proposing a preconitioner.  

\subsection{Projected Newton's method}
Various methods for constrained optimization have been developed over the years, some related to image deblurring can be found in \cite{Bardsley2003nonnegatively,Bonettini2009scaled,Hanke2000Quasi,More1991solution,Nagy2000Enforcing}. 
For our computations, we chose a projected Newton's method\footnote{In \cite{Haber2015Computational}, the method was derived as the Projected Gauss--Newton method. Here, since the evaluation of the Hessian does not represent a computational difficulty, we use it as a variant of Newton's method. Therefore, in the remainder of the text, the method is referred to as the Projected Newton's Method.}, combined with projected PCG to compute the search direction in each step, see \cite[sec. 6.4]{Haber2015Computational}. The convenience of this method lies in the fact that the projected PCG does not require any special form of the preconditioner and a generic conjugate gradient preconditioner can be employed. Besides lower bounds, upper bounds on the reconstruction can also be enforced. For completeness, we include the projected Newton method in Algorithm~\ref{alg:projGNCG}, and projected PCG in
Algorithm~\ref{alg:projPCG}.

\begin{algorithm}
\caption{Projected Newton's method \cite{Haber2015Computational}}
\label{alg:projGNCG}
\begin{algorithmic}
\STATE{$k = 0$}
\WHILE{not converged} 
\STATE{$\text{Active} = (x^{(k)} \leq 0$)} 
\STATE{$g = \text{grad}_{J_\lambda}(x^{(k)})$}
\STATE{$H = \text{Hess}_{J_\lambda}(x^{(k)})$}
\STATE{$M = \texttt{prec}(H)$} \COMMENT{setup preconditioner for the Hessian}
\STATE{$s = \texttt{projPCG}(H,-g,\text{Active},M)$} \COMMENT{compute the search direction for inactive cells}\STATE{$g_a = g(\text{Active})$} 
\IF{$\max(\text{abs}(g_a)) >  \max(\text{abs}(s))$} 
\STATE {$g_a = g_a\cdot\max(\text{abs}(s))/\max(\text{abs}(g_a))$} \COMMENT{rescaling needed} \ENDIF
\STATE{$s(\text{Active}) = g_a$} \COMMENT{take gradient direction in active cells}
\STATE{$x^{(k+1)} = \texttt{linesearch}(s,x^{(k)},J_\lambda,\text{grad}_{J_\lambda})$}
\STATE{$k = k+1$}
\ENDWHILE
\RETURN $x^{(k)}$
\end{algorithmic}
\end{algorithm}

\begin{algorithm}[H]
\caption{Projected PCG \cite{Haber2015Computational}}
\label{alg:projPCG}
\begin{algorithmic}
\STATE{input: $A$, $b$, Active, $M$}
\STATE{$x_0 = 0$}
\STATE{$D_\mathcal{I} = \text{diag}(1-\text{Active})$} \COMMENT{projection onto inactive set}
\STATE{$r_0 = D_\mathcal{I}b$}
\STATE{$z_0 = D_\mathcal{I}(M^{-1}r_0)$}
\STATE{$p_0 = z_0$}
\STATE{$k = 0$}
\WHILE{not converged} 
\STATE{$\alpha_k = \frac{r_k^Tz_k}{p^TD_\mathcal{I}Ap_k}$} 
\STATE{$x_{k+1} = x_k + \alpha_kp_k $}
\STATE{$r_{k+1} = x_k - \alpha_kD_\mathcal{I}Ap_k $}
\STATE{$z_{k+1} = D_\mathcal{I}(M^{-1}r_k)$}
\STATE{$\beta_{k+1} = \frac{z_{k+1}^Tr_{k+1}}{z_{k}^Tr_{k}}$}
\STATE{$p_{k+1} = z_{k+1} + \beta_kp_k$}
\STATE{$k = k+1$}
\ENDWHILE
\RETURN $x_{k}$
\end{algorithmic}
\end{algorithm}

\subsection{Solution of the linear subproblems}\label{sec:lin_systems}
Each step of the projected Newton method (Algorithm \ref{alg:projGNCG}) 
requires solving a linear system with the Hessian:
\begin{align}
\text{Hess}_{J_\lambda}(x^{(k)})s &= -\text{grad}_{J_\lambda}(x^{(k)})\nonumber \\
(A^TD^{(k)}A + \lambda L^TL)s &= -\left(A^Tz^{(k)} + \lambda L^TLx^{(k)}\right).\label{eq:lin_system}
\end{align}
For the objective functional \eqref{eq:functional}, the diagonal matrix $D^{(k)}$ and the vector $z^{(k)}$ have the form:
\begin{align}
z_{i} &= \left\{\begin{array}{ll}
\frac{1}{2}\left(1 - \frac{\left(b_i + \sigma^2\right)^2}{\left([Ax]_i + \sigma^2\right)^2}\right), & \left\vert\frac{[Ax]_i - b_i}{\sqrt{[Ax]_i + \sigma^2}}\right\vert\leq \beta,\\
0, & \text{otherwise}.
\end{array}
\right.\\
D_{ii} &= \left\{\begin{array}{ll}
\frac{\left(b_i + \sigma^2\right)^2}{\left([Ax]_i + \sigma^2\right)^3}, & \left\vert\frac{[Ax]_i - b_i}{\sqrt{[Ax]_i + \sigma^2}}\right\vert\leq \beta,\\
0, & \text{otherwise}.\label{eq:row_scaling}
\end{array}
\right.
\end{align}
Note that in case of constant weights, robust regression represents extra computational cost in comparison with standard least squares because it leads to a sequence of weighted least squares problems, while standard least squares problems are solved in one step. In our setting, the weights in \eqref{eq:wls} themselves have to be updated and therefore employing a different loss function does not change the type of the problem we need to solve. 

Without preconditioning, the convergence of projected PCG can be rather slow, and
it is therefore important to consider preconditioning.  The idea of many preconditioners, such as constraint \cite{Keller2000Constraint,Dollar2007Using}, constraint-type \cite{Dollar2007Constraint} or Hermitian and skew-Hermitian \cite{Benzi2006Preconditioned} preconditioners is based on the fact that in many cases it is possible to efficiently solve the 
linear system in \eqref{eq:lin_system} if the diagonal matrix $D^{(k)}$ is the identity matrix; that is, if the linear system
involves the matrix
\begin{equation}
A^TA + \lambda L^TL \label{eq:without_D}.
\end{equation}
For example, in the case of image deblurring, it is well known that linear systems involving the matrix
\eqref{eq:without_D} can be solved efficiently using fast trigonometric or 
fast Fourier transforms (FFT).

Although the constraint-type, and Hermitian and skew-Hermitian preconditioners seem to perform well for problems with a random 
matrix $D^{(k)}$ (i.e., a random row scaling), see \cite{Benzi2006Preconditioned}, they performed unsatisfactorily for problems of the form \eqref{eq:lin_system}, \eqref{eq:row_scaling}.

A preconditioner based on a similar idea of fast computations with matrices of type \eqref{eq:without_D} for imaging problems was proposed in \cite{Fessler1999Conjugate}. In this case, the row scaling is approximated by a column scaling; that is, 
we find $\hat{D}^{(k)}$ such that 
\begin{equation}
A^TD^{(k)}A \approx \hat{D}^{(k)}(A^TA)\hat{D}^{(k)},\label{eq:pullout_prec}
\end{equation}
where
\begin{equation}
\hat{D}_{ii}^{(k)} \equiv \sqrt{\frac{e_i^T(A^TD^{(k)}A)e_i}{e_i^T(A^TA)e_i}}.\label{eq:out_diag_def}
\end{equation}
Note that for $\hat{D}^{(k)}$ defined in \eqref{eq:out_diag_def}, the diagonals of the matrices on the two sides of approximation \eqref{eq:pullout_prec} are exactly equal. 

Since for large-scale problems, matrix $A$ is typically not formed explicitly, exact evaluation of the entries of $\hat{D}^{(k)}$ might become too expensive. To get around this restriction, note that
\begin{equation}
e_i^T(A^TD^{(k)}A)e_i = ((A^T)\,.^{2}\,\text{diag}(D^{(k)}))_i \quad \text{and} \quad  e_i^T(A^TA)e_i =((A^T)\,.^2\,\mathbf{1})_i,
\label{eq:AT_squared}
\end{equation}
where $\mathbf{1}$ a vector of all ones, and 
we use MATLAB notation $.^2$ to mean component-wise squaring.
In some cases it may be relatively easy to compute both the entries of $(A^T).^2$ and the vector
$(A^T).^2\mathbf{1}$; this is the case for image deblurring, and is
discussed in more detail in Section~\ref{sec:num_exp}.

Using \eqref{eq:pullout_prec}, we define the preconditioner for the linear system \eqref{eq:lin_system} as
\begin{equation}
M \equiv\hat{D}^{(k)}  \left( A^TA + \hat{\lambda} L^TL \right) \hat{D}^{(k)},\label{eq:prec}
\end{equation}
with 
\[\hat{\lambda} \equiv \lambda/\text{mean}\left(\text{diag}(\hat{D}^{(k)})\right)^2 .\]
More details on the computational costs involved in constructing and applying the preconditioner
in the case of image deblurring are provided in Section~\ref{sec:num_exp}.

\section{Numerical tests}\label{sec:num_exp}

The Poisson--Gaussian model arises naturally in image applications, so in this section
we present numerical examples from image deblurring.  Specifically, we consider the 
inverse problem \eqref{eq:inverse_problem} with data model \eqref{eq:noise}, where
vector $b$ is an observed image that is corrupted by blur and noise, matrix $A$ models
the blurring operation, vector $x_\text{true}$ is the true image, and $\eta$ is noise.
Although an image is naturally represented as an array of pixel values, when we 
refer to `vector' representations, we assume the pixel values have been reordered
as vectors.  For example, if we have a $p \times p$ image of pixel values, these can
be stored in a vector of length $n = p^2$ by, for example, lexicographical ordering
of the pixel values.

In many practical image deblurring applications, the blurring is spatially invariant,
and $A$ is structured matrix defined by a {\em point spread function} (PSF).
In this situation, image deblurring can also be referred to as image decovolution,
because the operation $Ax_\text{true}$ is the convolution of $x_\text{true}$ and the PSF.
Although the PSF may be given as an actual function, the more common situation is
to compute estimates of it by imaging point source objects.  Thus, the PSF can be 
represented as an image; we typically display the PSF as a mesh plot, which makes
it easier to visualize how a point in an image is spread to its neighbors because of
the blurring operation.
The precise structure of the matrix $A$ depends on the imposed boundary condition;
see \cite{Hansen2006Deblurring} for details.  In this section we impose periodic boundary
conditions, so that $A$ and $L$ are both diagonalizable by FFTs.

So far we have only described what we refer to as the {\em single-frame} situation, where
$b$ is a single observed image.  It is often the case, especially in astronomical imaging,
to have multiple observed images of the same object, but with each having
a different blurring matrix associated with it.  We refer to this as the
multi-frame image deblurring problem.  Here, $b$ represents all observed images, stacked
one on top of each other, and similarly $A$ is formed by stacking the various blurring matrices.

Before describing the test problems used in this section, we first summarize the computational
costs.
From the discussion around equation \eqref{eq:AT_squared}, to construct the preconditioner we need to 
be able to efficiently square all entries of the matrix $A^T$, or equivalently those of $A$; this can
easily be approximated by squaring the point-spread function component-wise before forming the operator, i.e.,
\begin{equation}
(A_\text{PSF}).^2 \approx A_{\text{PSF}.^2}\, .\label{eq:PSF_squaring}
\end{equation}
Using this approximation, in each Newton step we only need to perform one multiplication by a matrix,
one component-wise multiplication, and one component-wise square-root to obtain the entries of the diagonal matrix \eqref{eq:out_diag_def}. 
With the assumption that $A$ and $L$ are both diagonalizable by FFTs, 
efficient multiplication by the Hessian \eqref{eq:lin_system} involves two 
two-dimensional forward and inverse FFTs, which we refer to as \texttt{fft2} and \texttt{ifft2}, respectively. 
Solving systems with matrix \eqref{eq:prec} involves only  one \texttt{fft2} and one \texttt{ifft2}. 
In addition to the \texttt{fft2} requirements, multiplication by the Hessian \eqref{eq:lin_system} involves 4 pixel-wise multiplications and 1 addition. Solving systems with the preconditioner \eqref{eq:prec} involves 3 pixel-wise multiplications (component-wise reciprocals are assumed to be computed only once at the beginning). The total counts for each operation are shown in Table \ref{tab:op_counts}.
\begin{table}
\caption{Operation counts for single-frame case.}\label{tab:op_counts}
\centering
{\small
\begin{tabular}{llcccc}
\toprule
& operation& \texttt{fft2} & \texttt{ifft2} & mults & adds \\ 
\midrule
Hessian \eqref{eq:lin_system} &multiply & 2 & 2 & 4 & 1 \\ 
preconditioner \eqref{eq:prec} &solve & 1 & 1 & 3 & 0 \\ 
\bottomrule
\end{tabular}} 
\end{table}

\begin{figure}[!th]
\centering
\includegraphics[width=.2\textwidth]{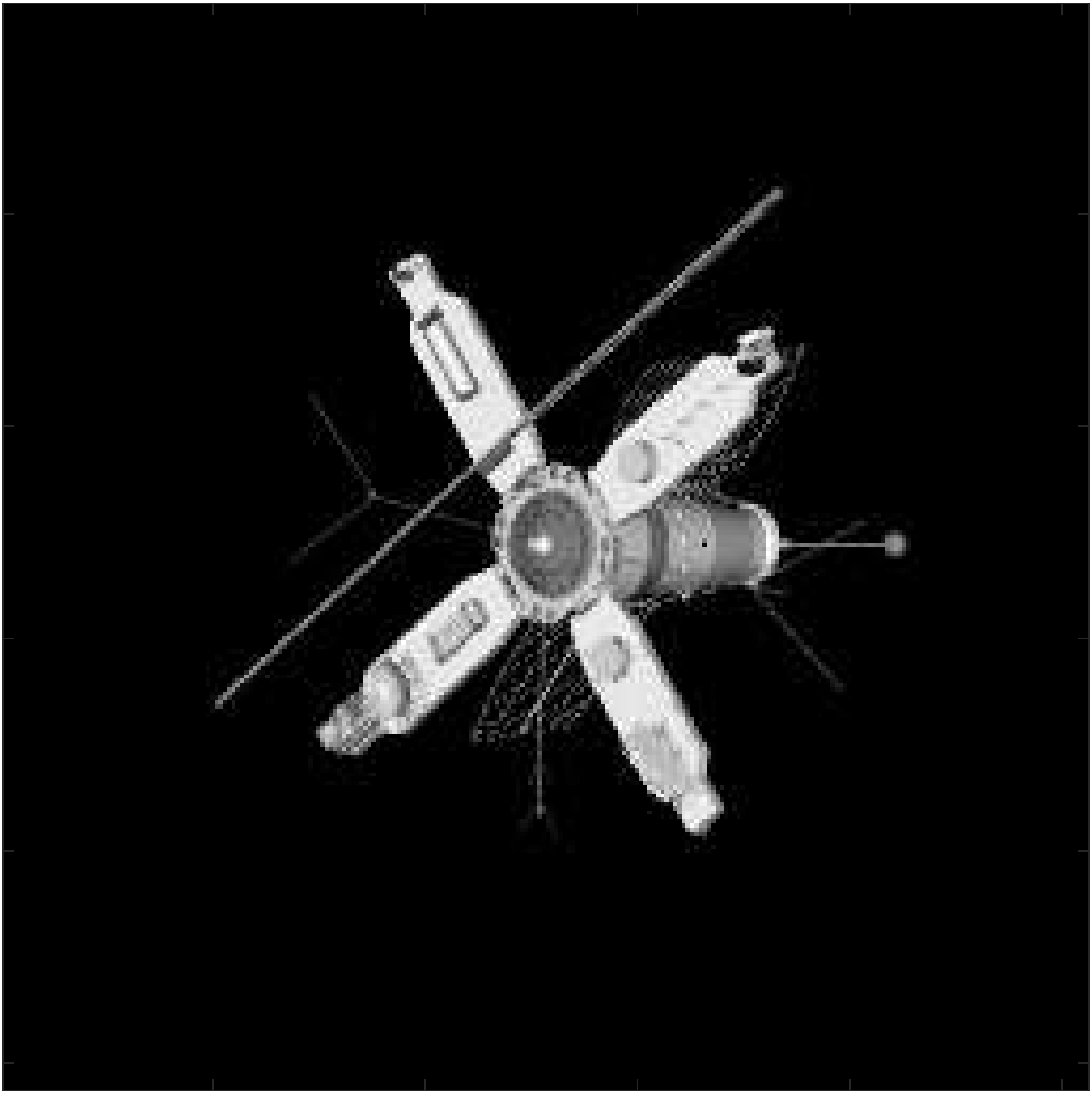}
\hspace*{.5cm}
\includegraphics[width=.2\textwidth]{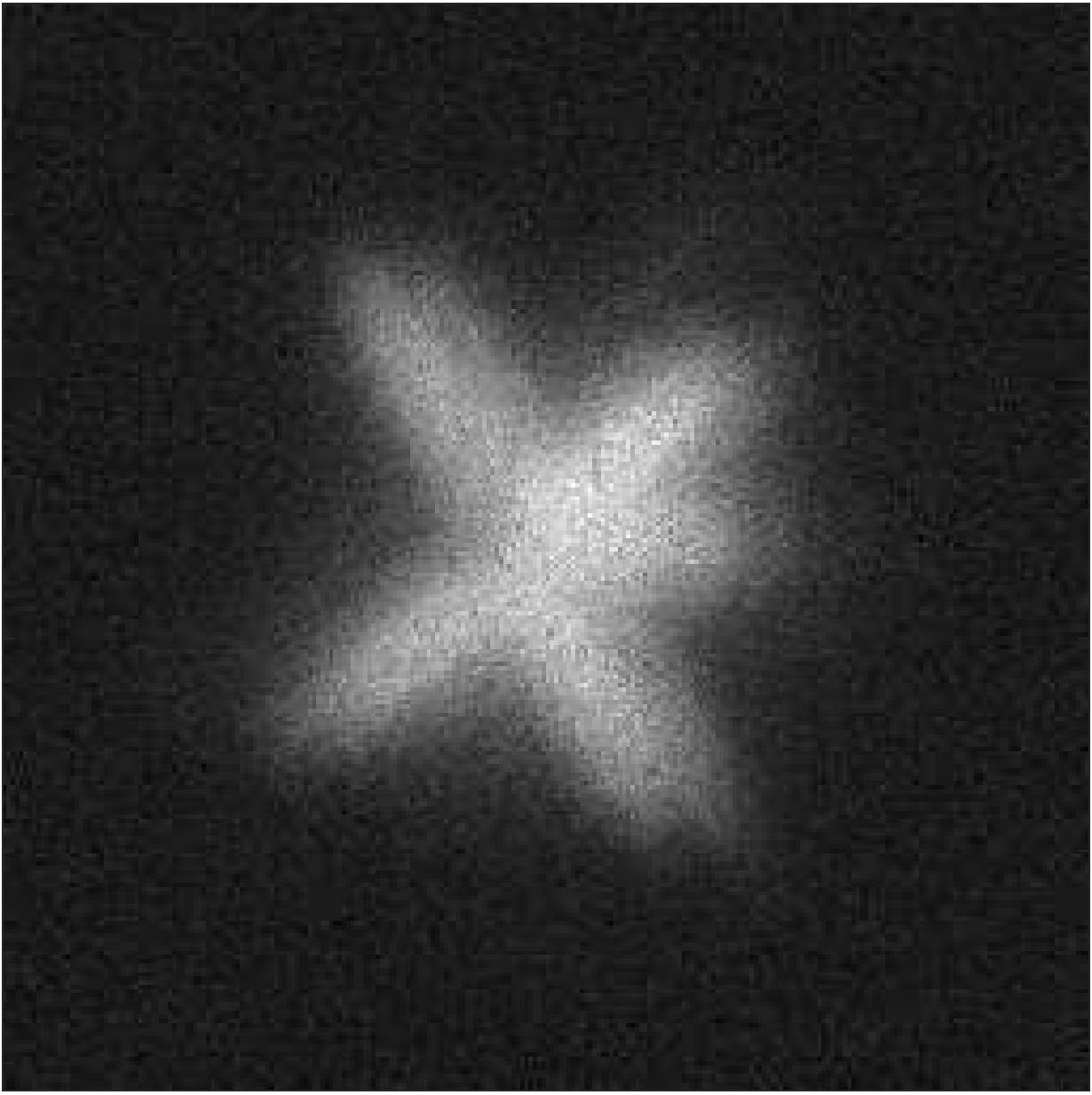}
\includegraphics[width=.2\textwidth]{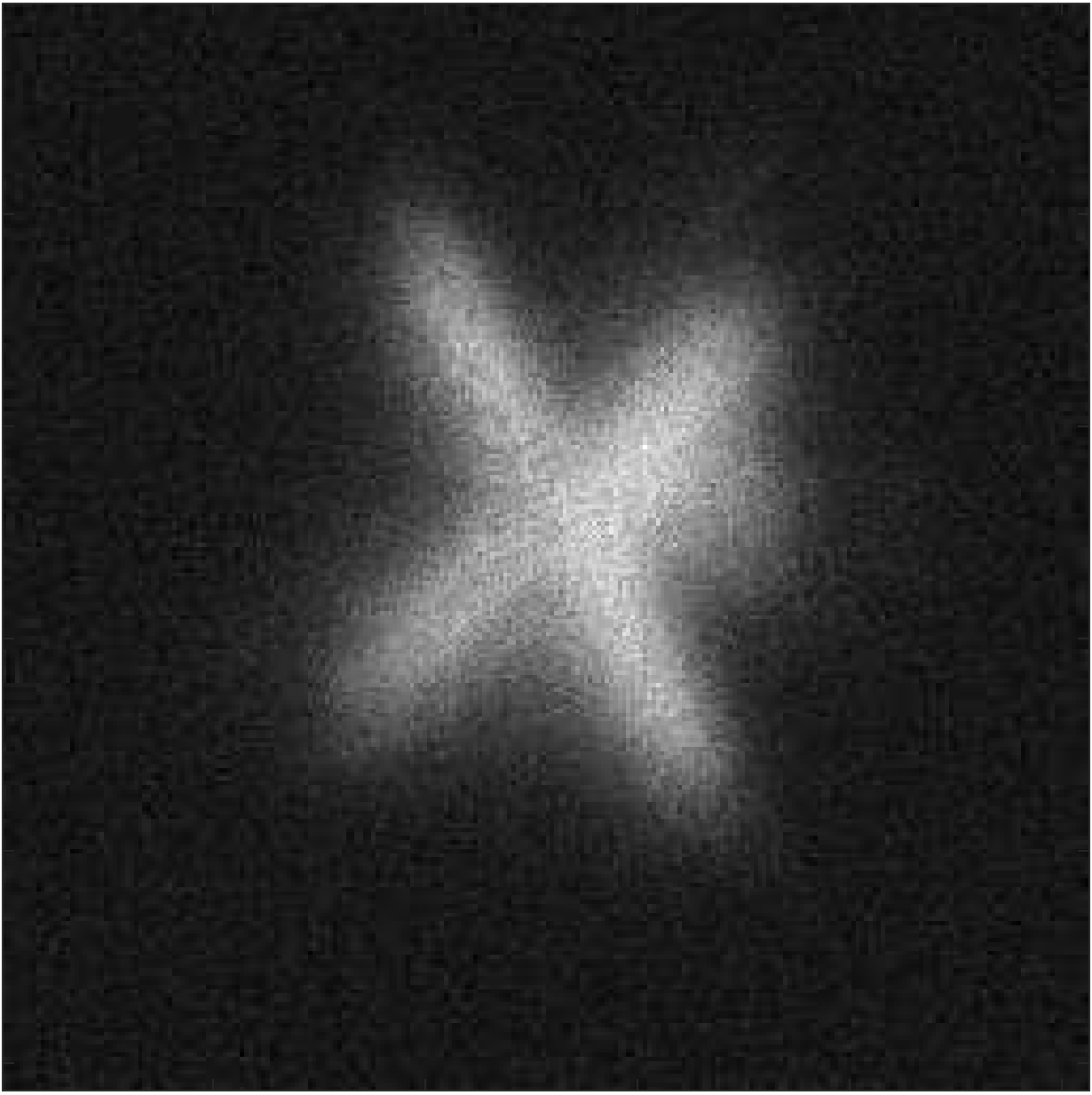}
\includegraphics[width=.2\textwidth]{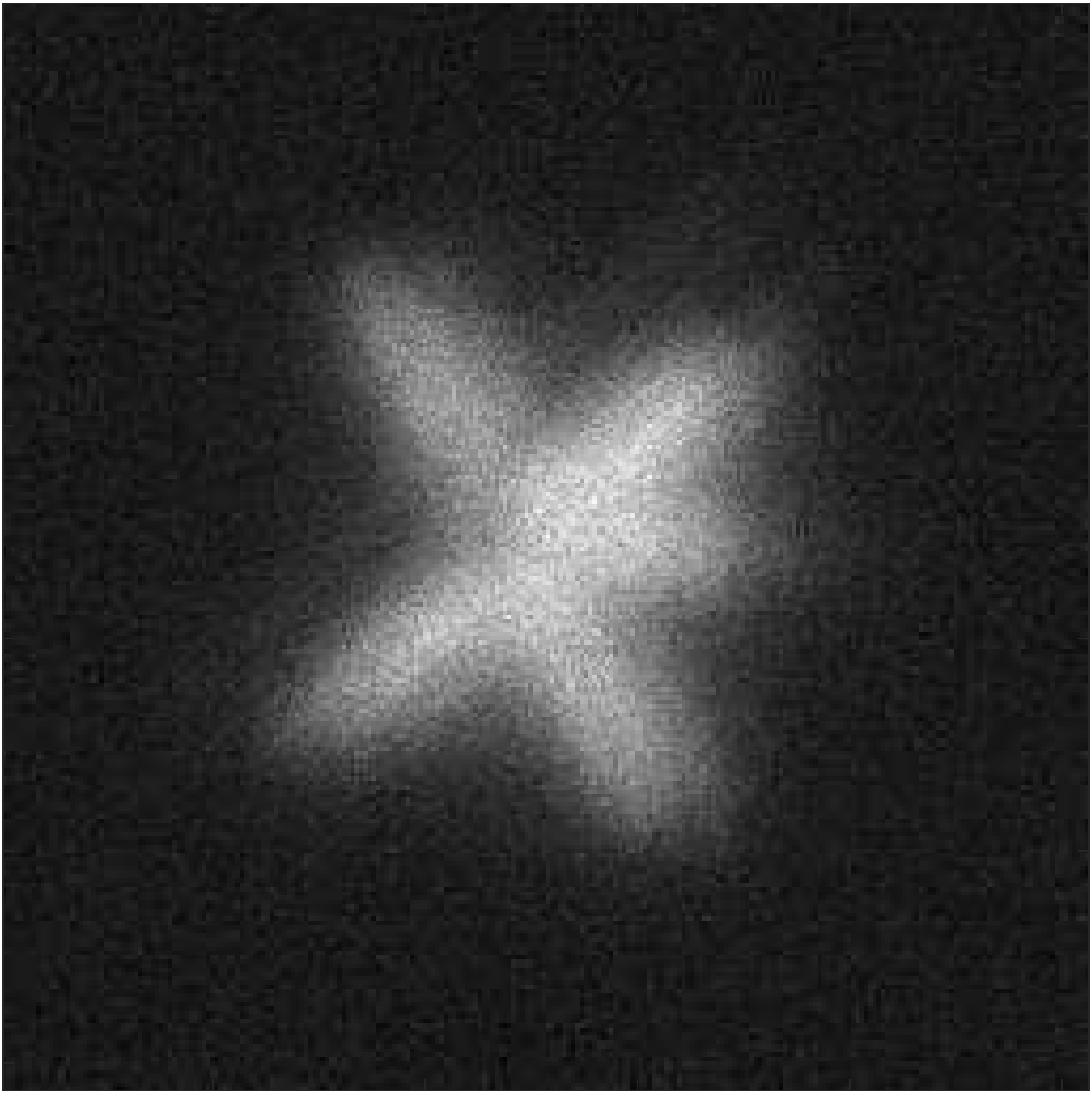}
\caption{Test problem Satellite: true image (left) together with three blurred noisy images (right). 
}\label{fig:satellite}
\end{figure}

\begin{figure}[!th]
\centering
\includegraphics[width=.2\textwidth]{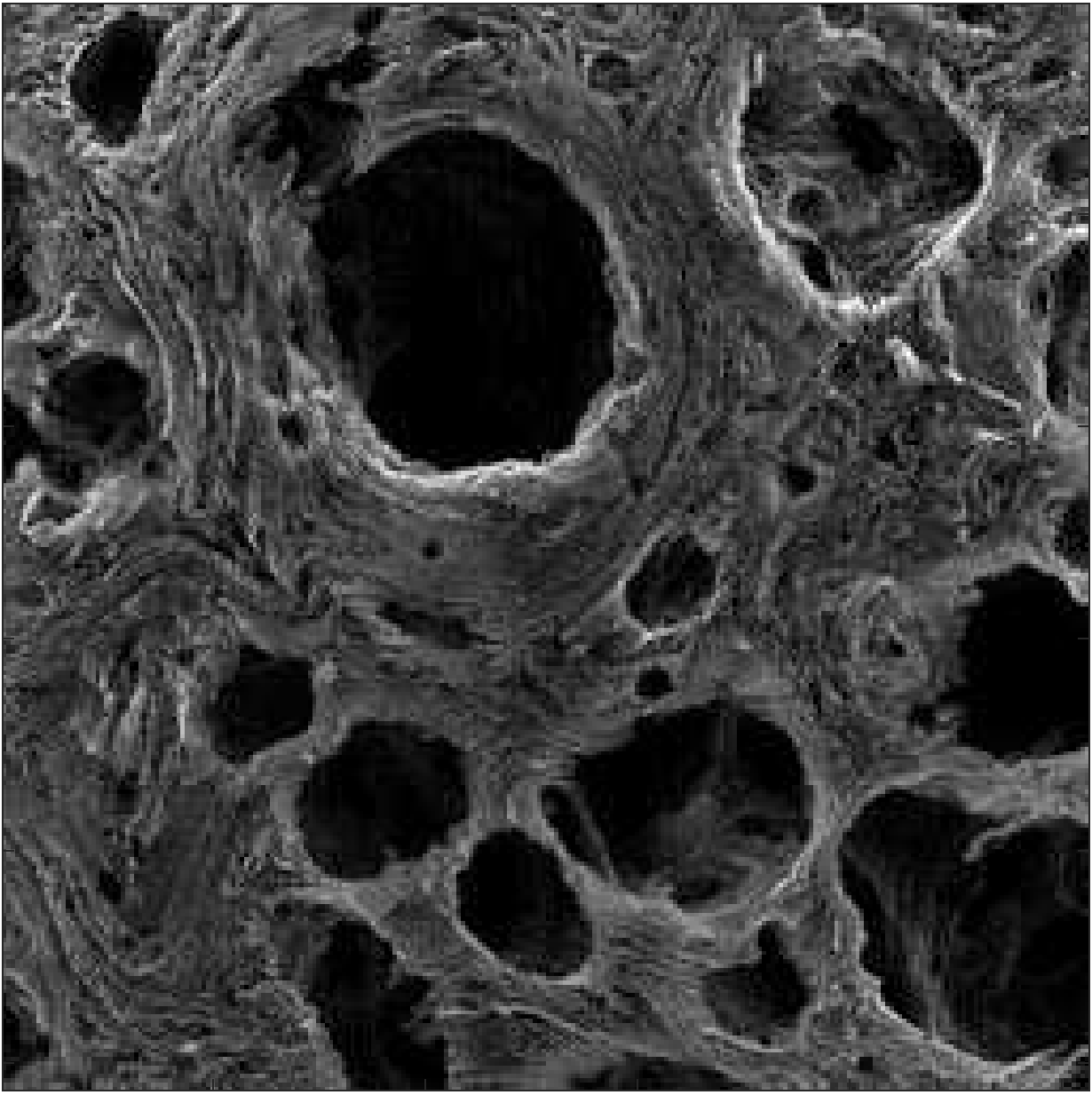}
\hspace*{.5cm}
\includegraphics[width=.2\textwidth]{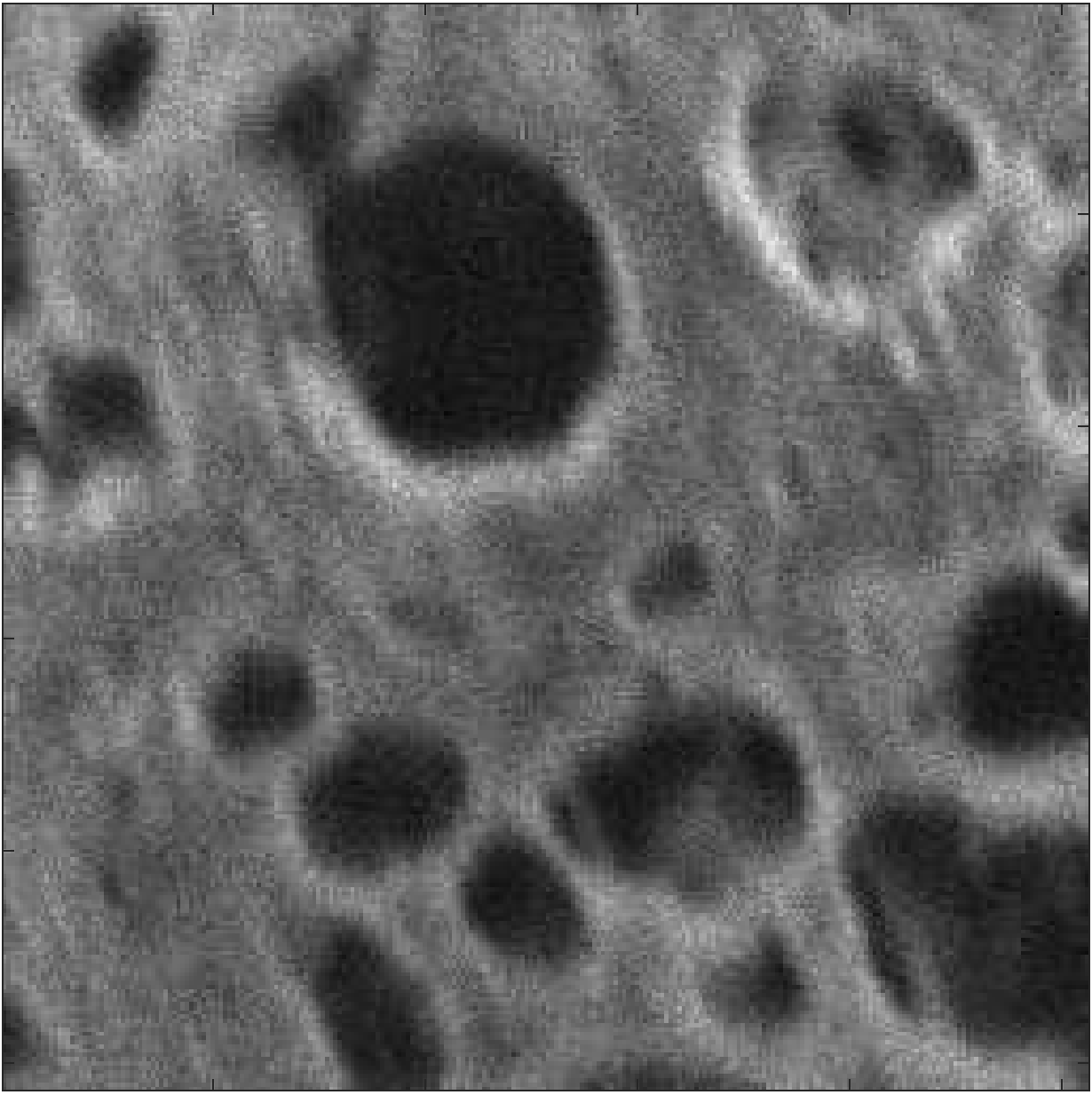}
\includegraphics[width=.2\textwidth]{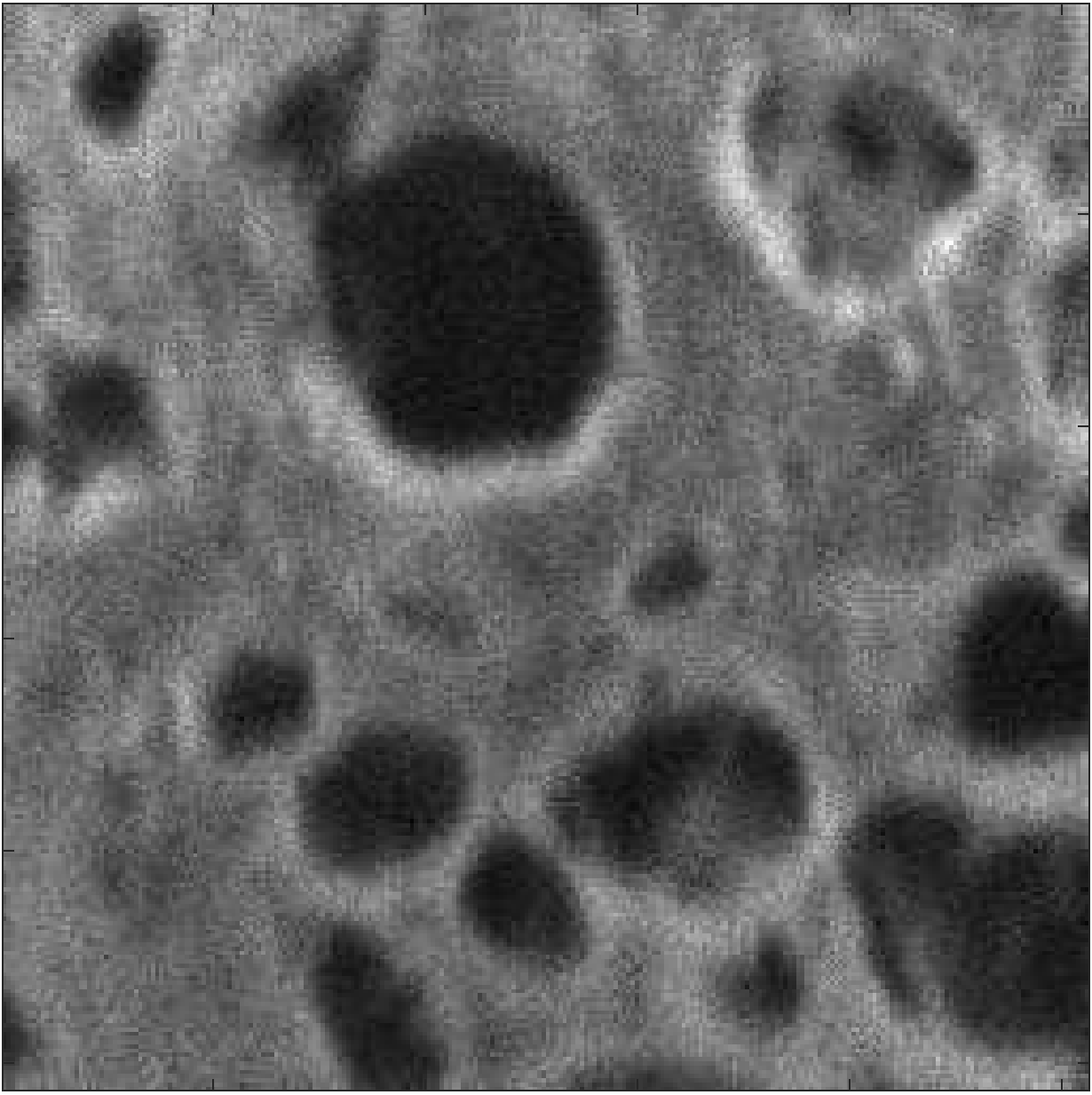}
\includegraphics[width=.2\textwidth]{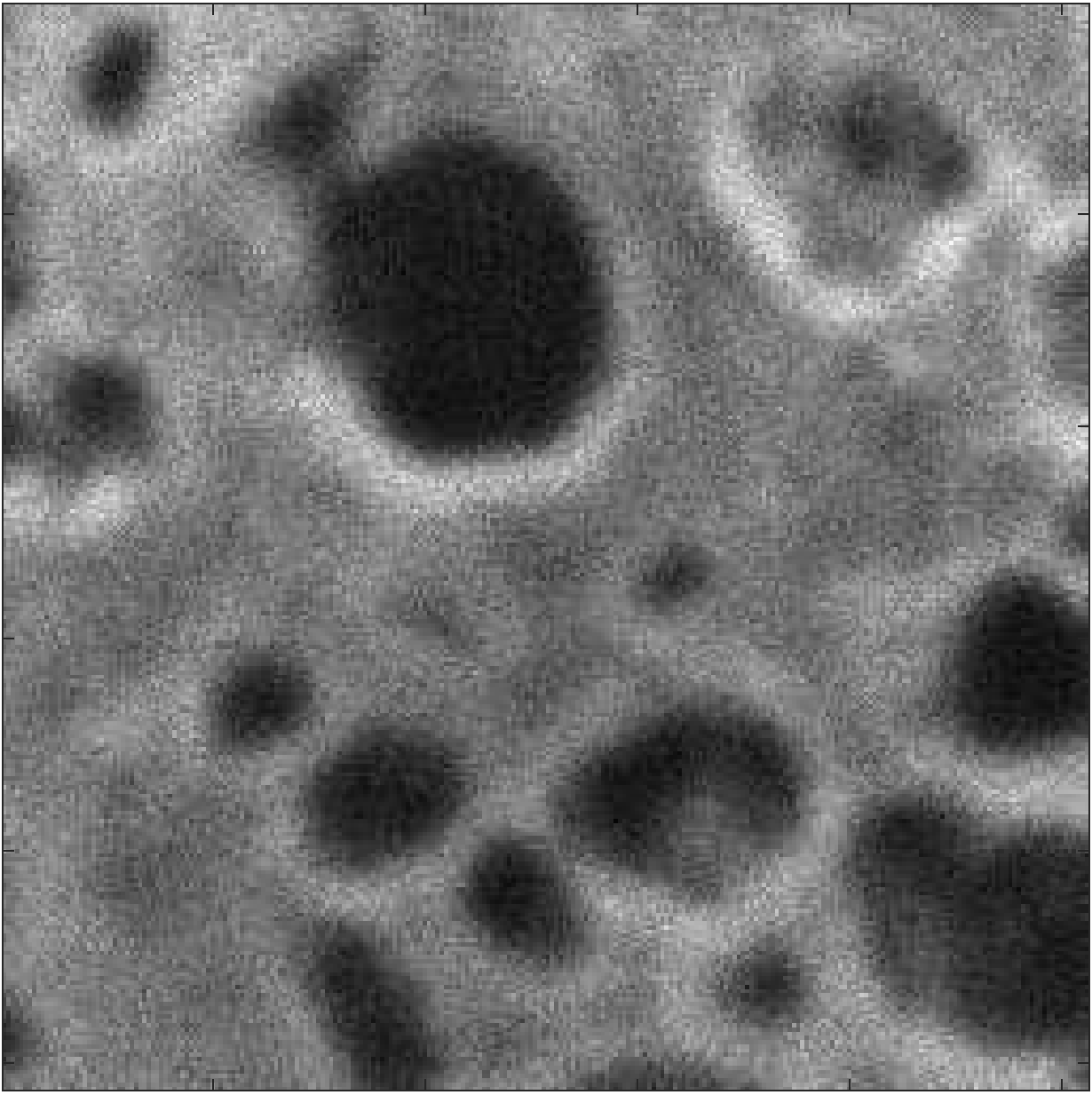}
\caption{Test problem Carbon ash: true image (left) together with three blurred noisy images (right). 
}\label{fig:carbon_ash}
\end{figure}

\begin{figure}[!th]
\centering
\begin{subfigure}[b]{.45\textwidth}
\includegraphics[width=.95\textwidth]{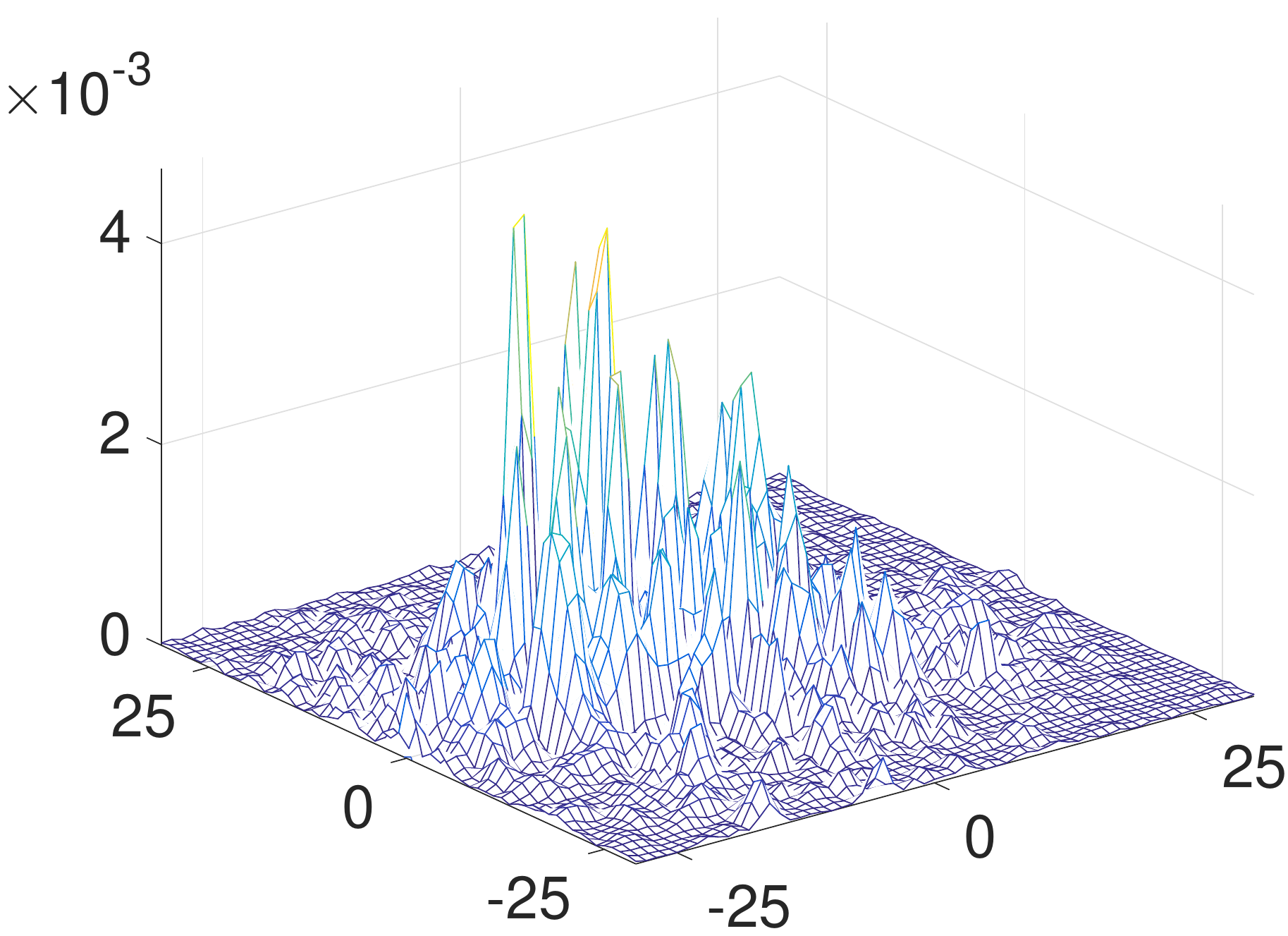}
\caption{Satellite}\label{fig:PSF_mesh_1}
\end{subfigure}
\begin{subfigure}[b]{.45\textwidth}
\includegraphics[width=.95\textwidth]{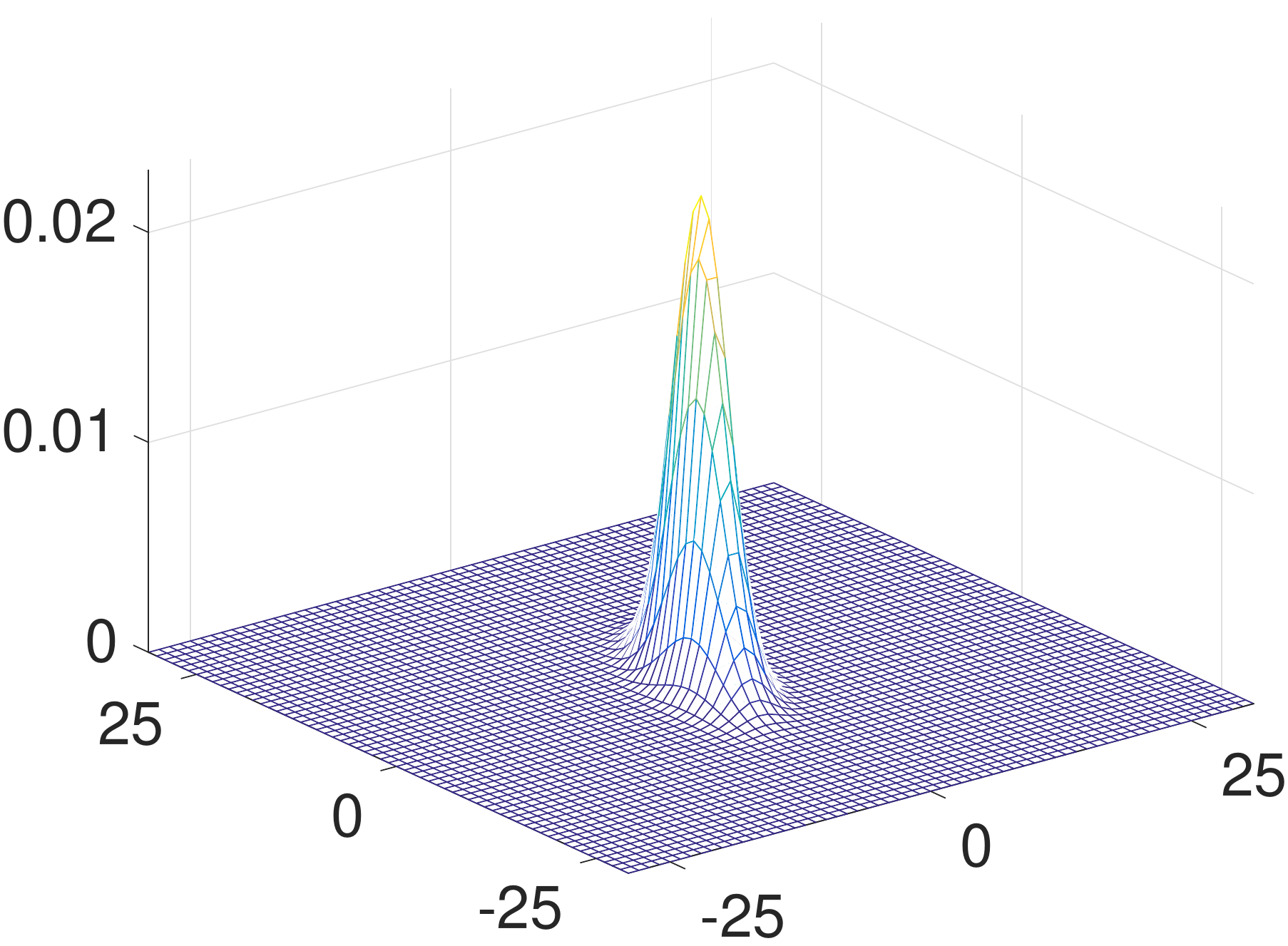}
\caption{Carbon ash}\label{fig:PSF_mesh_2}
\end{subfigure}
\caption{Point-spread functions for the first frame of each test problem.
}\label{fig:setting}
\end{figure}

The robustness and the efficiency of the proposed method is demonstrated on two test problems adopted from \cite{Nagy2004Iterative}:
\paragraph{Satellite}
An atmospheric seeing problem with spatially invariant atmospheric blur (moderate seeing conditions with the Fried parameter 30). We also consider a multi-frame case, where the same object is blurred by three different PSFs. 
These PSFs are generated by transposing and flipping the first PSF. 
The setting is shown in Figures \ref{fig:satellite} and \ref{fig:PSF_mesh_1}.
\paragraph{Carbon ash}
An image deblurring problem with spatially invariant non-separable Gaussian blur, where the PSF has the
functional definition
\begin{equation}
\label{eq:GaussianPSF}
  \mbox{PSF}(s,t) = \frac{1}{2\pi \sqrt{\gamma}} \exp\left\{
  -\frac{1}{2} \left[ \begin{array}{cc} s & t \end{array} \right] C^{-1}
   \left[ \begin{array}{c} s \\ t \end{array} \right] \right\}\, ,
\end{equation}
where 
$$
  C = \left[ \begin{array}{cc} \gamma_1^2 & \tau^2 \\[3pt] \tau^2 & \gamma_2^2 \end{array} \right]\,, \quad
  \mbox{and} \quad \gamma_1^2 \gamma_2^2 - \tau^4 > 0\,.
$$
The shape of the Gaussian PSF depends on the parameters $\gamma_1, \gamma_2$ and $\tau$;
we use
$\gamma_1 = 4$, $\gamma_2 = 2$; $\tau = 2$. We also consider a multi-frame case, where the same object is blurred by three different PSFs. The other two PSFs are Gaussian blurs with parameters $\gamma_1 = 4$, $\gamma_2 = 2$, $\tau = 0$, and $\gamma_1 = 4$, $\gamma_2 = 2$, $\tau = 0$. The setting is shown in Figures \ref{fig:carbon_ash} and \ref{fig:PSF_mesh_2}.

As previously mentioned, in the multi-frame case, the vector $b$ in (\ref{eq:inverse_problem}) is concatenation of the vectorized blurred noisy images, the matrix $A$ is concatenation of the blurring operators, i.e., $A\in\mathbb{R}^{3n\times n}$. For the test problems all true images are $256 \times 256$ arrays of pixels (with intensities scaled to $[0,255]$), and thus
$n = 65536$.

Computation was performed in MATLAB R2015b. Noise is generated artificially using MATLAB functions \texttt{poissrnd} and \texttt{randn}. Unless specified otherwise, the standard deviation $\sigma$ is set to $5$. We use the discretized Laplacian, see \cite[p. 95]{Hansen2006Deblurring}, as the regularization matrix $L$. The Projected Newton method (Algorithm \ref{alg:projGNCG}) is terminated when the relative size of the projected gradient
\begin{equation}
\mathcal{P}(\text{grad}_{J_{\lambda}}(x^{(k)})), \quad \text{where} \quad \mathcal{P}(v) \equiv v.*(1 - \text{Active}) + \text{Active}.*(v < 0),
\label{eq:proj_grad_def}
\end{equation}
reaches the tolerance $10^{-4}$ or after 40 iterations. We use MATLAB notation $.^*$ to mean component-wise multiplication. Projected PCG (Algorithm \ref{alg:projPCG}) is terminated when the relative size of the projected residual (denoted in Algorithm \ref{alg:projPCG} by $r_i$) reaches $10^{-1}$, or the number of iterations reaches $100$. 
We use the preconditioner given in \eqref{eq:prec} as the default preconditioner. Given a search direction $s_k$, we apply a projected backtracking linesearch, with the initial step length equal to 1, which we terminate when 
\[J_\lambda(x^{(k+1)}) < J_\lambda(x^{(k)}).\]  

\subsection{Robustness with respect to various types of outliers}

In this section, we consider several types of outliers, whose choice was motivated by \cite{Calef2013Iteratively}, and demonstrate the robustness of the proposed method. Note that the difference between \cite{Calef2013Iteratively} and the proposed approach lies, among others, in the fact that while in \cite{Calef2013Iteratively}, the approximation of the solution is computed in order to update the outer (robust) weights associated with the components of residual. Here, the weights are represented by the loss function $\rho$ and are updated implicitly in each Newton step and therefore our approach does not involve any outer iteration.

\subsubsection*{Random corruptions}
\begin{figure}[!ht]
\centering
\begin{subfigure}[b]{\textwidth}
\includegraphics[width=.31\textwidth]{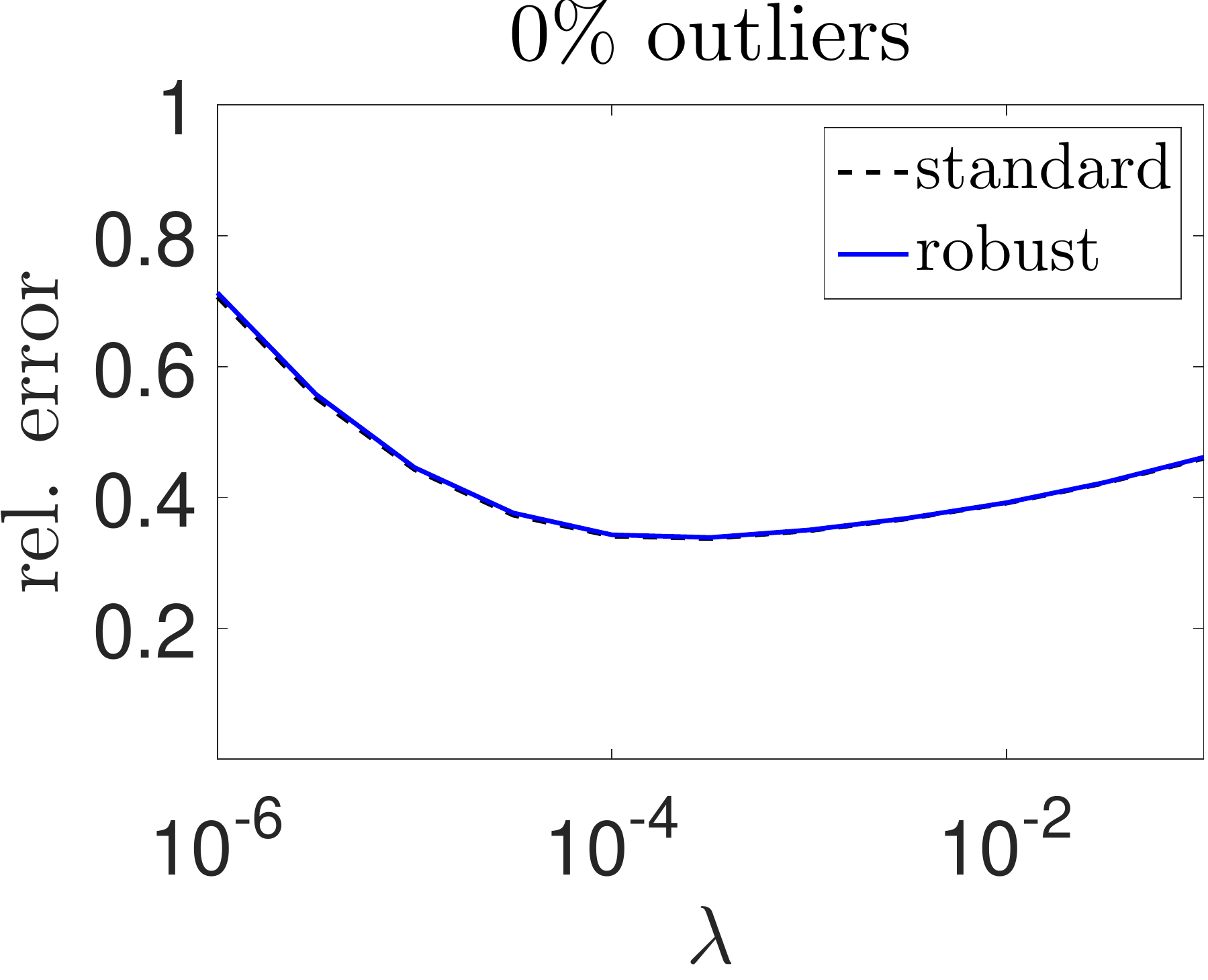}
\hspace*{.2cm}
\includegraphics[width=.31\textwidth]{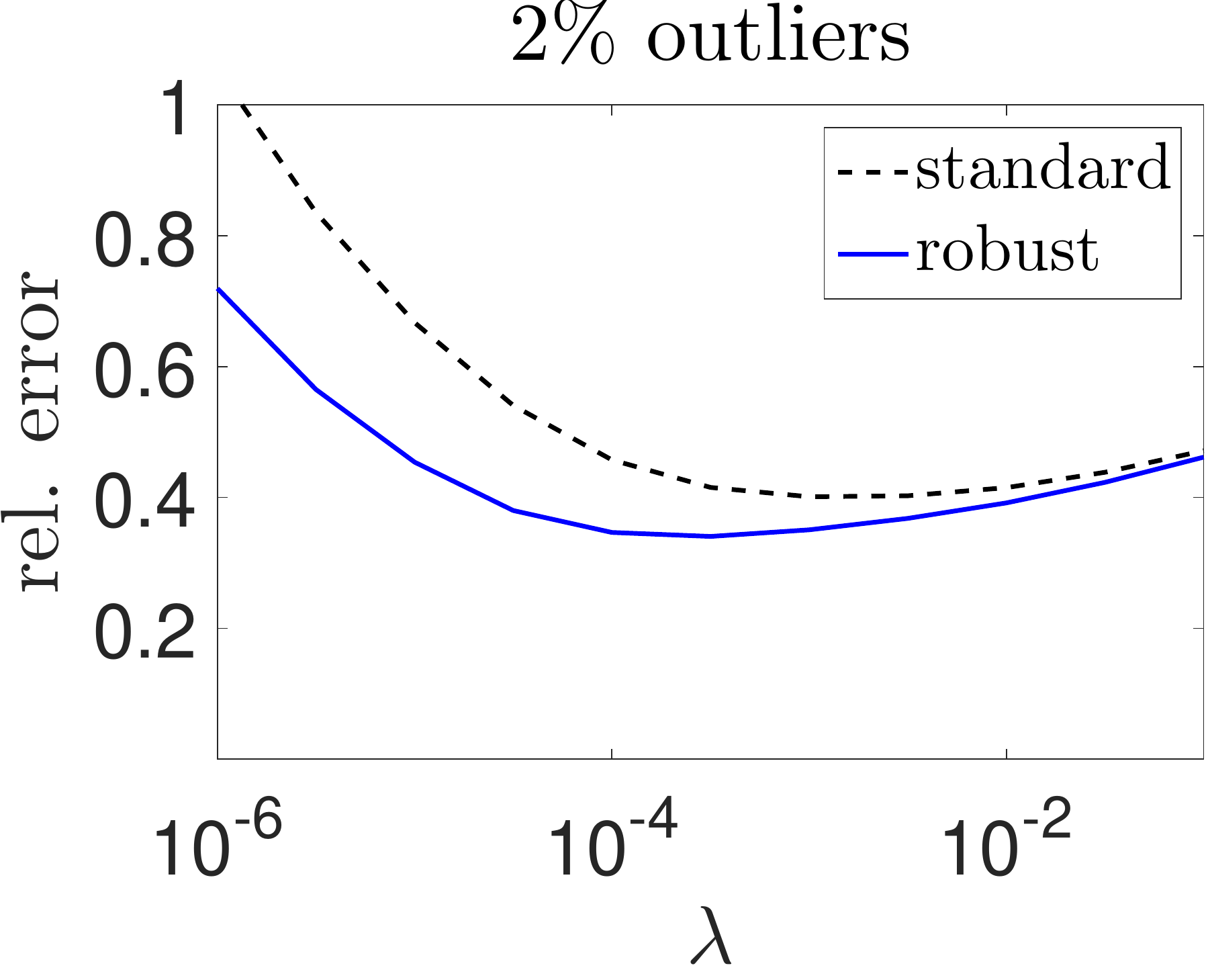}
\hspace*{.2cm}
\includegraphics[width=.31\textwidth]{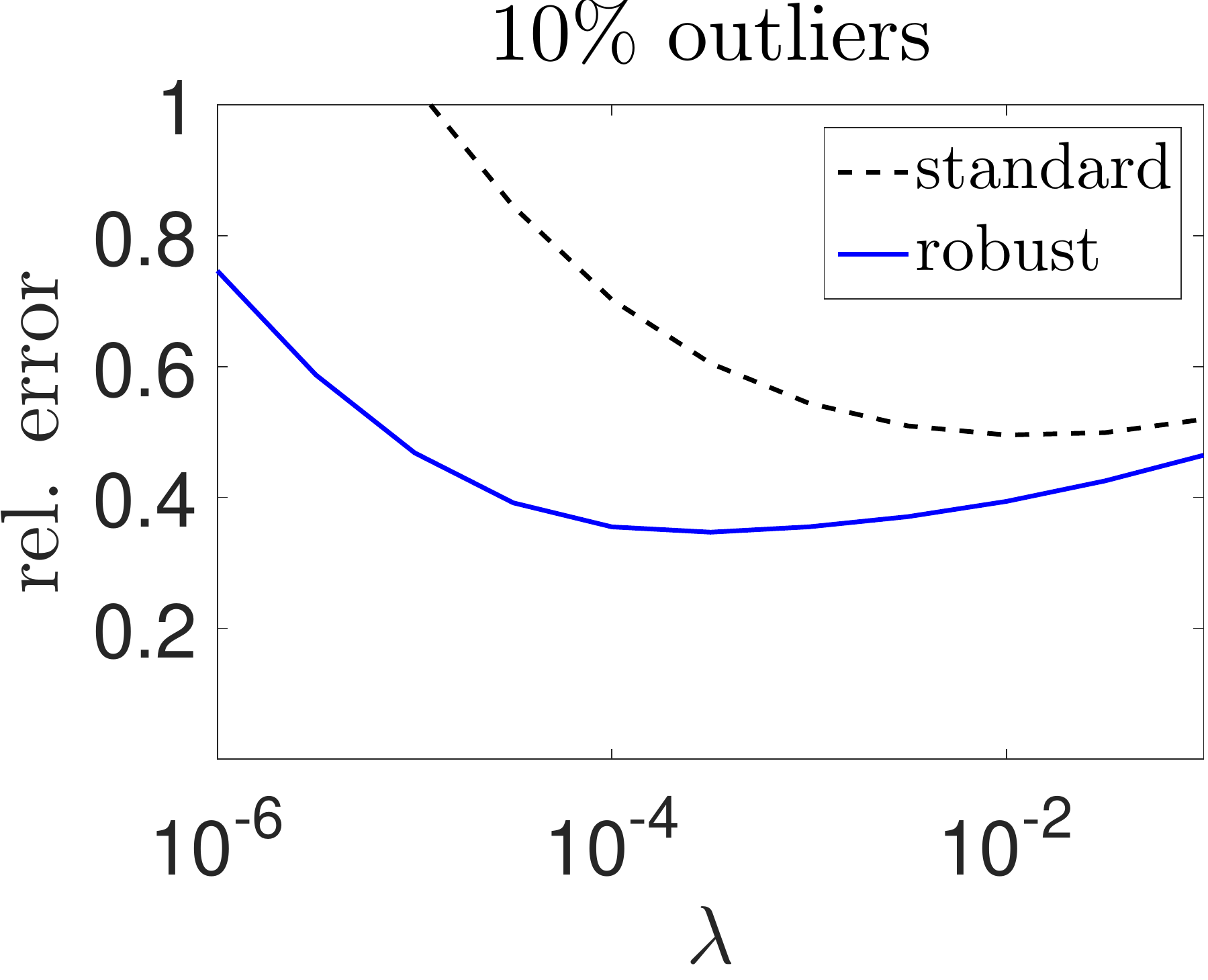}
\caption{Satellite single-frame}\label{fig:sliding_curves_1}
\end{subfigure}

\vspace*{.3cm}

\begin{subfigure}[b]{\textwidth}
\includegraphics[width=.31\textwidth]{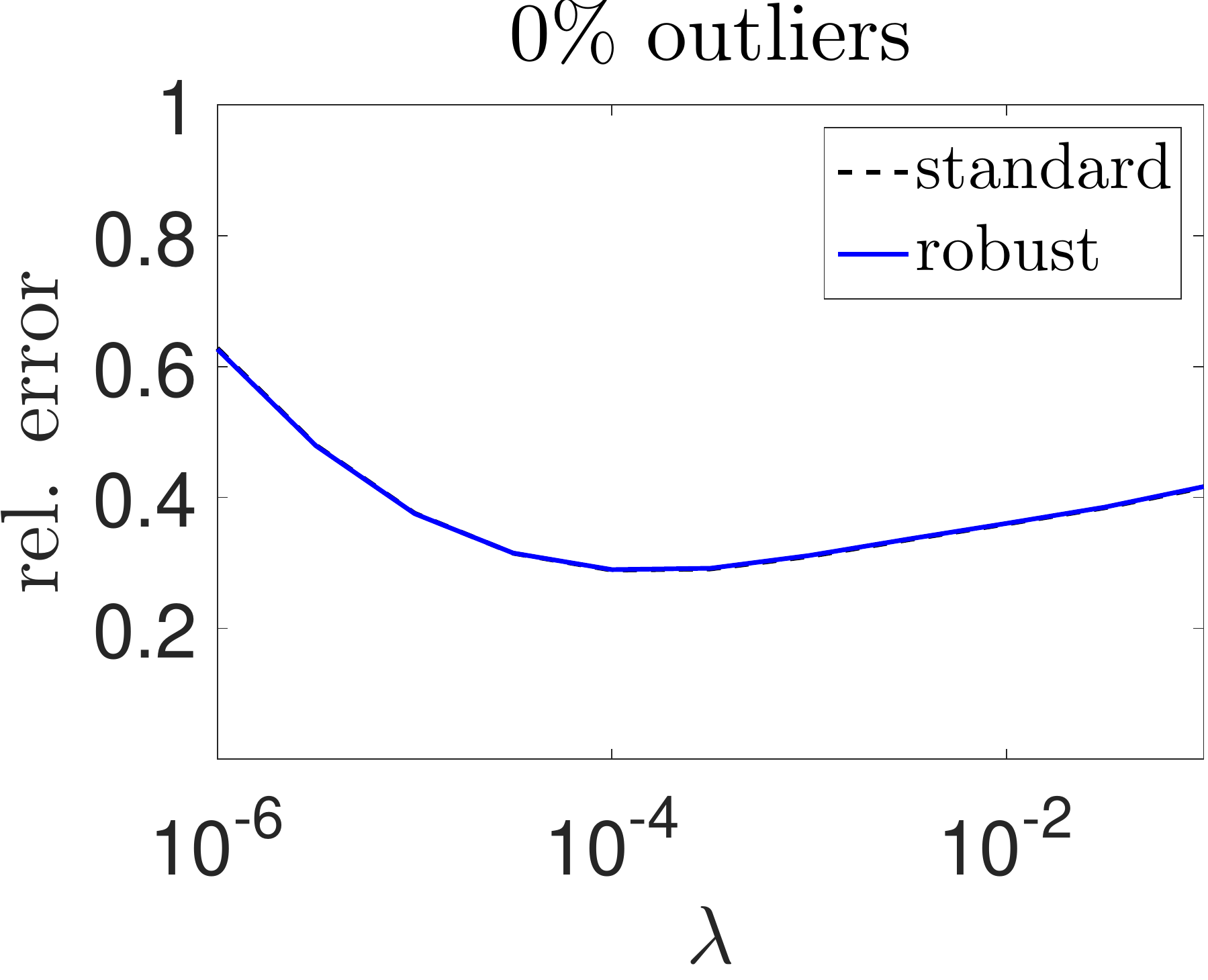}
\hspace*{.2cm}
\includegraphics[width=.31\textwidth]{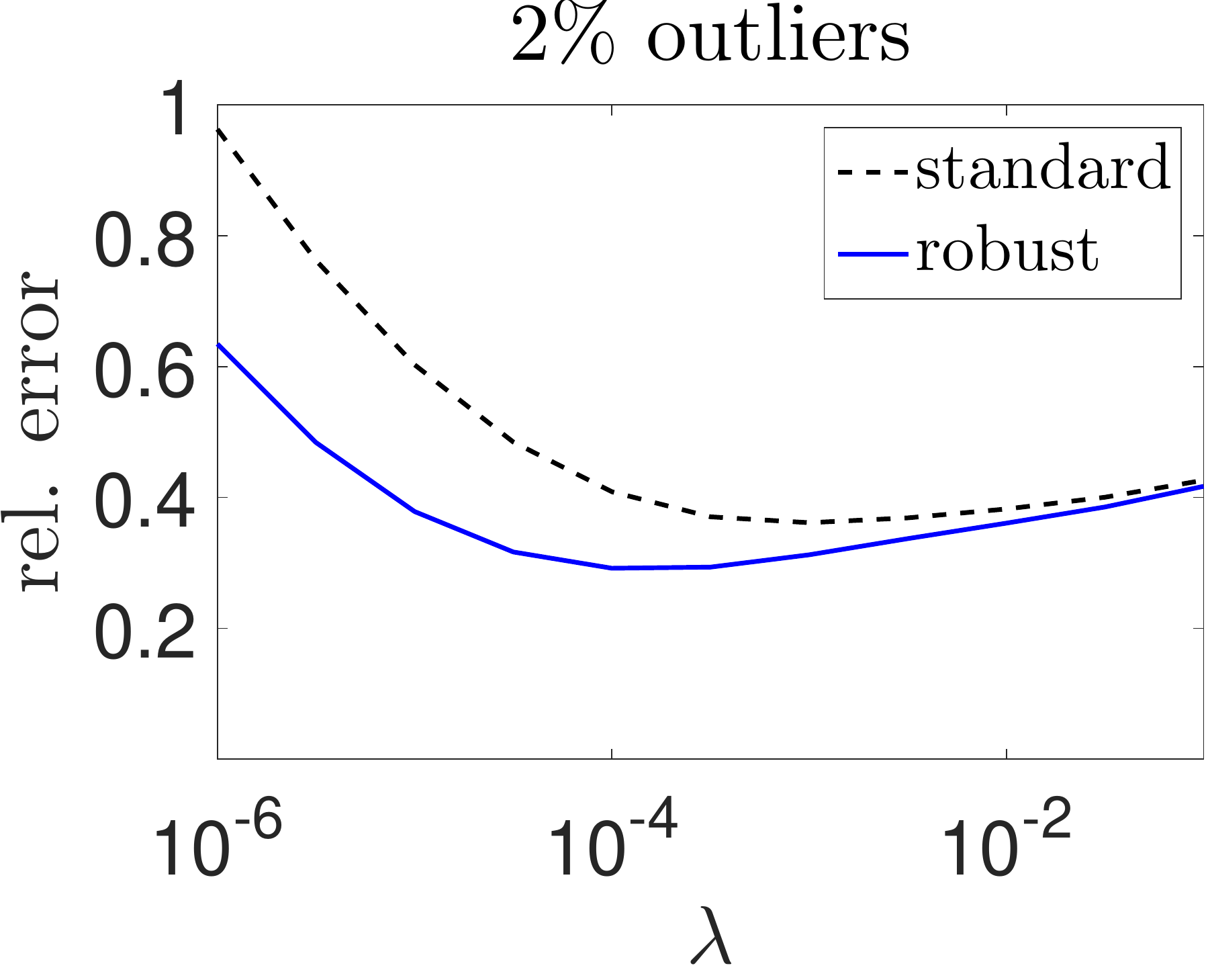}
\hspace*{.2cm}
\includegraphics[width=.31\textwidth]{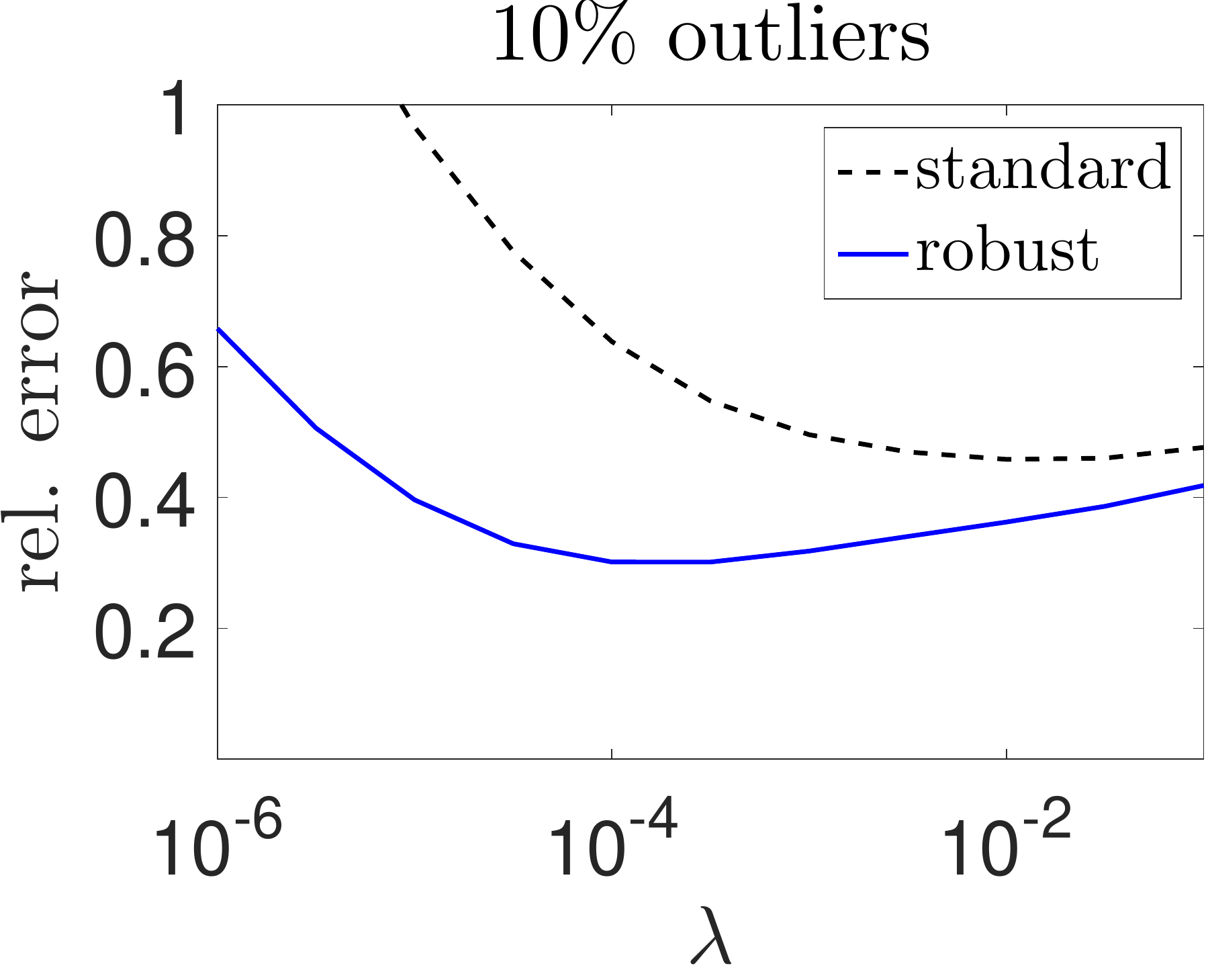}
\caption{Satellite multi-frame}\label{fig:sliding_curves_2}
\end{subfigure}

\vspace*{.3cm}

\begin{subfigure}[b]{\textwidth}
\includegraphics[width=.31\textwidth]{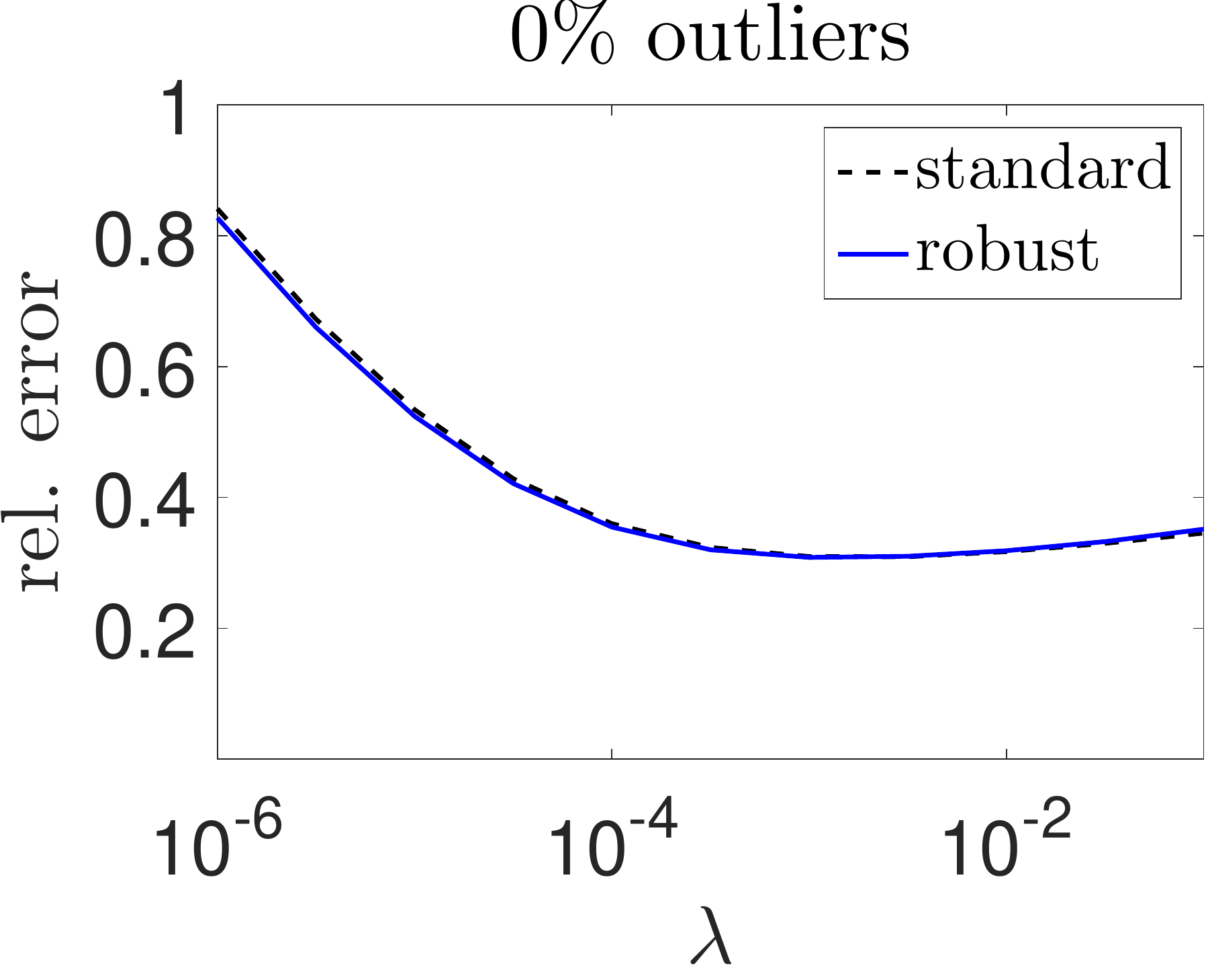}
\hspace*{.2cm}
\includegraphics[width=.31\textwidth]{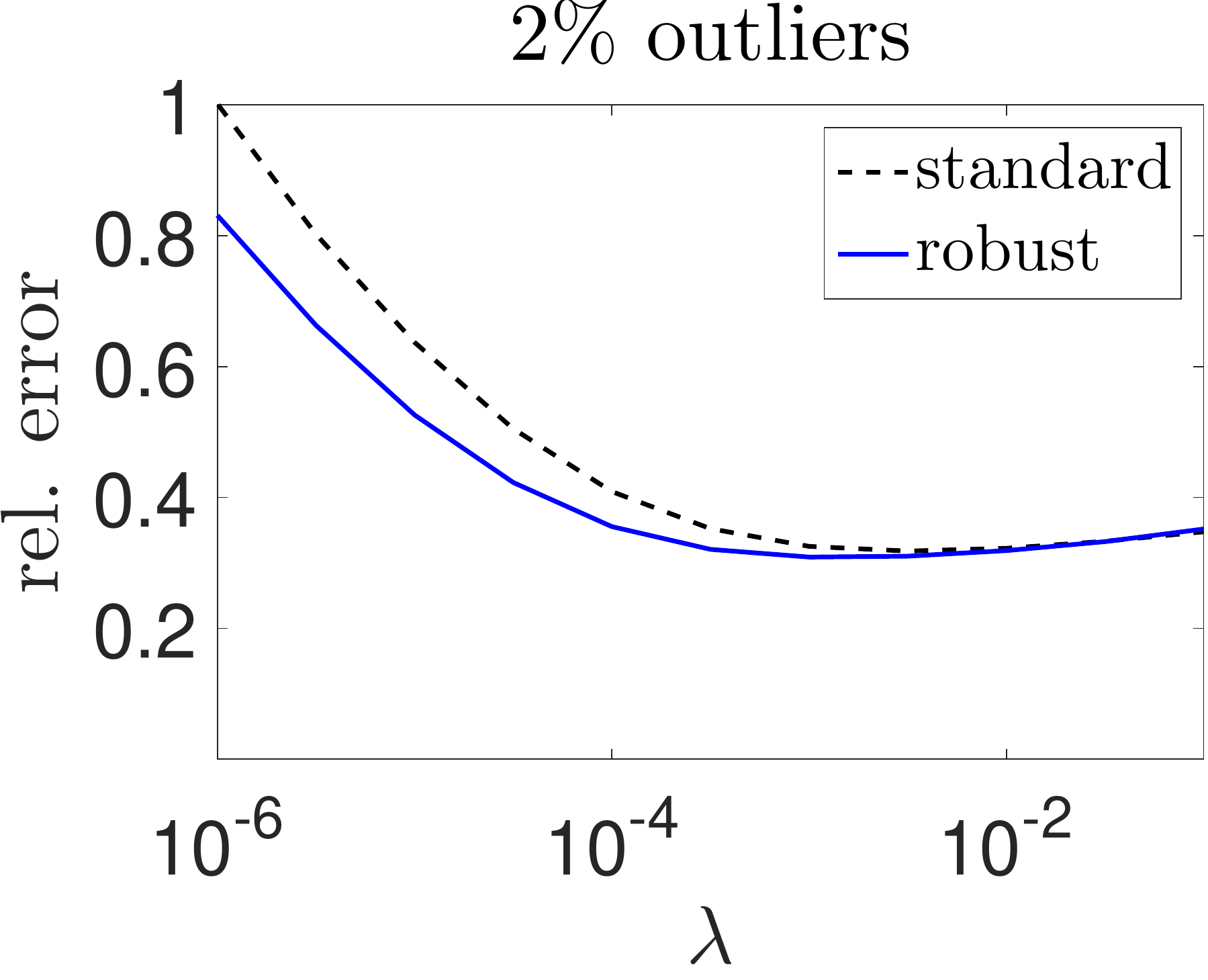}
\hspace*{.2cm}
\includegraphics[width=.31\textwidth]{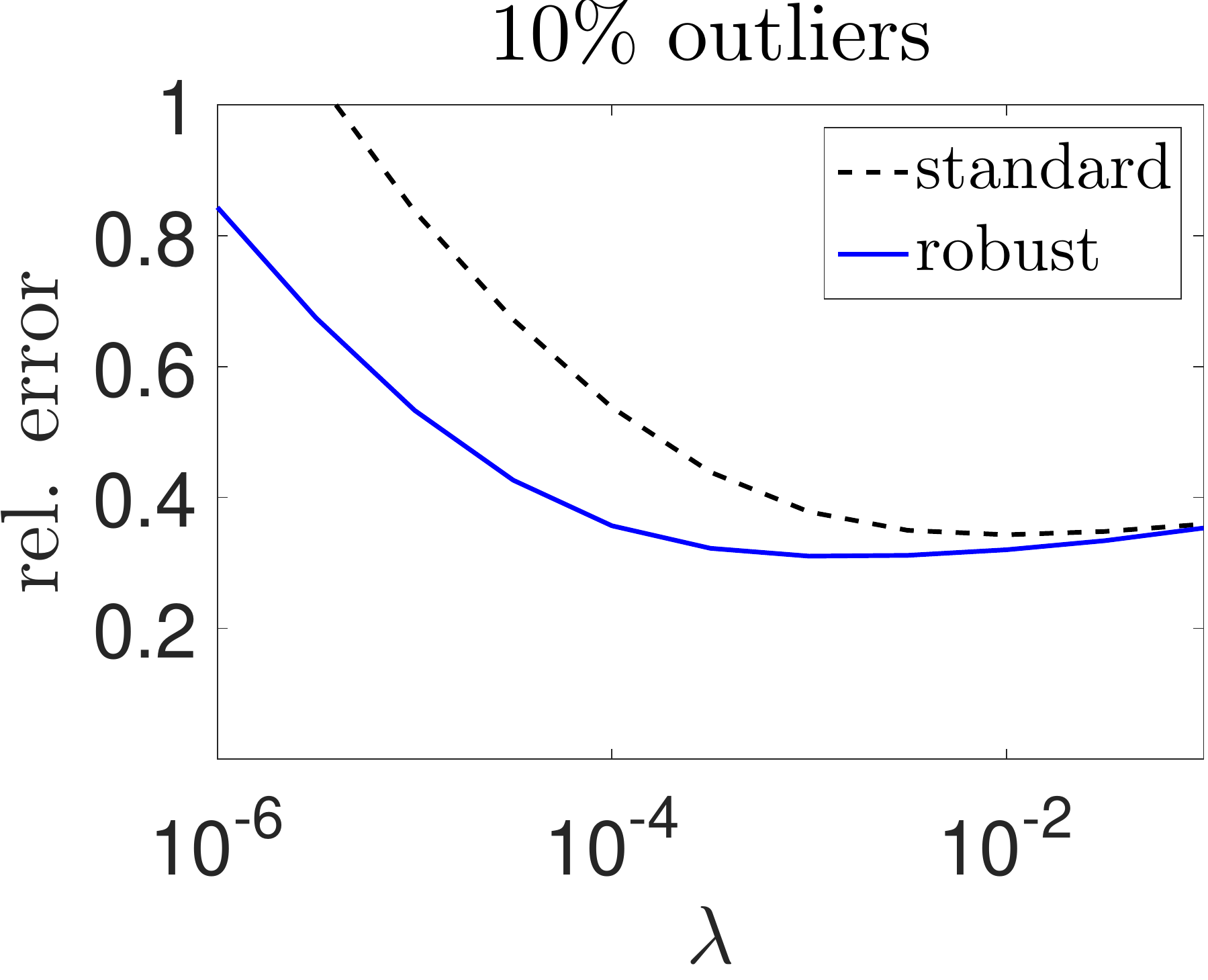}
\caption{Carbon ash single-frame}\label{fig:sliding_curves_3}
\end{subfigure}

\vspace*{.3cm}

\begin{subfigure}[b]{\textwidth}
\includegraphics[width=.31\textwidth]{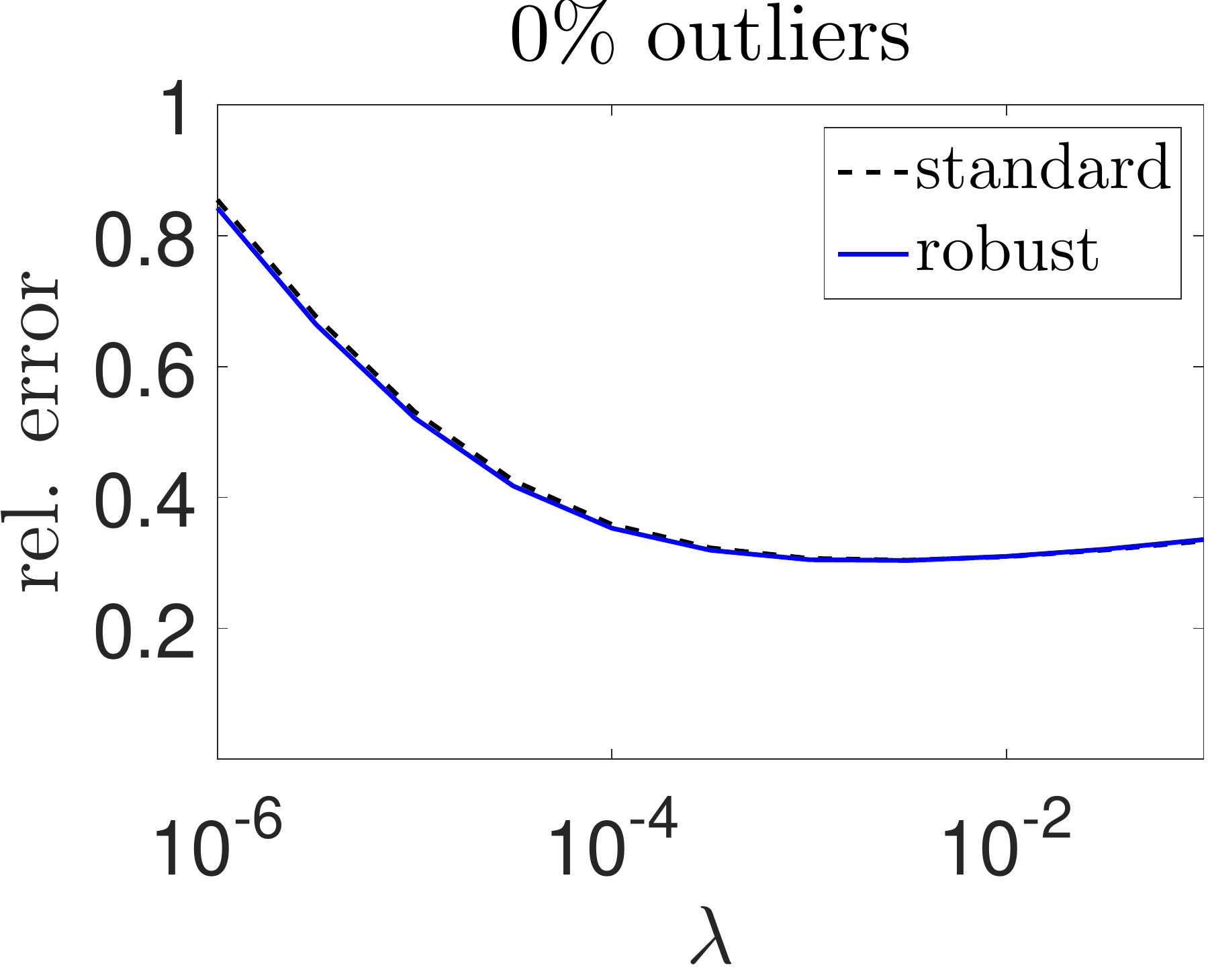}
\hspace*{.2cm}
\includegraphics[width=.31\textwidth]{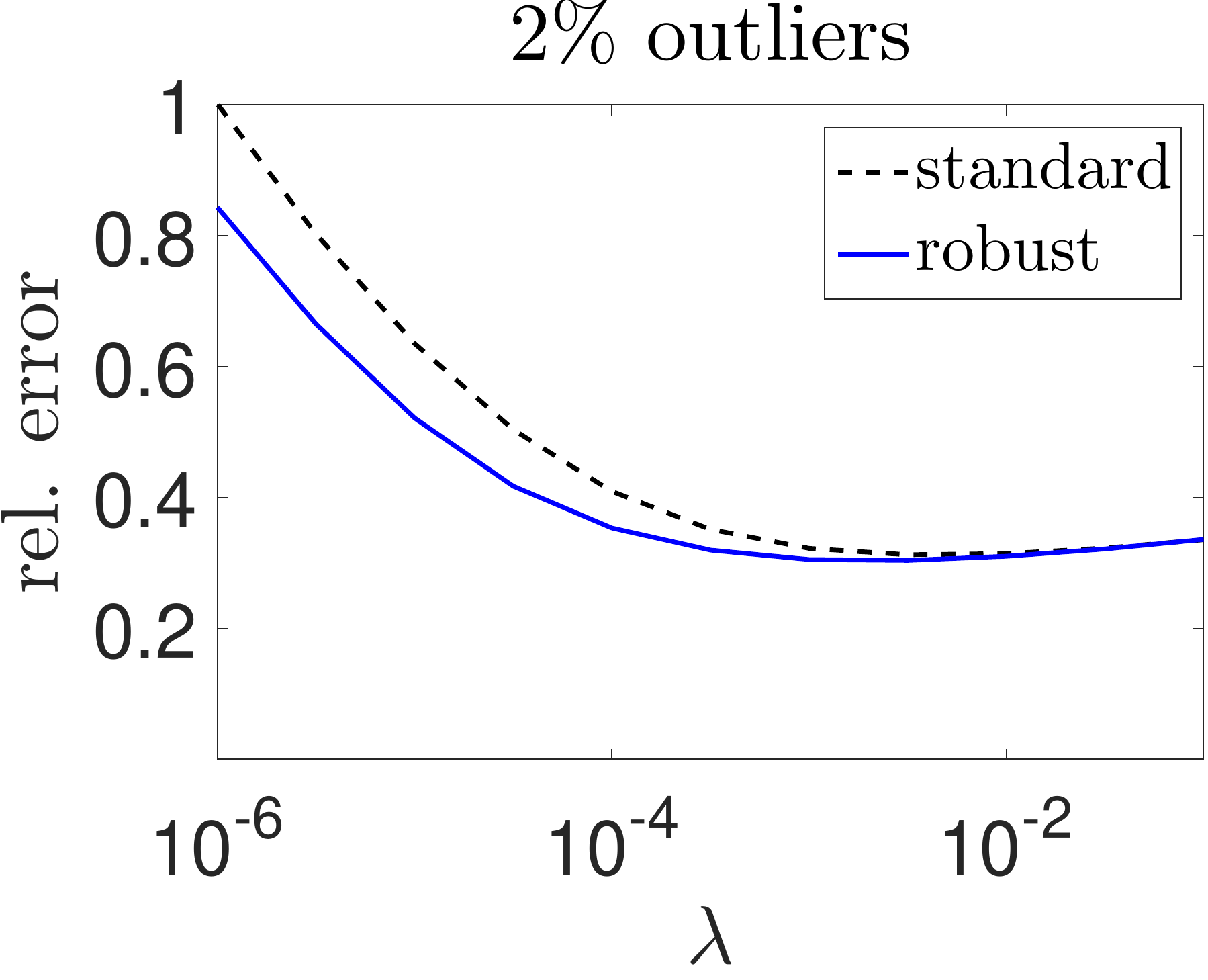}
\hspace*{.2cm}
\includegraphics[width=.31\textwidth]{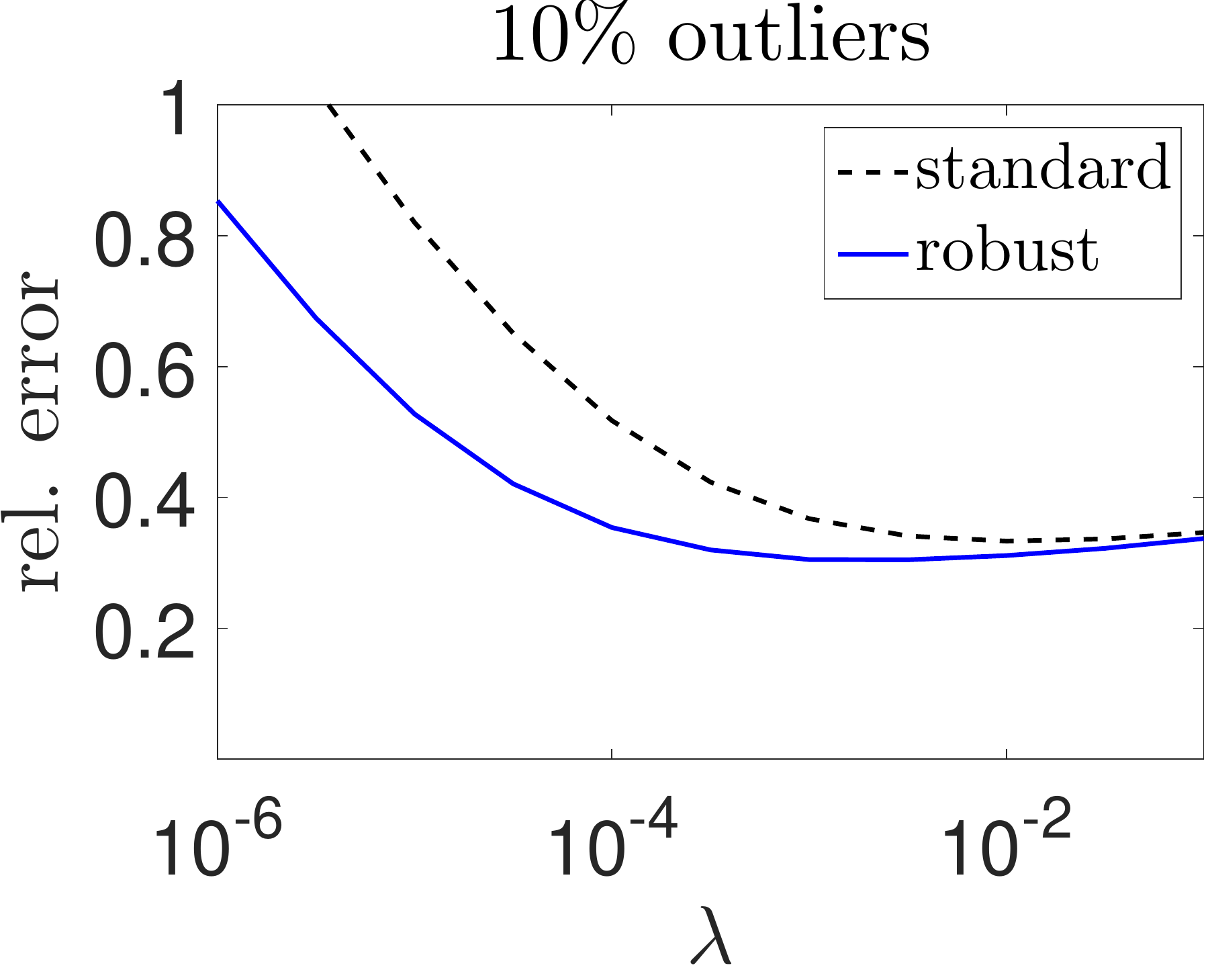}
\caption{Carbon ash multi-frame}\label{fig:sliding_curves_4}
\end{subfigure}
\caption{Semiconvergence curves -- dependence of the relative error of the reconstruction on the size of the regularization parameter $\lambda$ for various percentages of outliers: Talwar \eqref{eq:robust_wls} - \eqref{eq:beta_opt} (solid line) and the standard data fidelity function \eqref{eq:wls} (dashed line). 
}
\end{figure}

\begin{figure}[!th]
\centering
\begin{subfigure}[b]{.18\textwidth}\captionsetup{justification=centering}
\includegraphics[width=.98\textwidth]{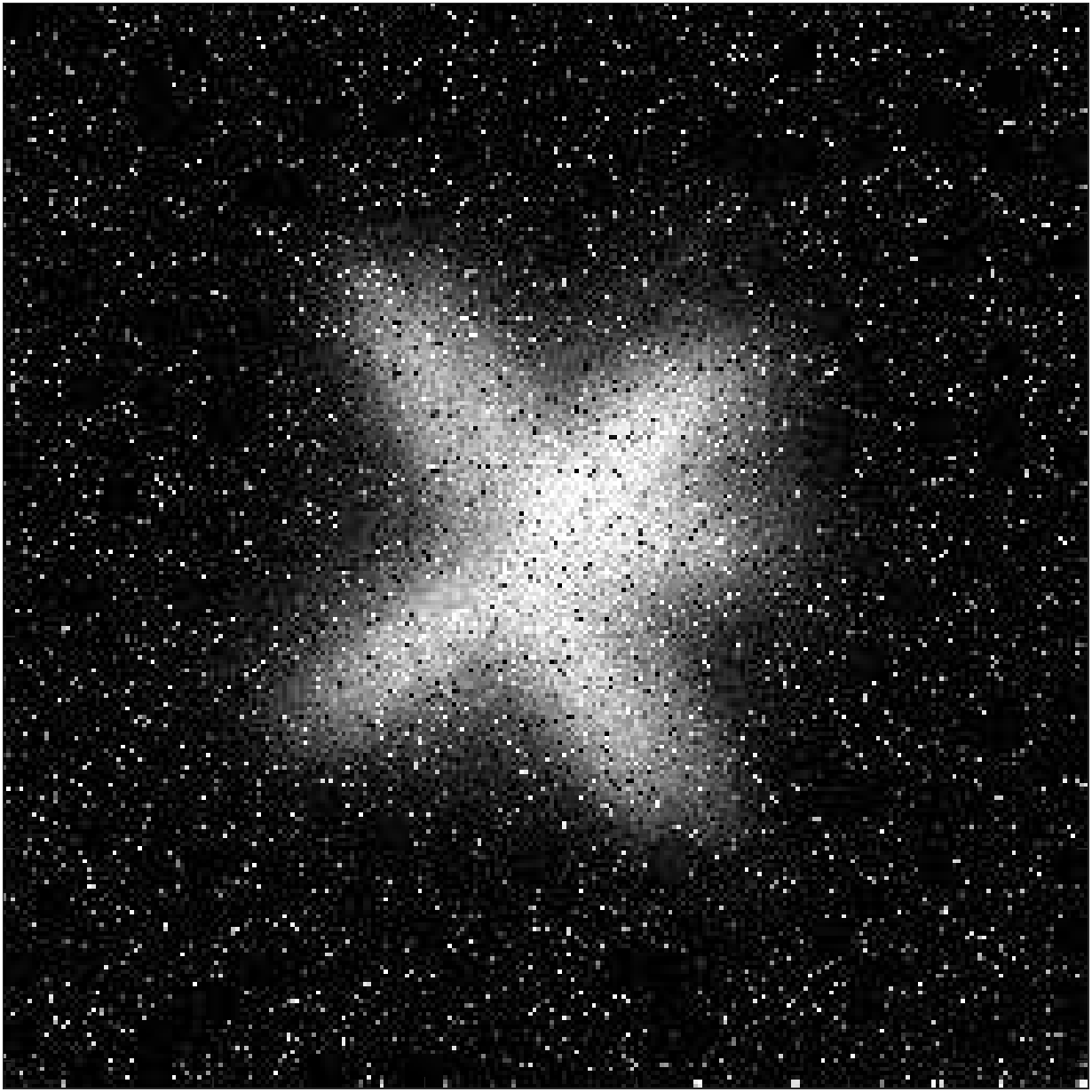} 
\caption{data\newline}
\end{subfigure}
\hspace*{.1cm}
\begin{subfigure}[b]{.18\textwidth}\captionsetup{justification=centering}
\includegraphics[width=.98\textwidth]{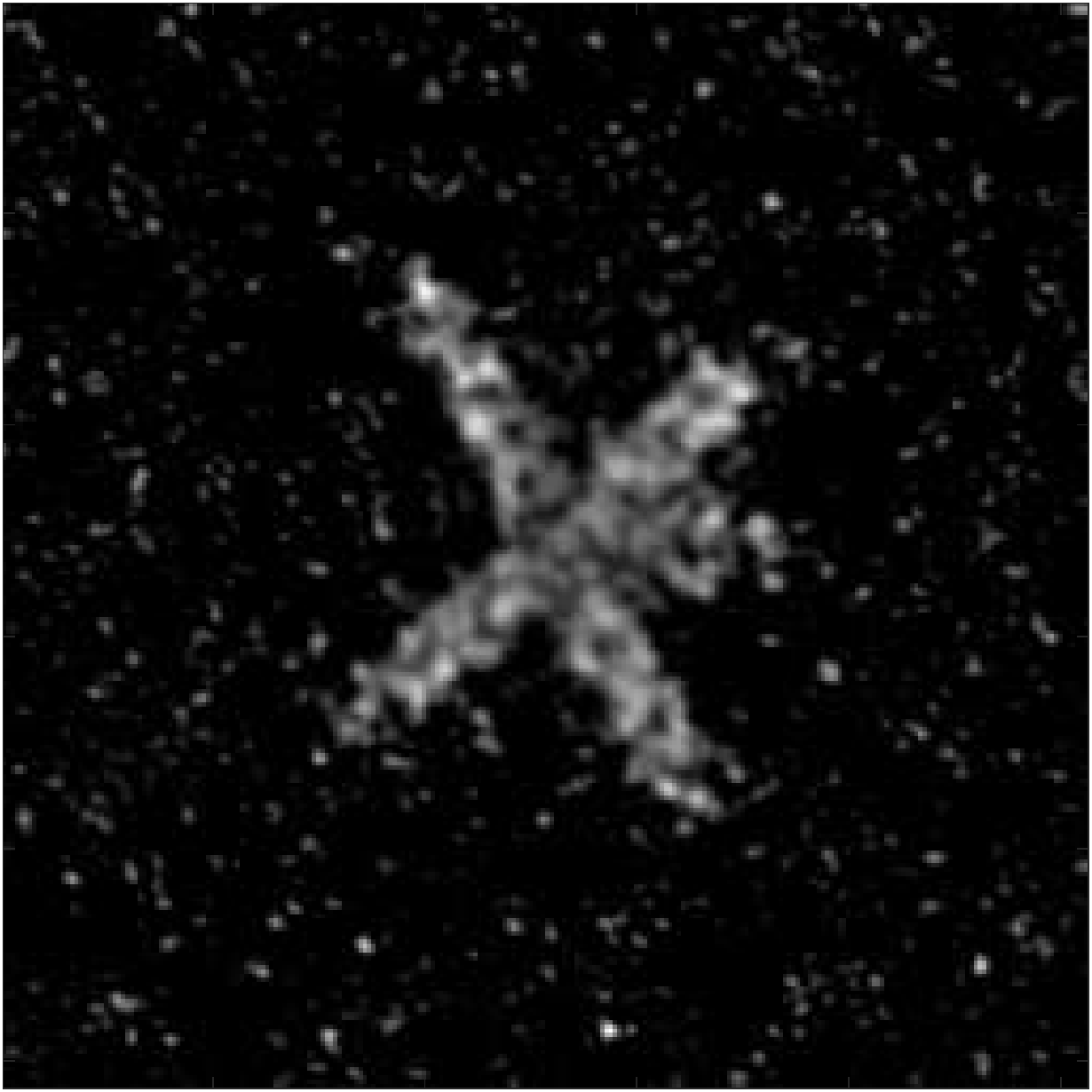} 
\caption{single-frame\\ standard}
\end{subfigure}
\begin{subfigure}[b]{.18\textwidth}\captionsetup{justification=centering}
\includegraphics[width=.98\textwidth]{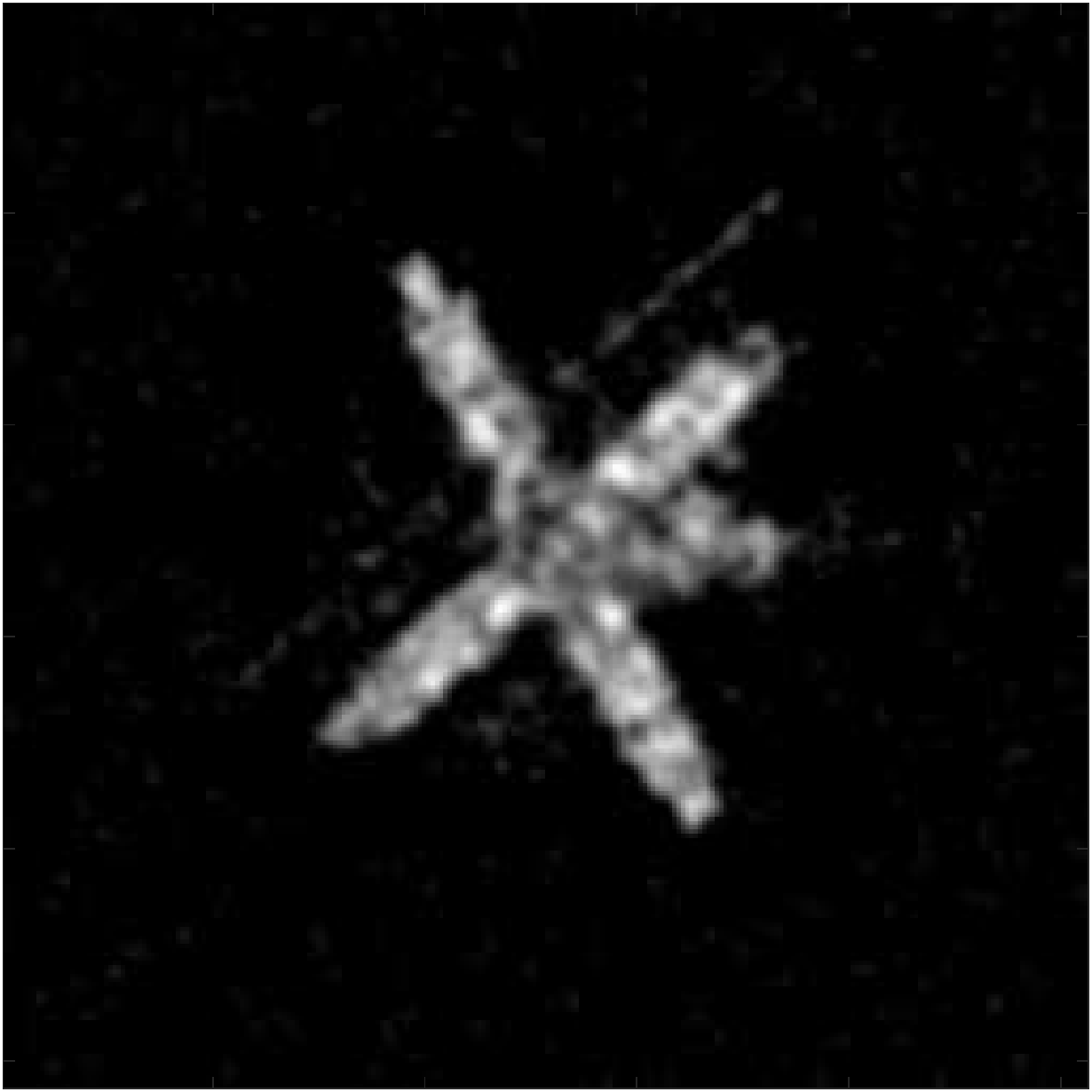} 
\caption{single-frame\\ robust}
\end{subfigure}
\hspace*{.1cm}
\begin{subfigure}[b]{.18\textwidth}\captionsetup{justification=centering}
\includegraphics[width=.98\textwidth]{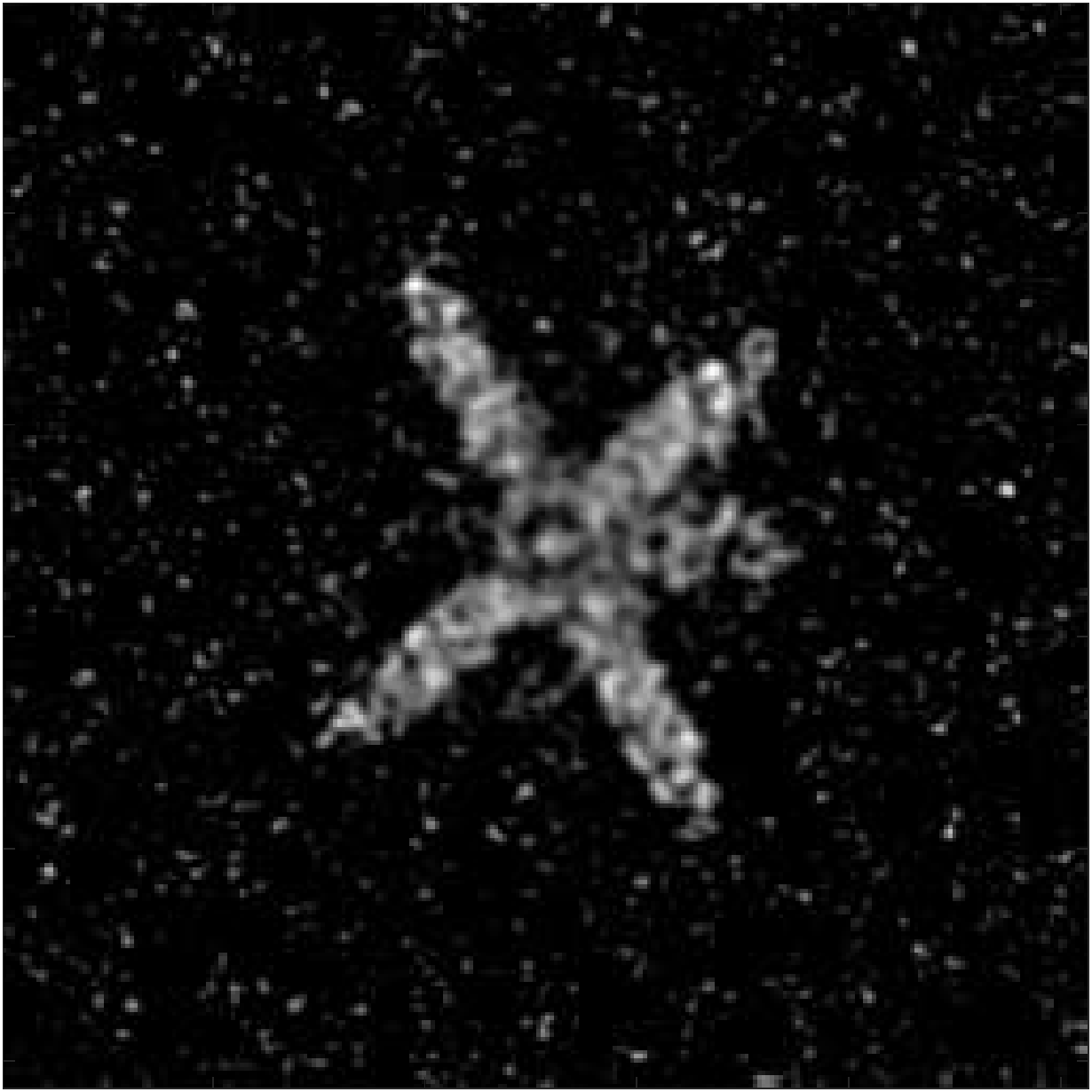} 
\caption{multi-frame,\\ standard}
\end{subfigure}
\begin{subfigure}[b]{.18\textwidth}\captionsetup{justification=centering}
\includegraphics[width=.98\textwidth]{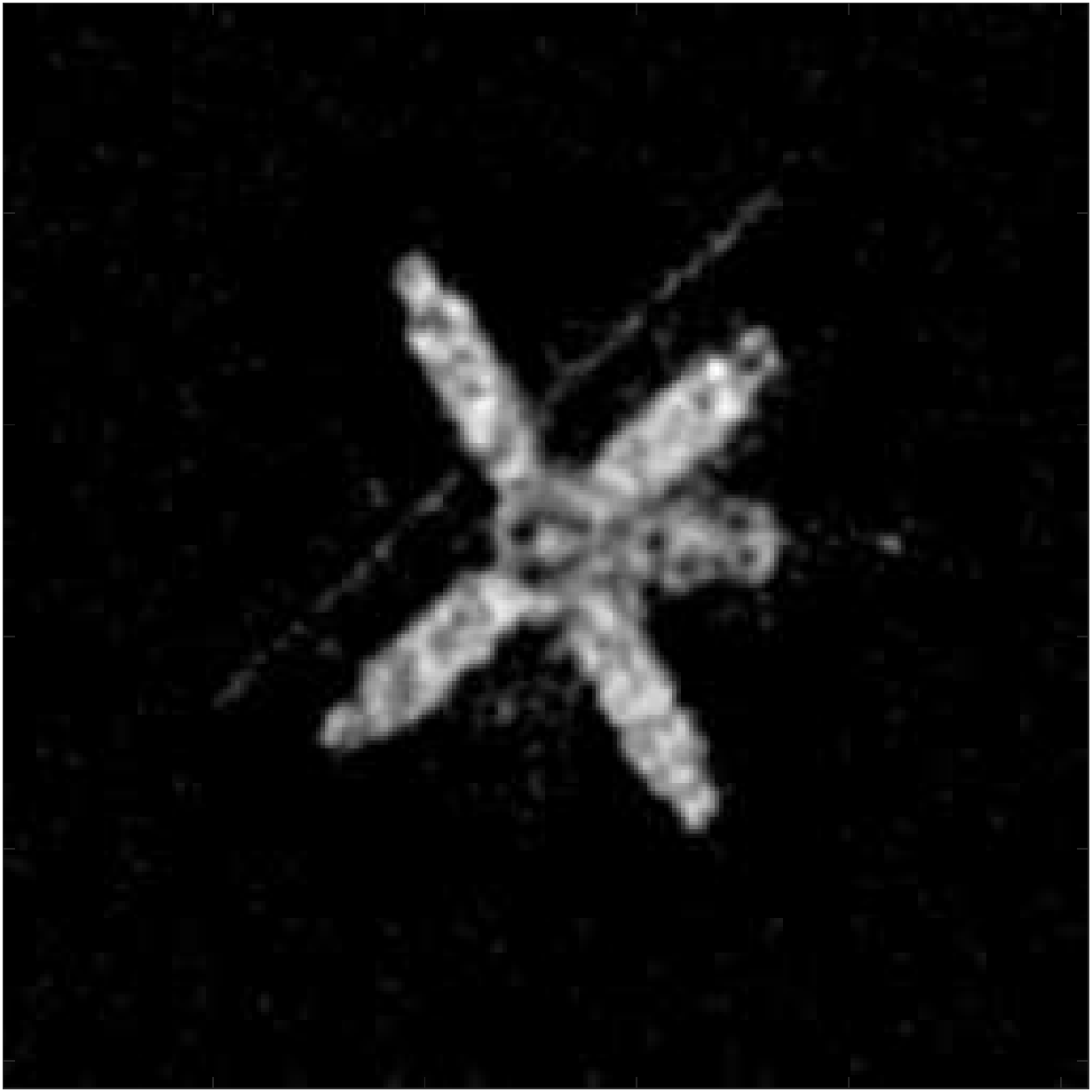} 
\caption{multi-frame \\robust}
\end{subfigure}
\caption{Random corruptions: (a) blurred noisy image with $10 \%$ corrupted pixels (only first frame is shown); (b) - (d) reconstructions corresponding to $\lambda = 10^{-4}$. 
}\label{fig:reconstructions_satellite}
\end{figure}

\begin{figure}[!th]
\centering
\begin{subfigure}[b]{.18\textwidth}\captionsetup{justification=centering}
\includegraphics[width=.98\textwidth]{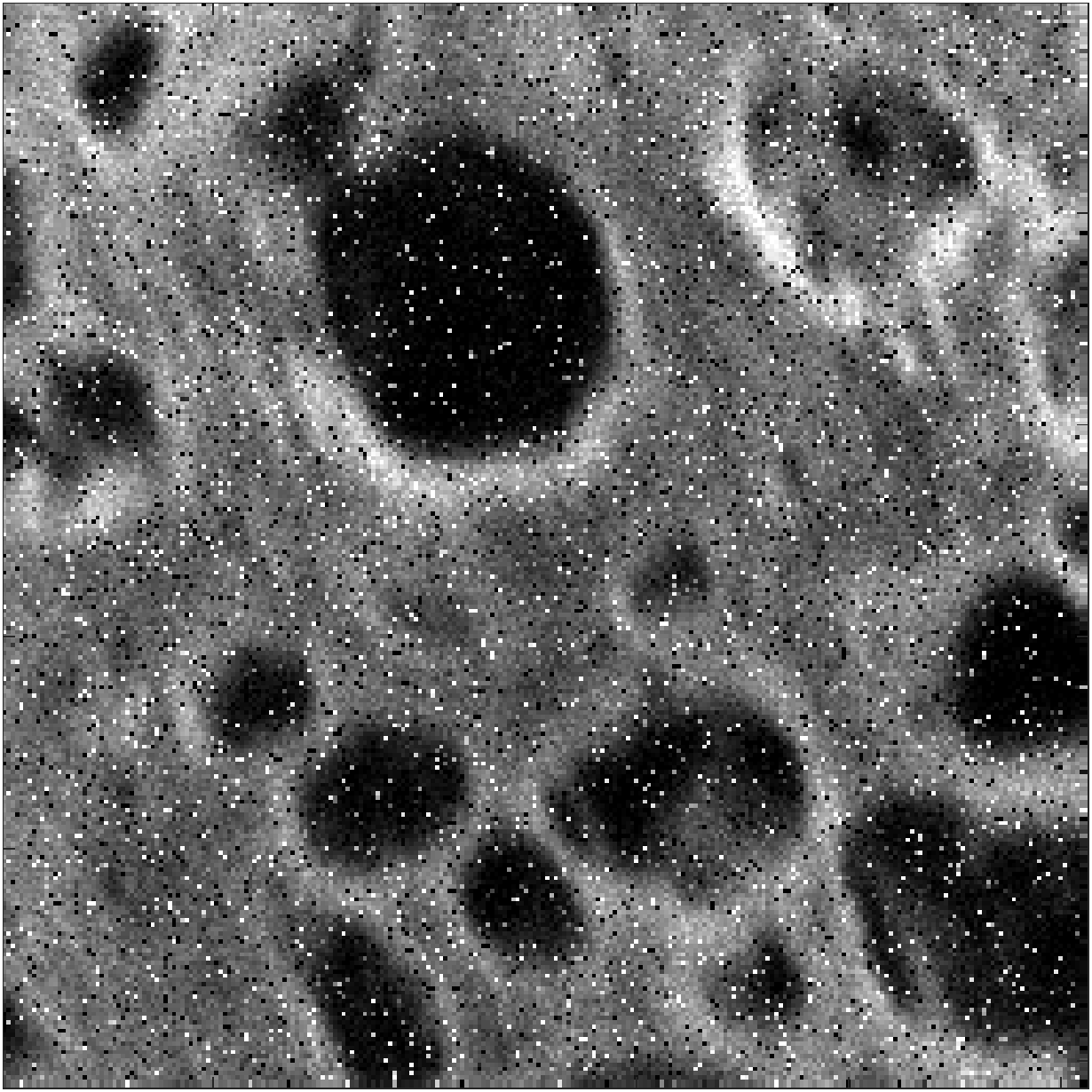} 
\caption{data\newline}
\end{subfigure}
\hspace*{.1cm}
\begin{subfigure}[b]{.18\textwidth}\captionsetup{justification=centering}
\includegraphics[width=.98\textwidth]{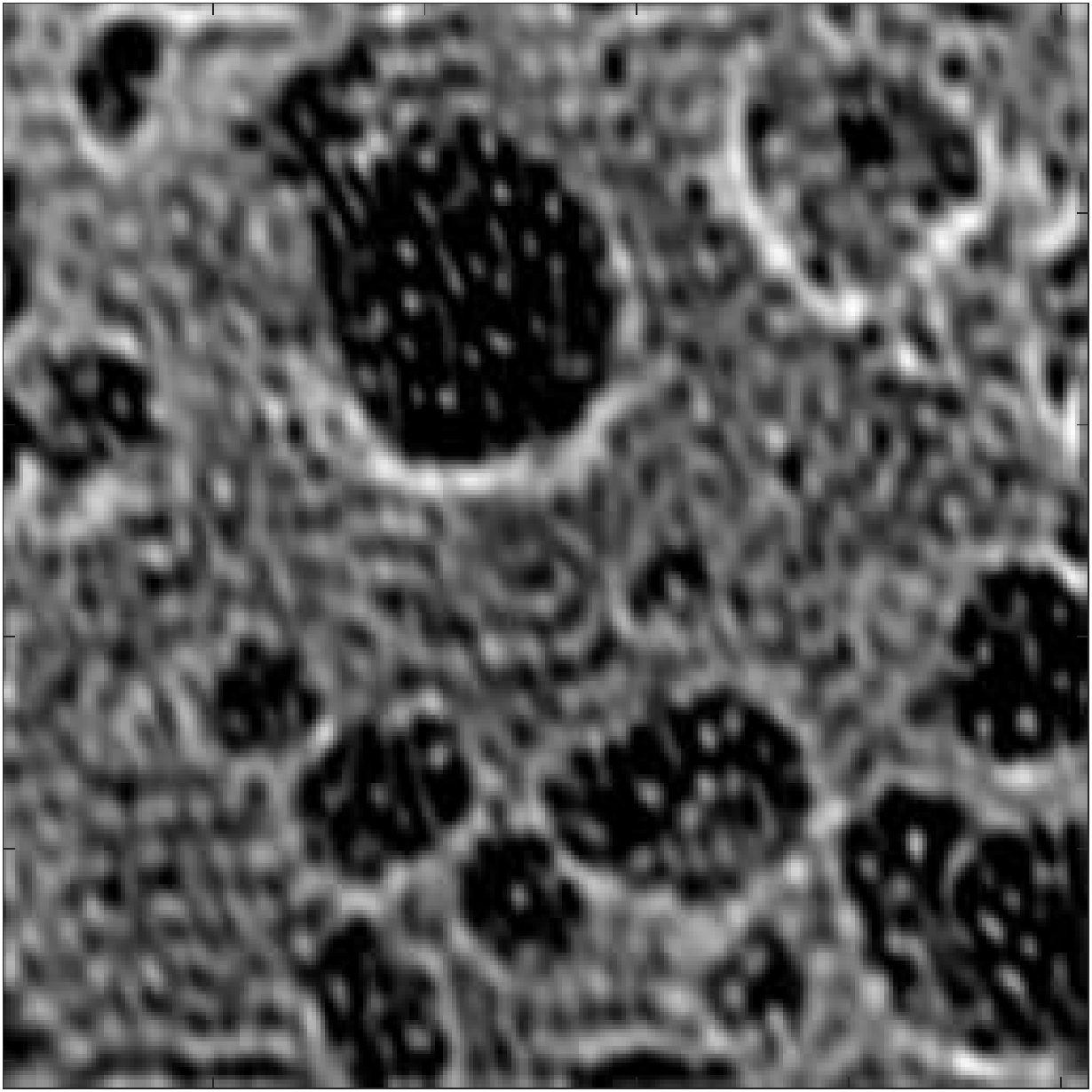} 
\caption{single-frame\\ standard}
\end{subfigure}
\begin{subfigure}[b]{.18\textwidth}\captionsetup{justification=centering}
\includegraphics[width=.98\textwidth]{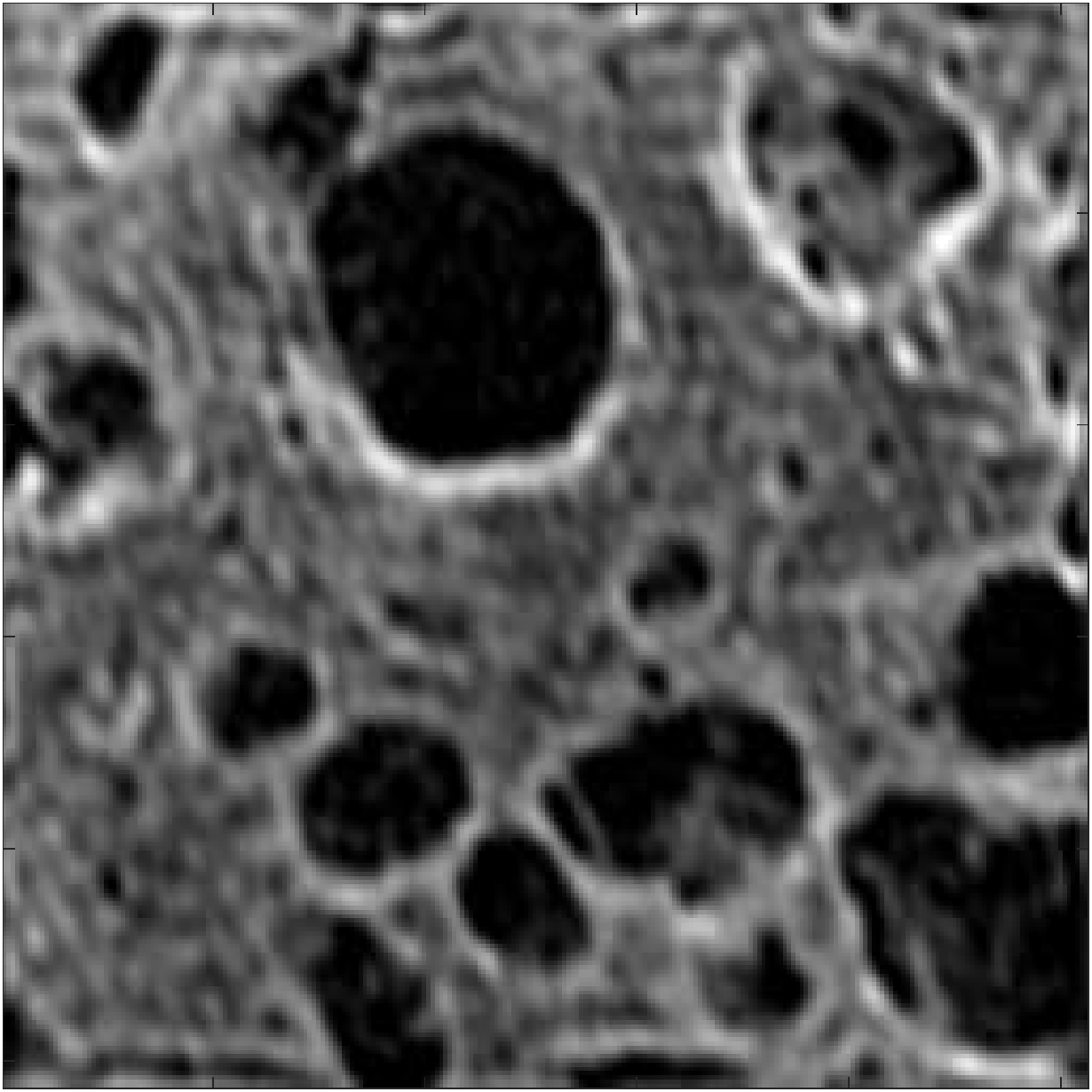} 
\caption{single-frame\\ robust}
\end{subfigure}
\hspace*{.1cm}
\begin{subfigure}[b]{.18\textwidth}\captionsetup{justification=centering}
\includegraphics[width=.98\textwidth]{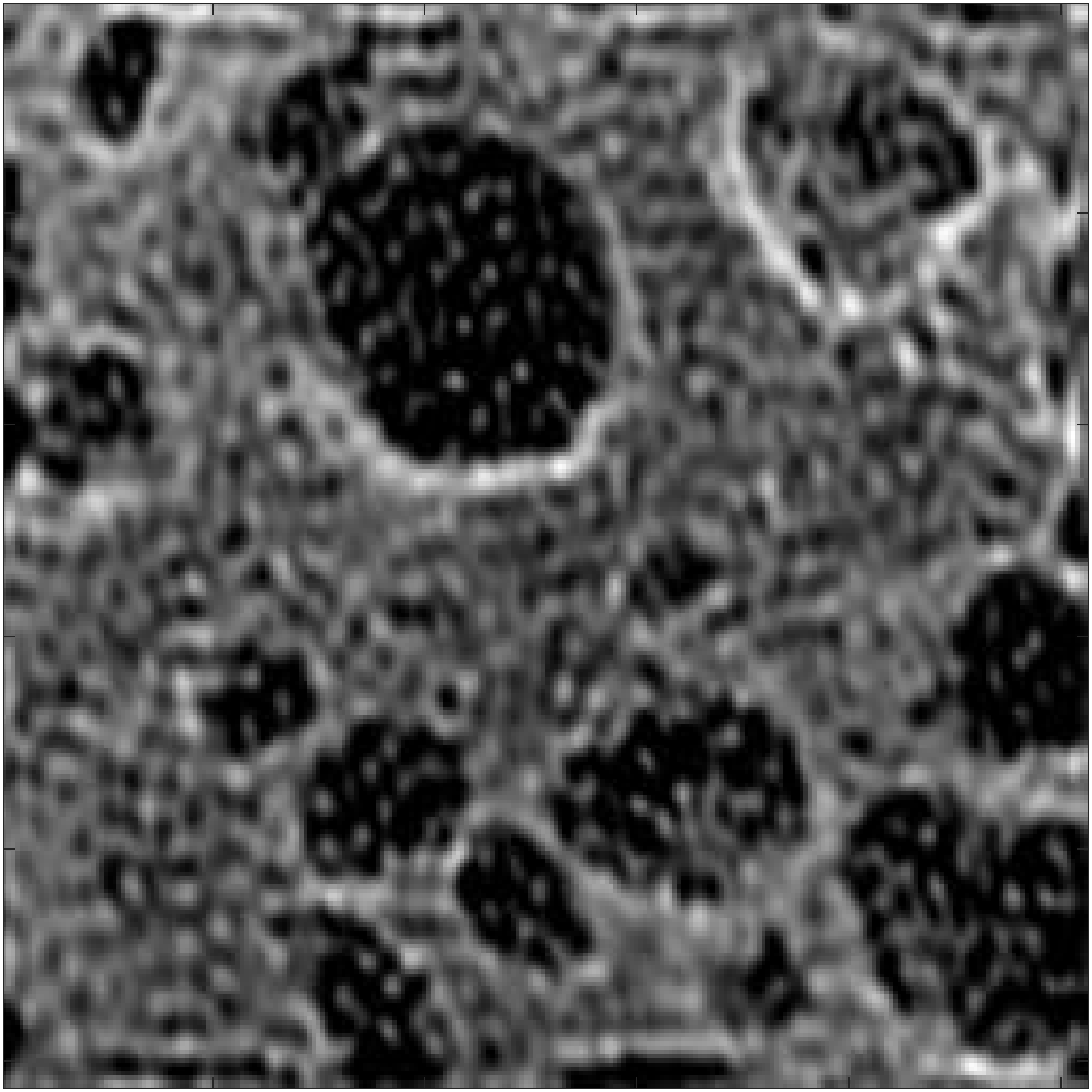} 
\caption{multi-frame,\\ standard}
\end{subfigure}
\begin{subfigure}[b]{.18\textwidth}\captionsetup{justification=centering}
\includegraphics[width=.98\textwidth]{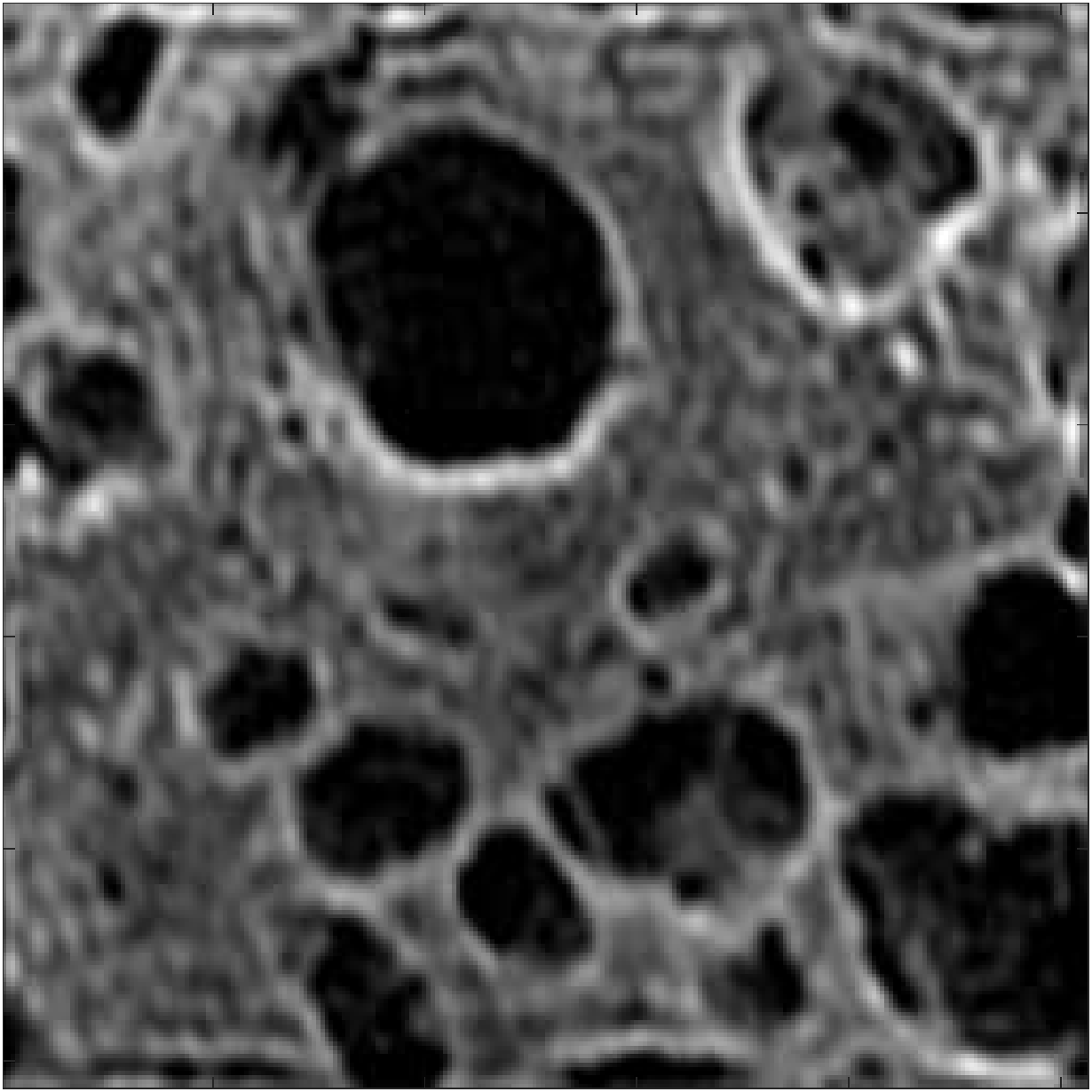} 
\caption{multi-frame \\robust}
\end{subfigure}
\caption{Random corruptions: (a) blurred noisy image with $10 \%$ corrupted pixels (only first frame is shown); (b) - (d) reconstructions corresponding to $\lambda = 10^{-3}$. 
}\label{fig:reconstructions_carbonash}
\end{figure}

First we consider the most simple case of the outliers -- a given percentage of pixels is corrupted at random. These corruptions are generated artificially by adding a value randomly chosen between $0$ and $\max(Ax_\text{true})$ to the given percentage of pixels. The location of these pixels is also chosen randomly. Figures \ref{fig:sliding_curves_1}, \ref{fig:sliding_curves_2}, \ref{fig:sliding_curves_3}, and \ref{fig:sliding_curves_4}  show semiconvergence curves\footnote{For ill-posed problems,
the relative error of an iterative method generally does not decrease monotonically. Instead, unless the problem is highly over-regularized, the relative errors decrease in the early iterations, but at later iterations the noise and
other errors tend to corrupt the approximations.  This behavior, where the relative errors decrease to a certain
level and then increase at later iterations, is referred to as {\em semiconvergence}; for
more information, we refer readers to \cite{Engl2000Regularization,Hansen2010Discrete,Mueller2012Linear,Vogel2002Computational}.}, representing the dependence of the error on the regularization parameter $\lambda$, when we increase the percentage of corrupted pixels.  
It is no surprise that when outliers occur, more regularization is needed in order to obtain a reasonable approximation of the true image $x_\text{true}$. This is however not the case if we use loss function Talwar, for which the semiconvergence curve remains the same even with increasing percentage of outliers, and therefore no adjustment of the regularization parameter is needed. In Figures \ref{fig:reconstructions_satellite} and \ref{fig:reconstructions_carbonash}, we show the reconstructions corresponding to $10 \%$ outliers. The regularization paramter is chosen as a close-to-the-optimal regularization parameter for the same problem with no outliers. Note that Figures \ref{fig:reconstructions_satellite} and \ref{fig:reconstructions_carbonash} show only one frame for illustration. In the multi-frame case, the corruptions look similar for all frames, except that the random locations of the outliers is different. For random outliers like this, robust regression is clearly superior to standard weighted least squares. The influence of the outliers in the multi-frame case is less severe, due to intrinsic regularization of the overdetermined system (\ref{eq:inverse_problem}). A more comprehensive comparison of the standard and robust approach is shown in Table \ref{tab:percent}, giving the percentage of cases in which the robust approach provides better reconstruction. The robust approach provides better reconstruction in all cases except for the test problem Satellite with no outliers, where the standard approach gave sometimes slightly better reconstructions. However, even in these cases we observed that the difference between the errors of the reconstructions is rather negligible, about $3\%$.
\begin{table}[!th]
\centering
\caption{Comparison of the quality of reconstruction for the standard vs. the robust approach. For each test problem and each percentage of outliers, the results are averaged over 100 independent positions and sizes of random corruptions. Regularization parameters are chosen identically as in Figures \ref{fig:reconstructions_satellite} and \ref{fig:reconstructions_carbonash}. Reconstructions are considered to be of the same quality if the difference between the corresponding relative errors is smaller than 1\%. 
}\label{tab:percent}
{\small \begin{tabular}{lcccc} \toprule 
\multicolumn{5}{c}{better reconstruction: robust/same/standard} \\ 
problem/\% outliers & \multicolumn{1}{c}{0\%} & \multicolumn{1}{c}{1\%} & \multicolumn{1}{c}{2\%} & \multicolumn{1}{c}{5\%} \\ \midrule
Satellite single-frame& 0/93/7& 100/0/0& 100/0/0& 100/0/0\\ 
Satellite multi-frame& 0/94/6& 100/0/0& 100/0/0& 100/0/0\\ 
Carbon ash single-frame& 0/100/0& 100/0/0& 100/0/0& 100/0/0\\ 
Carbon ash  multi-frame& 0/100/0& 100/0/0& 100/0/0& 100/0/0\\ 
\bottomrule\end{tabular}}
\end{table}

\subsubsection*{Added object with different blurring}
\begin{figure}[!th]
\centering
\begin{subfigure}[b]{.18\textwidth}\captionsetup{justification=centering}
\includegraphics[width=.98\textwidth]{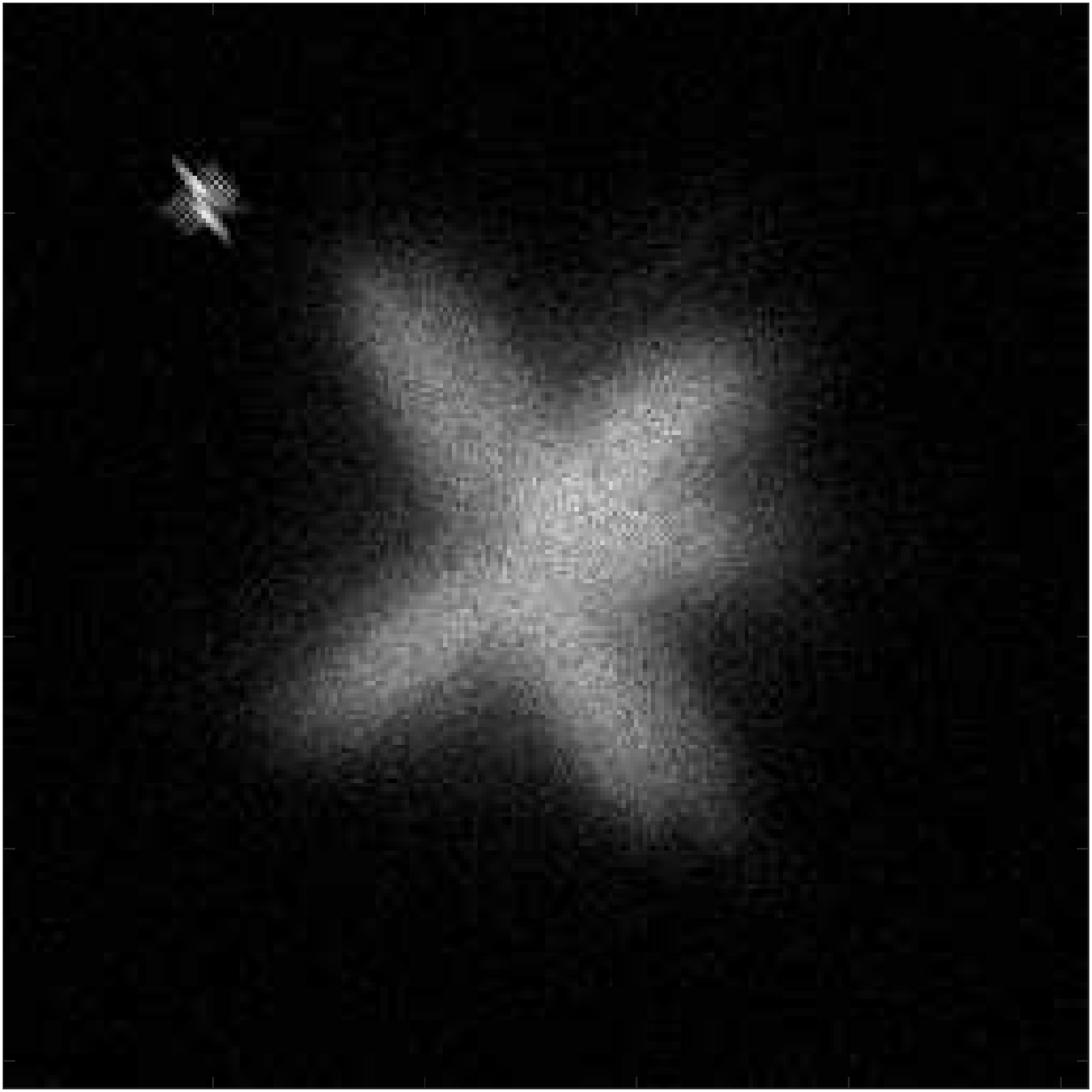} 
\caption{data\newline}
\end{subfigure}
\hspace*{.1cm}
\begin{subfigure}[b]{.18\textwidth}\captionsetup{justification=centering}
\includegraphics[width=.98\textwidth]{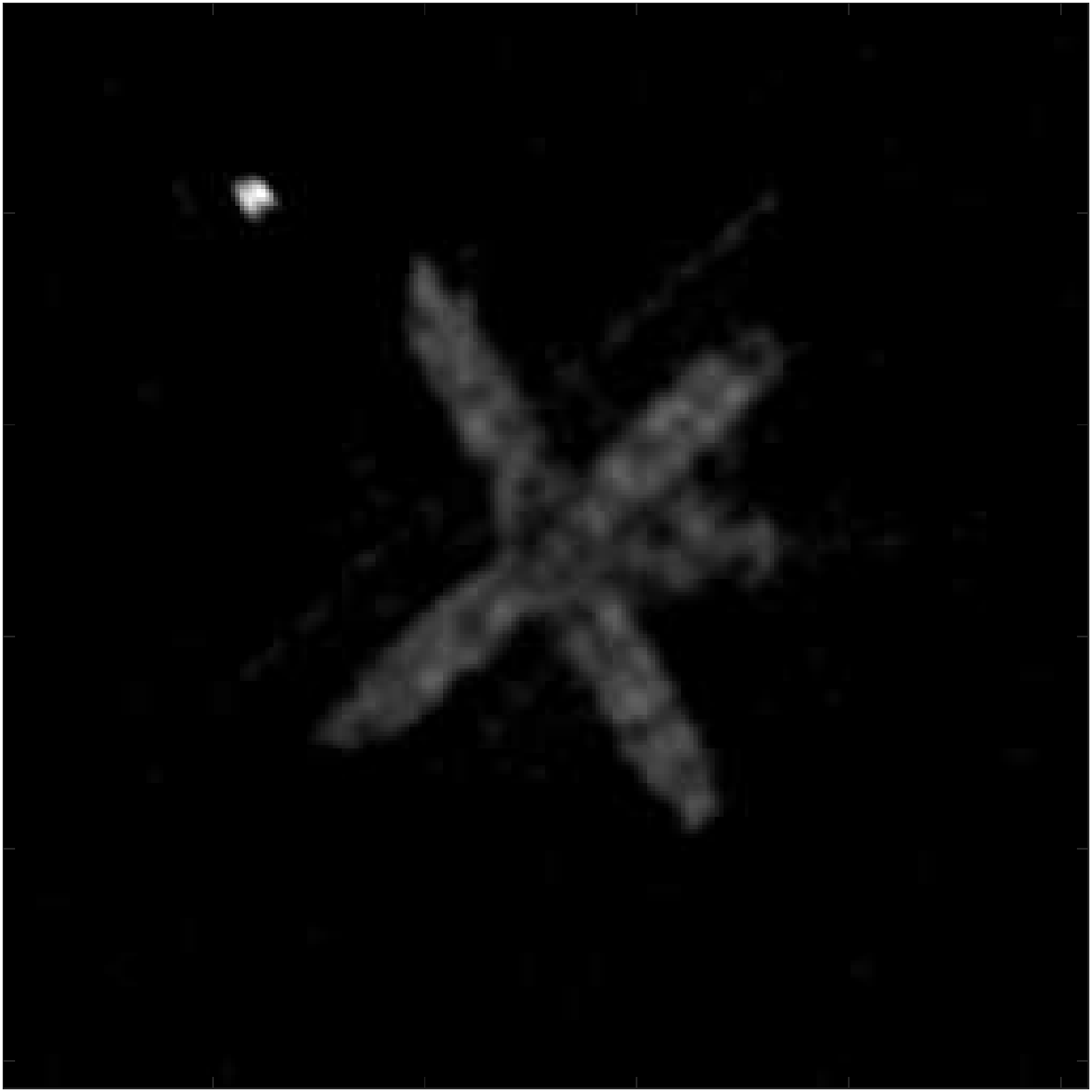} 
\caption{single-frame\\ standard}
\end{subfigure}
\begin{subfigure}[b]{.18\textwidth}\captionsetup{justification=centering}
\includegraphics[width=.98\textwidth]{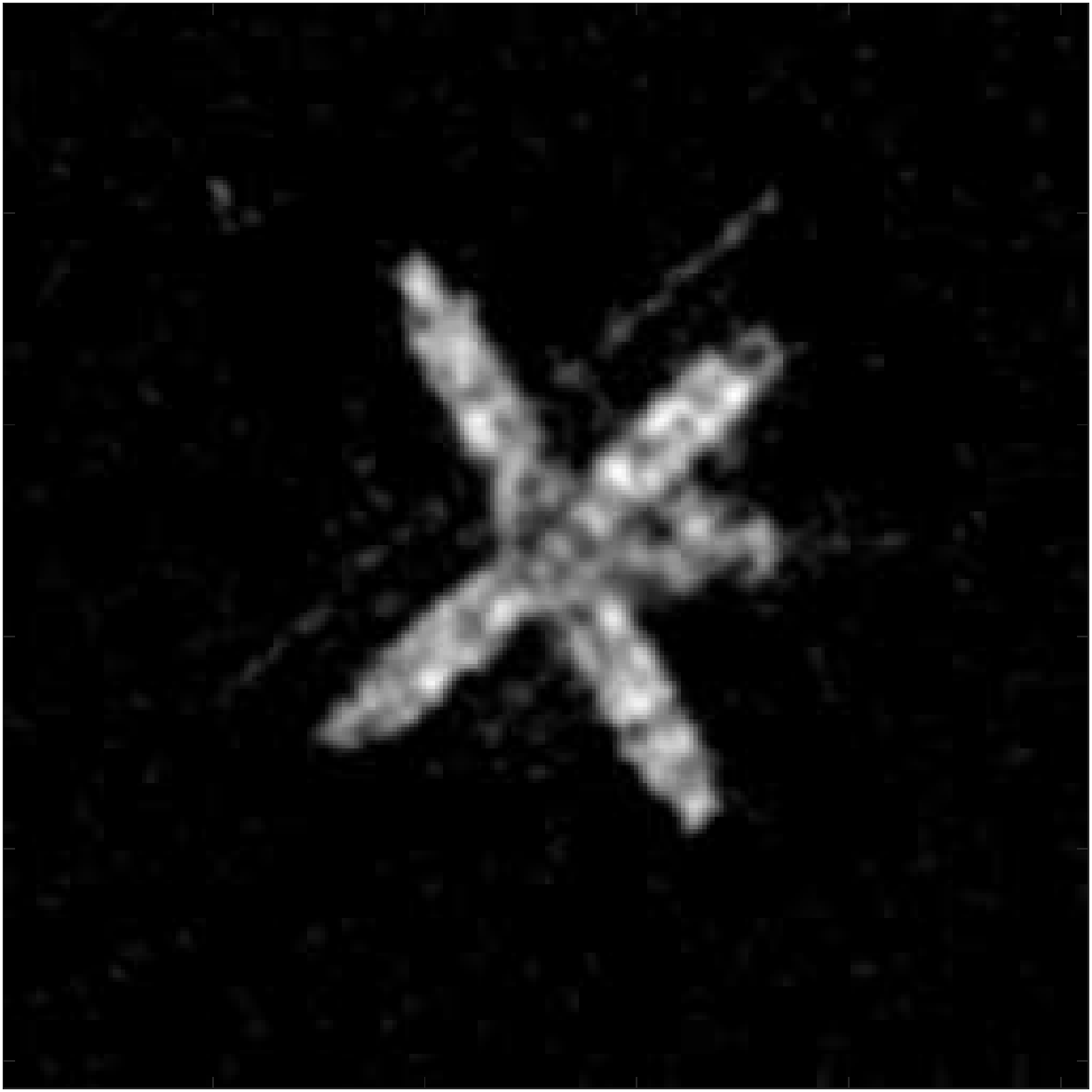} 
\caption{single-frame\\ robust}
\end{subfigure}
\hspace*{.1cm}
\begin{subfigure}[b]{.18\textwidth}\captionsetup{justification=centering}
\includegraphics[width=.98\textwidth]{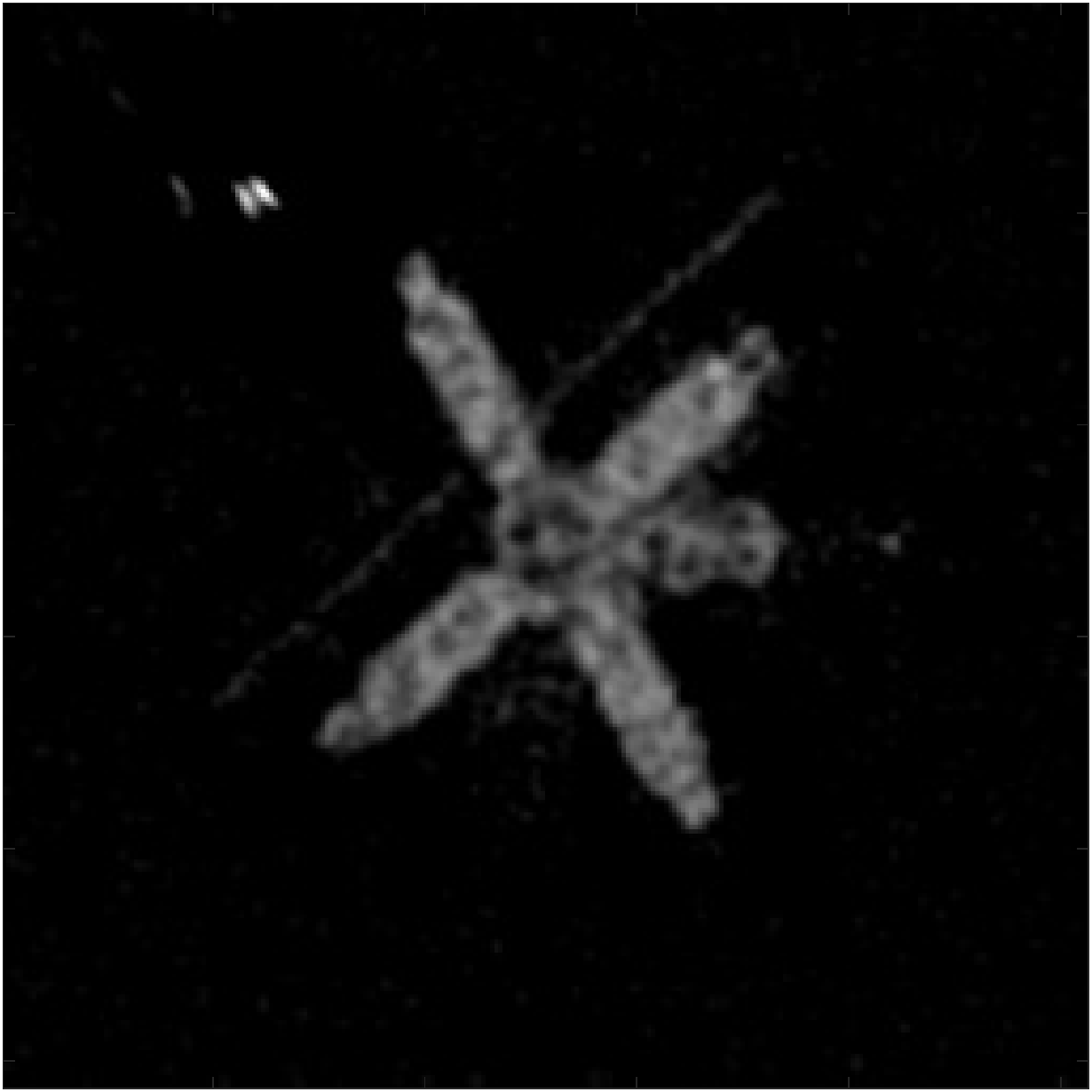} 
\caption{multi-frame,\\ standard}
\end{subfigure}
\begin{subfigure}[b]{.18\textwidth}\captionsetup{justification=centering}
\includegraphics[width=.98\textwidth]{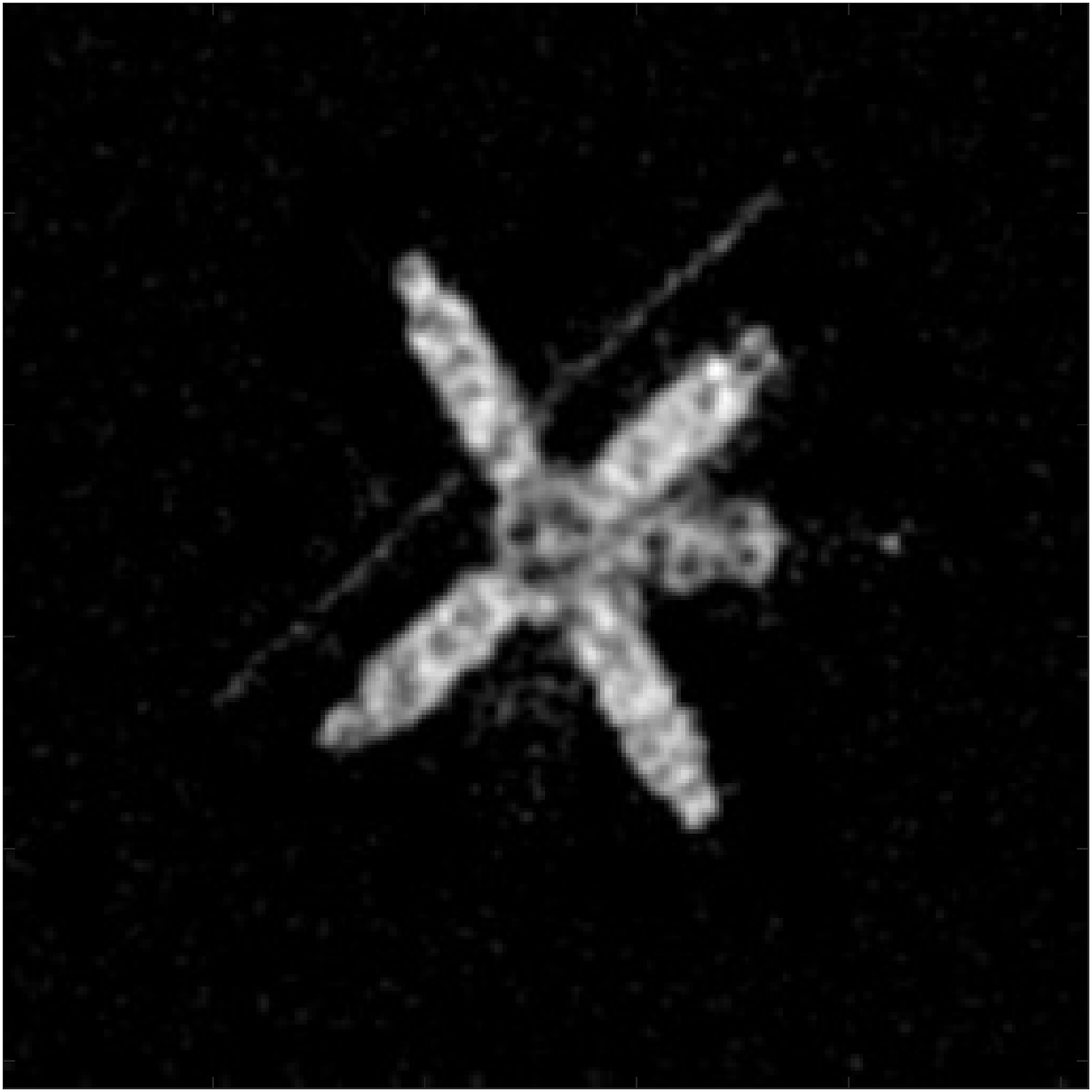} 
\caption{multi-frame \\robust}
\end{subfigure}
\caption{Added object: (a) blurred noisy image with a small object added to the first frame (only first frame is shown); (b) - (d) reconstructions corresponding to $\lambda = 10^{-4}$. 
}\label{fig:alien}
\end{figure}

We also consider a situation when a small object appears in the scene, but is blurred by a different PSF than the main object (satellite). The aim is to recover the main object, while suppressing the influence of the added one. In our case, the added object is a small satellite in the left upper corner that is blurred by a small motion blur. In the multi-frame case, the small satellite is added to the first frame only. The difference between the reconstructions using standard and robust approach is shown in Figure \ref{fig:alien}. For the single frame problem, reconstructions obtained using the standard loss function is fully dominated by the small added object. For the multi-frame situation, the influence of the outlier is somewhat compensated by the two frames without outliers. In both cases, however, robust regression provides better reconstruction, comparable to the reconstruction from the data without outliers.

\subsubsection*{Outliers introduced by boundary conditions}
\begin{figure}[!th]
\centering
\begin{subfigure}[b]{.18\textwidth}\captionsetup{justification=centering}
\includegraphics[width=.98\textwidth]{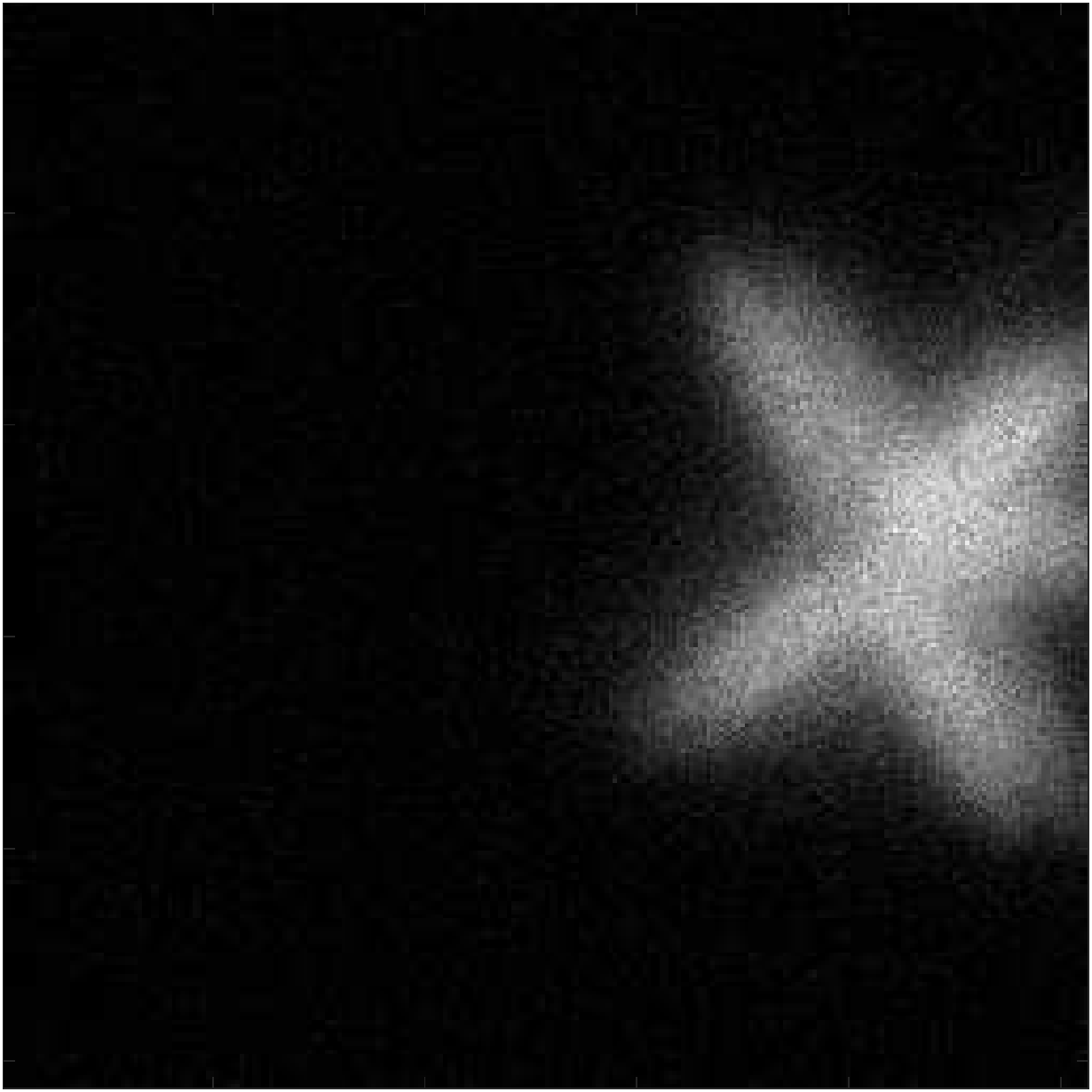} 
\caption{data\newline}
\end{subfigure}
\hspace*{.1cm}
\begin{subfigure}[b]{.18\textwidth}\captionsetup{justification=centering}
\includegraphics[width=.98\textwidth]{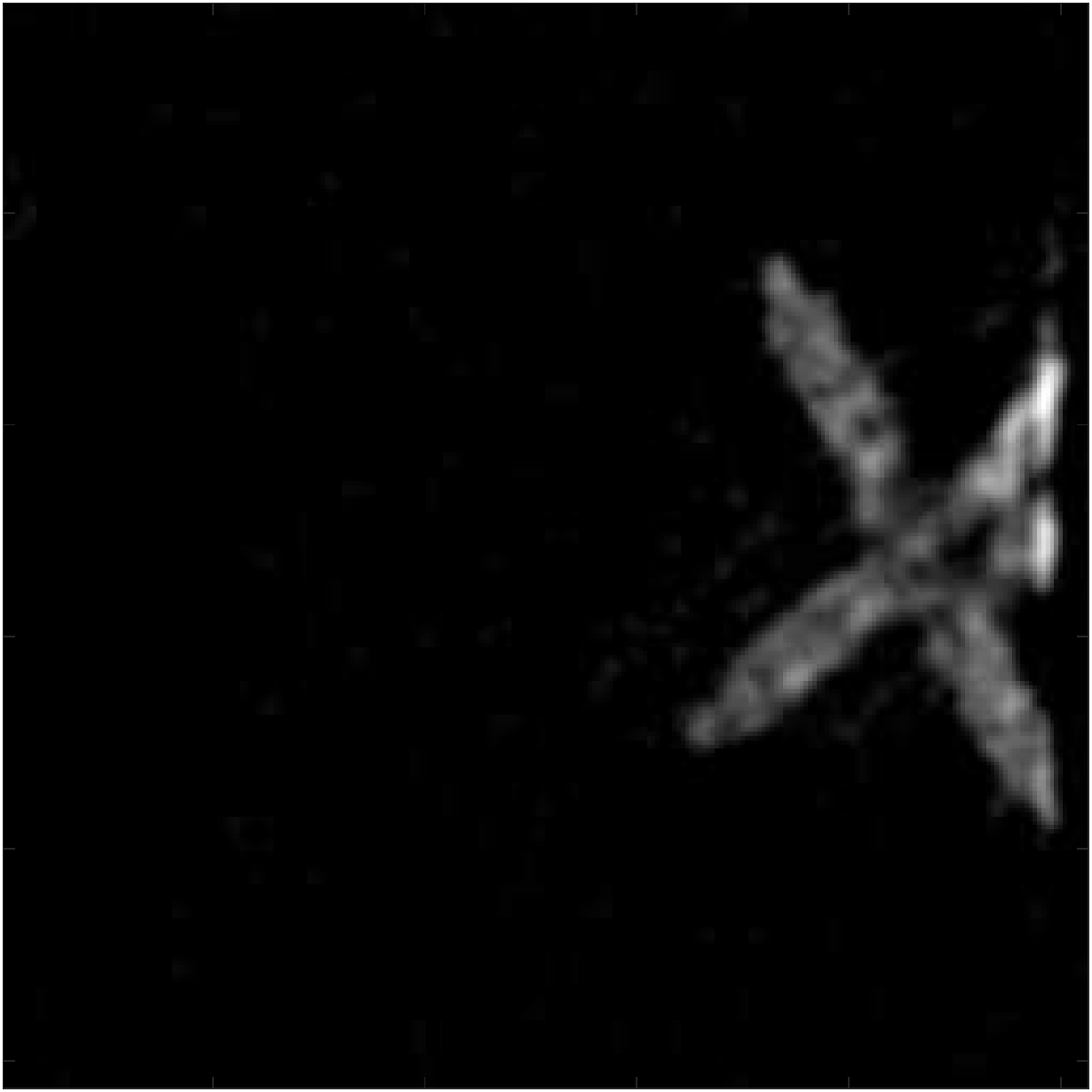} 
\caption{single-frame\\ standard}
\end{subfigure}
\begin{subfigure}[b]{.18\textwidth}\captionsetup{justification=centering}
\includegraphics[width=.98\textwidth]{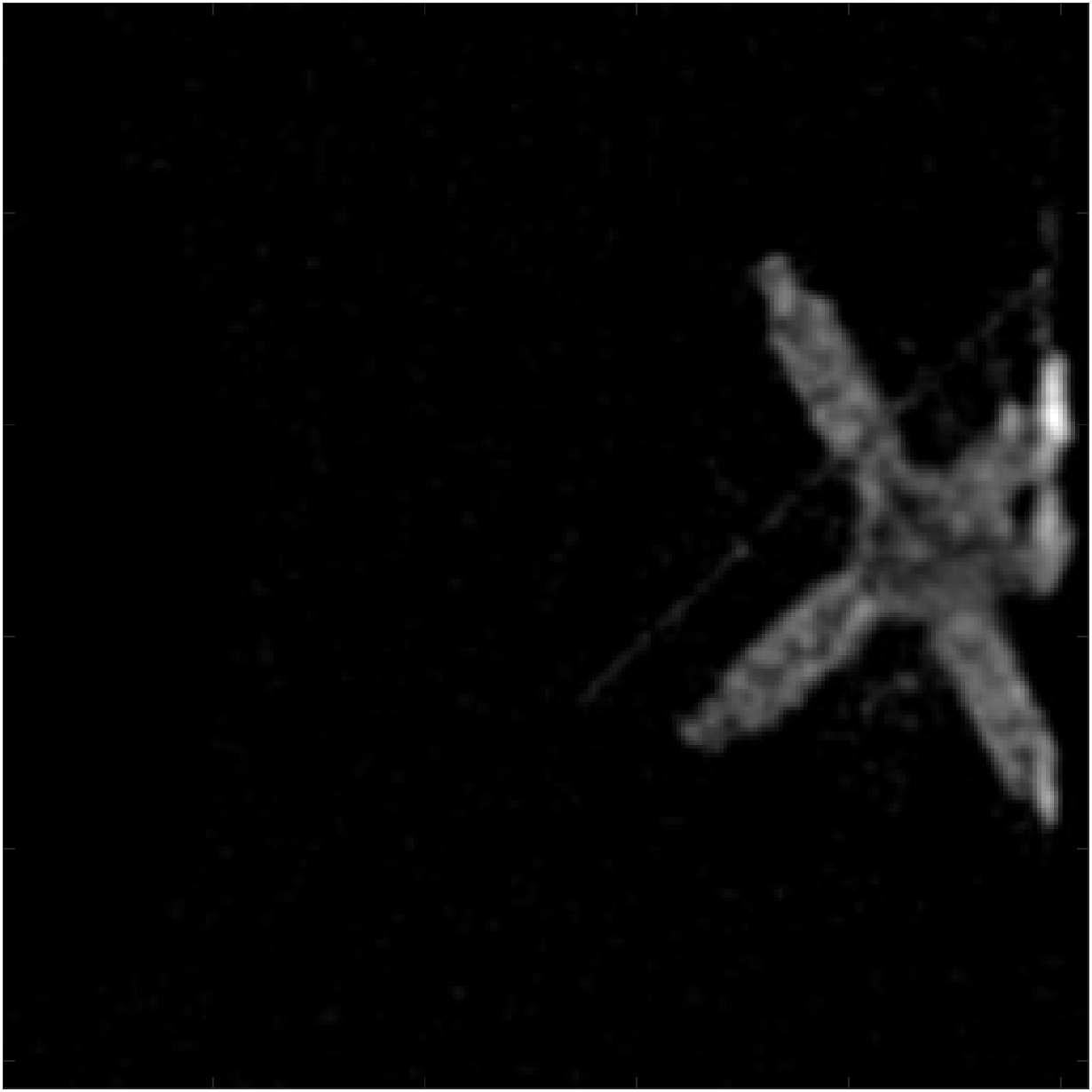} 
\caption{single-frame\\ robust}
\end{subfigure}
\hspace*{.1cm}
\begin{subfigure}[b]{.18\textwidth}\captionsetup{justification=centering}
\includegraphics[width=.98\textwidth]{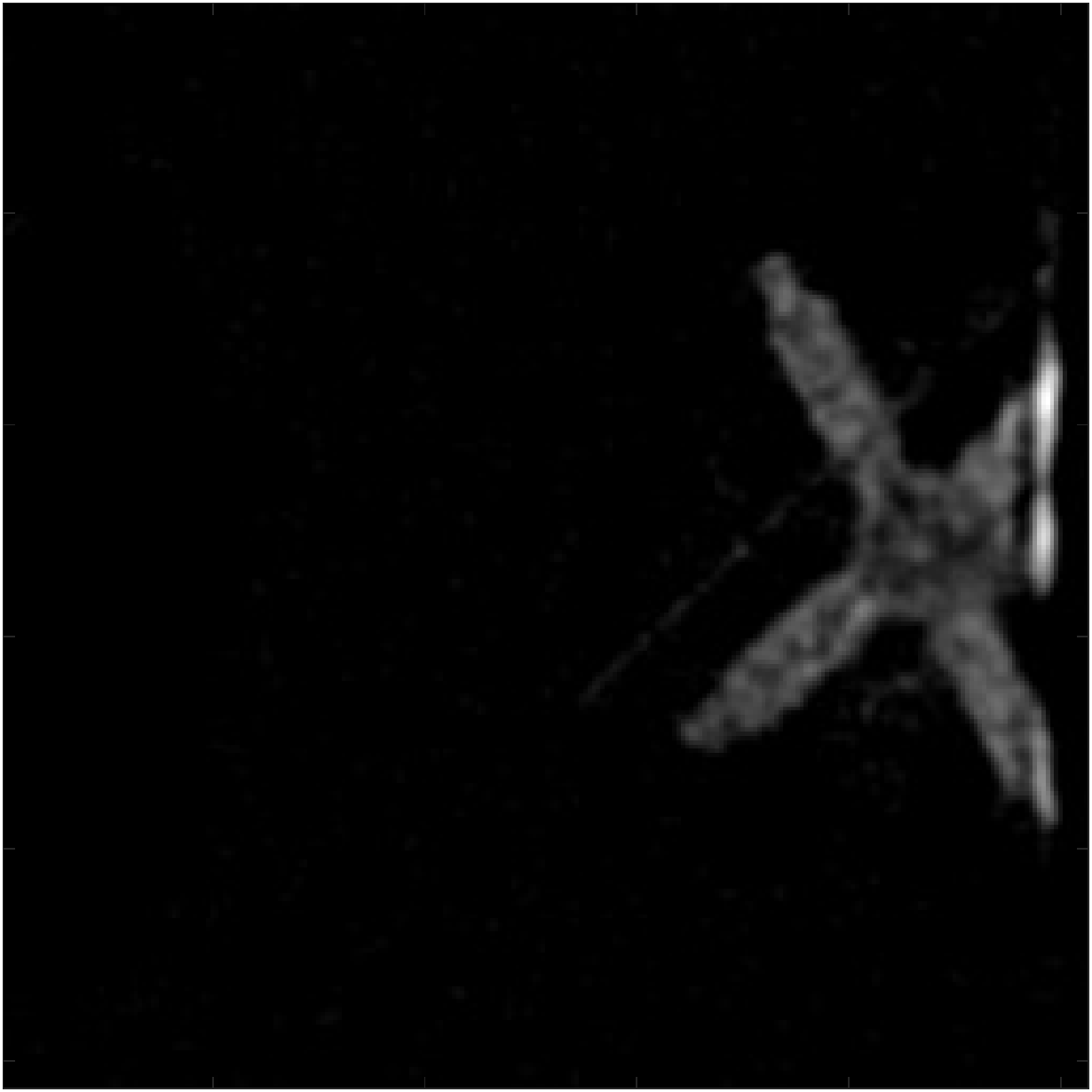} 
\caption{multi-frame,\\ standard}
\end{subfigure}
\begin{subfigure}[b]{.18\textwidth}\captionsetup{justification=centering}
\includegraphics[width=.98\textwidth]{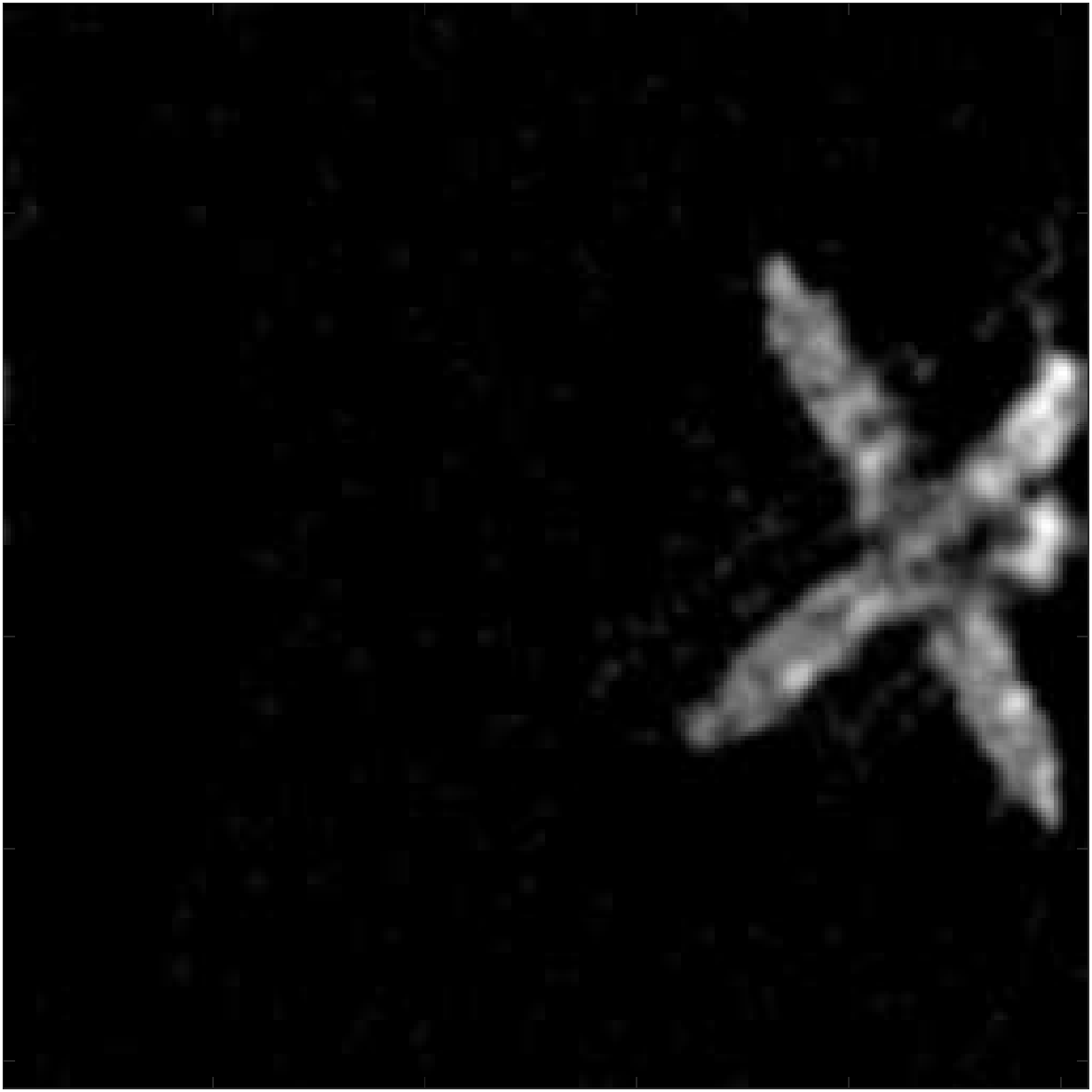} 
\caption{multi-frame \\robust}
\end{subfigure}
\caption{Incorrectly imposed periodic boundary conditions: (a) blurred noisy image close to the edge (only first frame is shown); (b) - (d) reconstructions corresponding to $\lambda = 10^{-4}$. 
}\label{fig:bc}
\end{figure}

Defining the boundary conditions plays an important role in solving image deblurring problems. As is well known, see e.g. \cite{Hansen2006Deblurring}, unless some strong a priori information about the scene outside the borders is available, any choice of the boundary conditions may lead to artefacts around edges in the reconstruction.  Similarly as in \cite{Calef2013Iteratively}, we may expect that the robust objective functional (\ref{eq:functional}) can to some extent compensate for these edge artifacts, i.e., the outliers are represented by the `incorrectly' imposed boundary conditions. In our model we assume periodic boundary conditions, which is computationally very appealing, since it allows evaluating the multiplication by $A$ very efficiently using the fast Fourier transform. However, if any of the objects in the scene is close to the boundary, these boundary conditions will most probably cause artifacts. In order to demonstrate the ability of the proposed scheme to eliminate influence of this type of outlier, we shifted the satellite to the right edge of the image. Other settings remain unchanged. Reconstructions using standard and robust approach are shown in Figure \ref{fig:bc}. We see that, although not spectacular, robust regression can reduce the artifacts caused by incorrectly imposed boundary conditions and therefore provide better reconstruction of the true image. Quantitative results for this and all the previous types of outliers are shown in Tables~\ref{tab:2a} and \ref{tab:2b}.

\begin{table}[!th]
\caption{Comparison of the standard and robust approach in terms of relative error of the reconstruction. Each row contains results for the standard and robust approach. Abbreviation `\# it' stands for the number of Newton steps performed before the relative size of the projected gradient reached the tolerance $10^{-4}$. Corresponding reconstructions are shown in Figures~\ref{fig:reconstructions_satellite}--\ref{fig:bc}. 
}
\centering
{\small
\begin{subtable}[h]{\textwidth}
\subcaption{single-frame}
\centering
\begin{tabular}{lcccc}\toprule 
 & \multicolumn{2}{c}{standard} & \multicolumn{2}{c}{robust} \\ 
\multicolumn{1}{c}{problem} & \multicolumn{1}{c}{\# it} & \multicolumn{1}{c}{rel. error} &  \multicolumn{1}{c}{\# it} & \multicolumn{1}{c}{rel. error} \\ \midrule
Satellite& 15& 3.40$\times 10^{-1}$& 16& 3.42$\times 10^{-1}$\\ 
Satellite random corr. 10\%& 14& 6.78$\times 10^{-1}$& 14& 3.57$\times 10^{-1}$\\ 
Carbon ash& 10& 3.10$\times 10^{-1}$& 11& 3.08$\times 10^{-1}$\\ 
Carbon ash random corr. 10\% & 11& 3.80$\times 10^{-1}$& 14& 3.10$\times 10^{-1}$\\ 
Satellite added object& 15& 4.72$\times 10^{-1}$& 15& 3.43$\times 10^{-1}$\\ 
Satellite boundary conditions& 15& 5.48$\times 10^{-1}$& 25& 4.51$\times 10^{-1}$\\ 
\bottomrule\end{tabular}\label{tab:2a}
\end{subtable}

\medskip

\begin{subtable}[h]{1\textwidth}
\subcaption{multi-frame}
\centering
\begin{tabular}{lcccc} \toprule
 & \multicolumn{2}{c}{standard} & \multicolumn{2}{c}{robust} \\ 
\multicolumn{1}{c}{problem} & \multicolumn{1}{c}{\# it} & \multicolumn{1}{c}{rel. error} & \multicolumn{1}{c}{\# it} & \multicolumn{1}{c}{rel. error} \\ \midrule
Satellite& 12& 2.89$\times 10^{-1}$& 11& 2.89$\times 10^{-1}$\\ 
Satellite random corr. 10\%& 11& 6.45$\times 10^{-1}$& 13& 3.00$\times 10^{-1}$\\ 
Carbon ash& 12& 3.07$\times 10^{-1}$& 11& 3.05$\times 10^{-1}$\\ 
Carbon ash random corr. 10\%& 9& 3.70$\times 10^{-1}$& 19& 3.06$\times 10^{-1}$\\ 
Satellite added object& 13& 3.33$\times 10^{-1}$& 11& 2.90$\times 10^{-1}$\\ 
Satellite boundary conditions& 14& 5.26$\times 10^{-1}$& 14& 4.27$\times 10^{-1}$\\ 
\bottomrule\end{tabular}\label{tab:2b}
\end{subtable}}
\end{table} 

\subsection{Generalized cross-validation}\label{sec:gcv_tests}
For the remainder of this section we will only assume the robust approach, i.e., functional \eqref{eq:functional} with the loss function Talwar. In Section~\ref{sec:gcv} we described a regularization parameter selection rule based on leave-one-out cross validation. Since GCV belongs to standard methods, we focus here mainly on the influence of the outliers on its reliability. To obtain various noise levels, we scale the original true scene (with maximum intensity = 255) by 10 and by 100, which results in a decrease of the relative size of Poisson noise. The standard deviation $\sigma$ for the additive Gaussian noise is scaled accordingly by $\sqrt{10}$ and $10$. We compute the resulting signal-to-noise ratio as the reciprocal of the coefficient of variation, i.e.,
\[
\text{SNR} = \frac{\|Ax\|}{\sqrt{\sum_{i=1}^n([Ax]_i + \sigma^2)}}.
\]

For our computations, we use CG to solve \eqref{eq:gcv_lin_syst}, which we terminate if the relative size of the residual reaches $10^{-4}$ or if the number of iterations reaches 150. 
To minimize the GCV functional, we use the MATLAB built-in function \texttt{fminbnd}, for which we set the lower bound to $0$ and the upper bound to $10^{-1}$, $10^{-2}$, $10^{-4}$, depending on the maximum intensity of the image. The tolerance \texttt{TolX} was set to $10^{-8}$. 

\begin{figure}[!ht]
    \centering
    \begin{subfigure}[b]{\textwidth}
    \includegraphics[width=.28\textwidth]{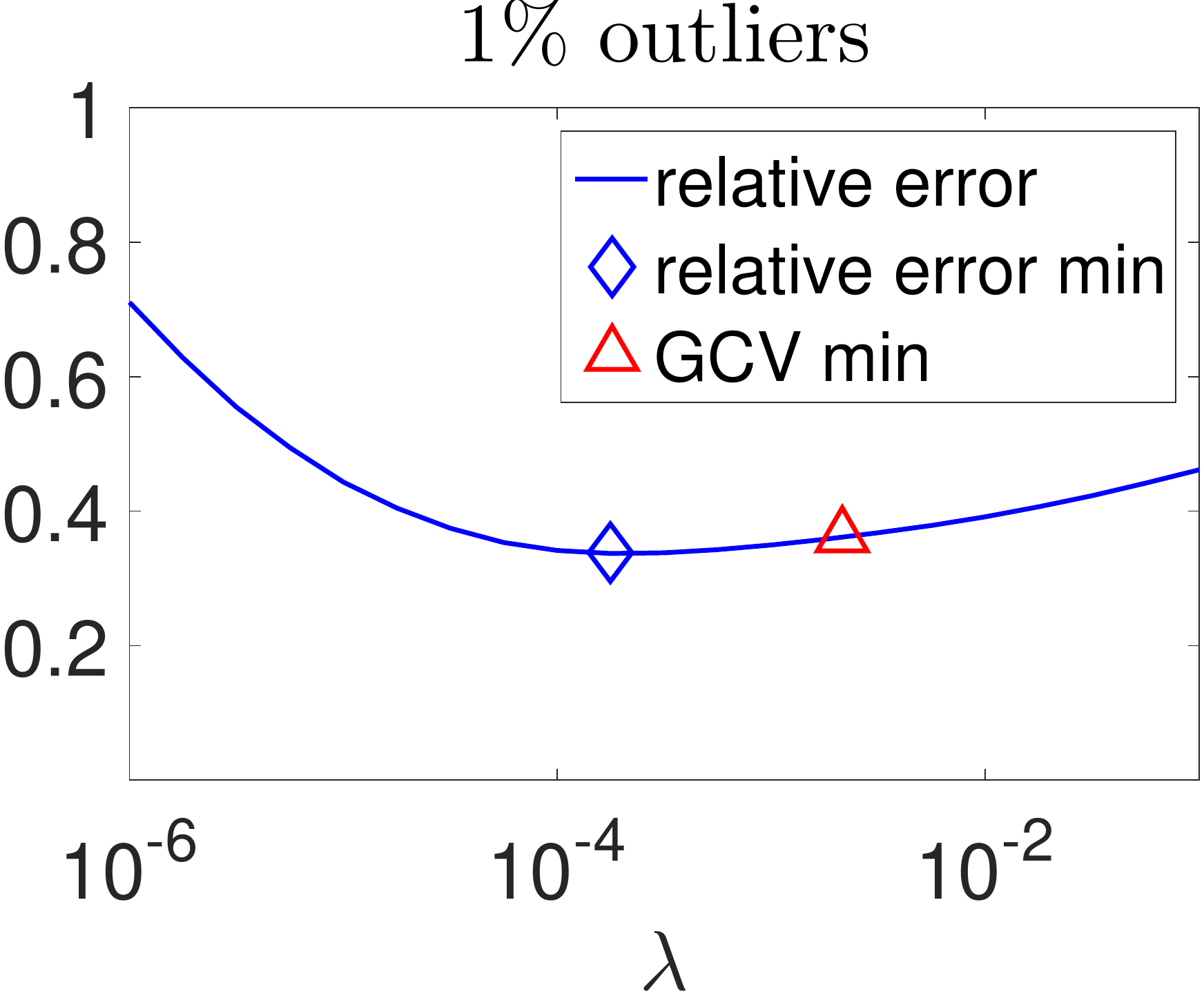} 
    \quad
	\includegraphics[width=.28\textwidth]{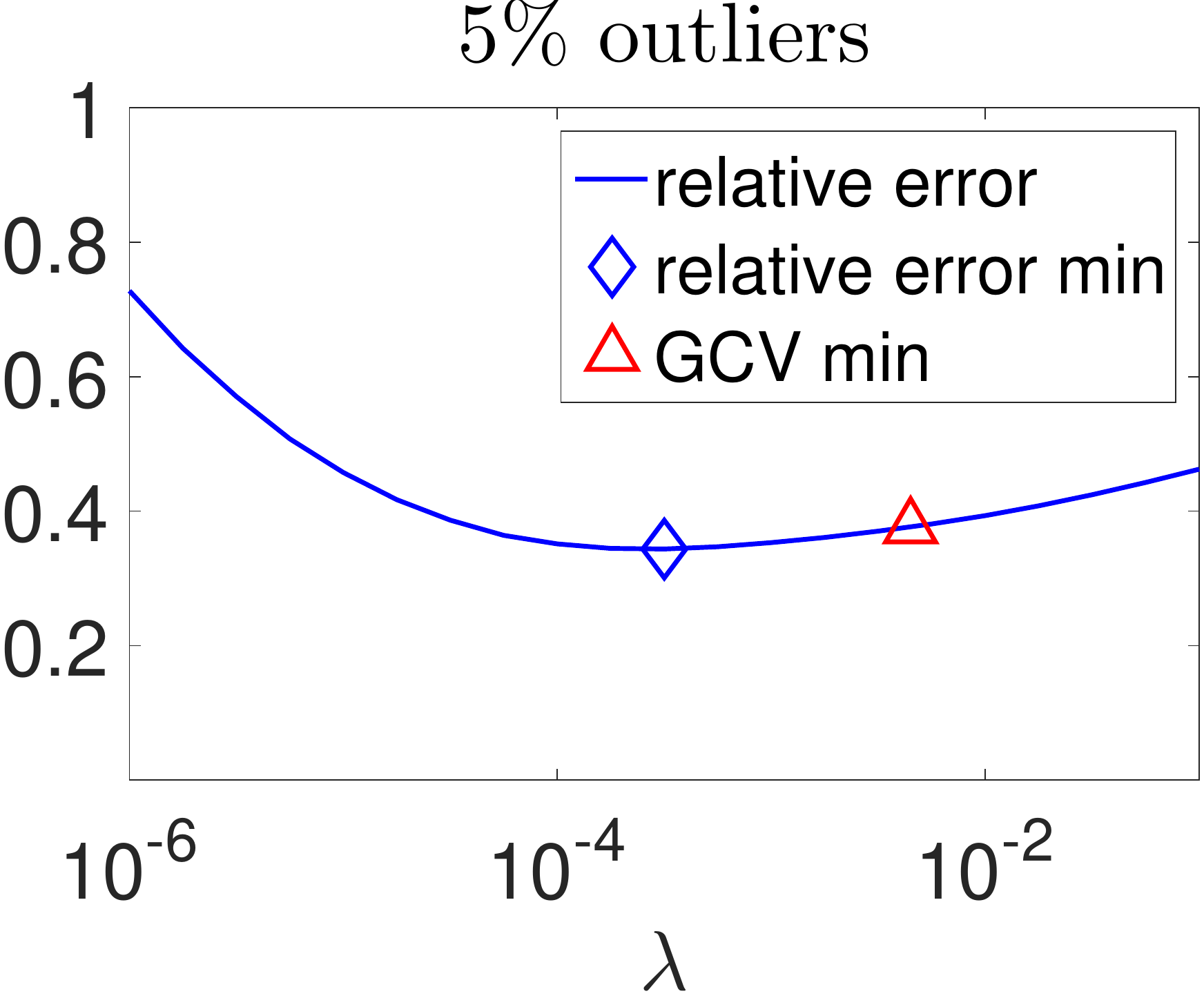} 
    \quad
	\includegraphics[width=.28\textwidth]{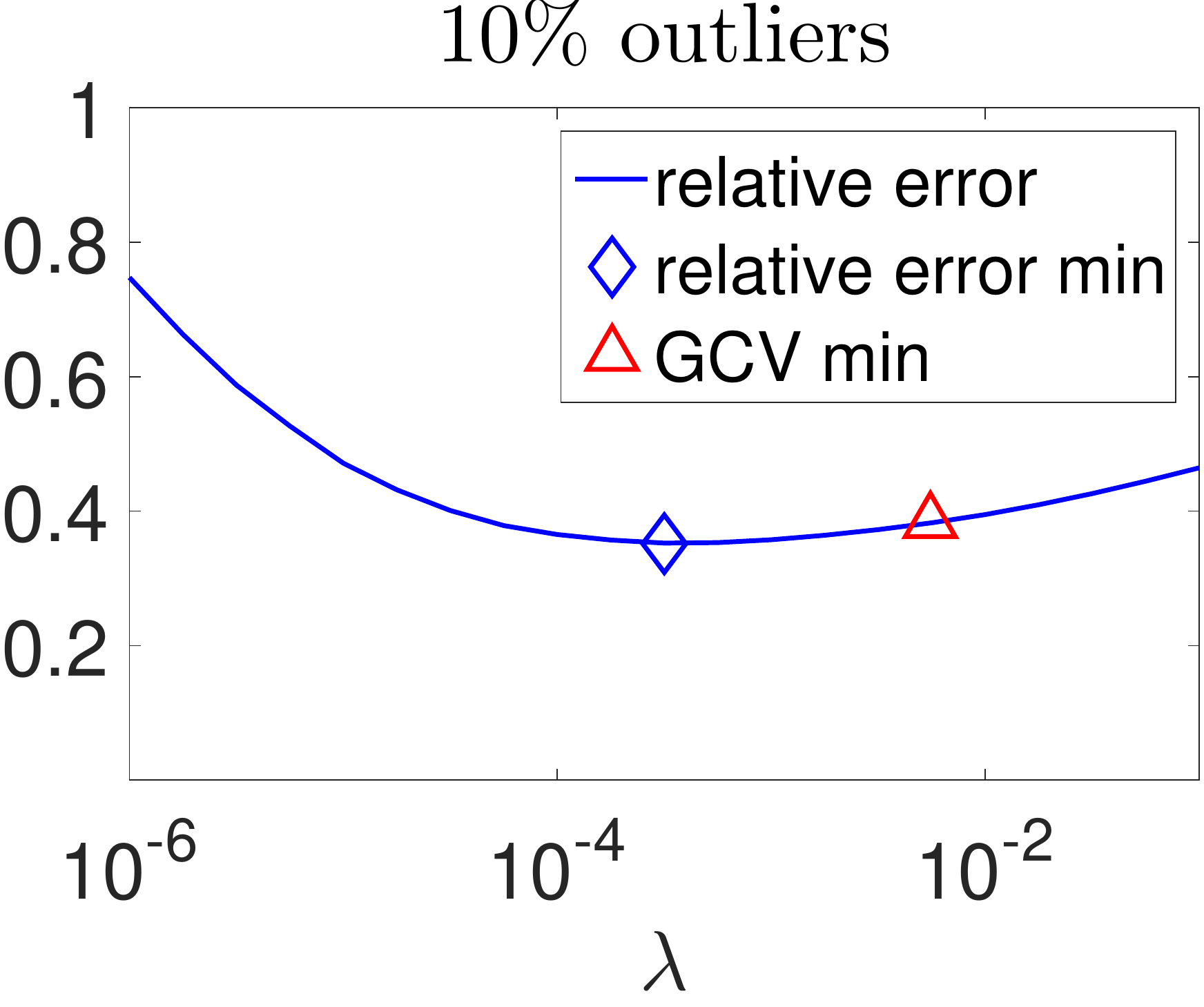} 
       	\caption{Satellite single-frame, max. intensity 255 (SNR = 5).}\end{subfigure}
    
\vspace*{.3cm}

    \begin{subfigure}[b]{\textwidth}
    \includegraphics[width=.28\textwidth]{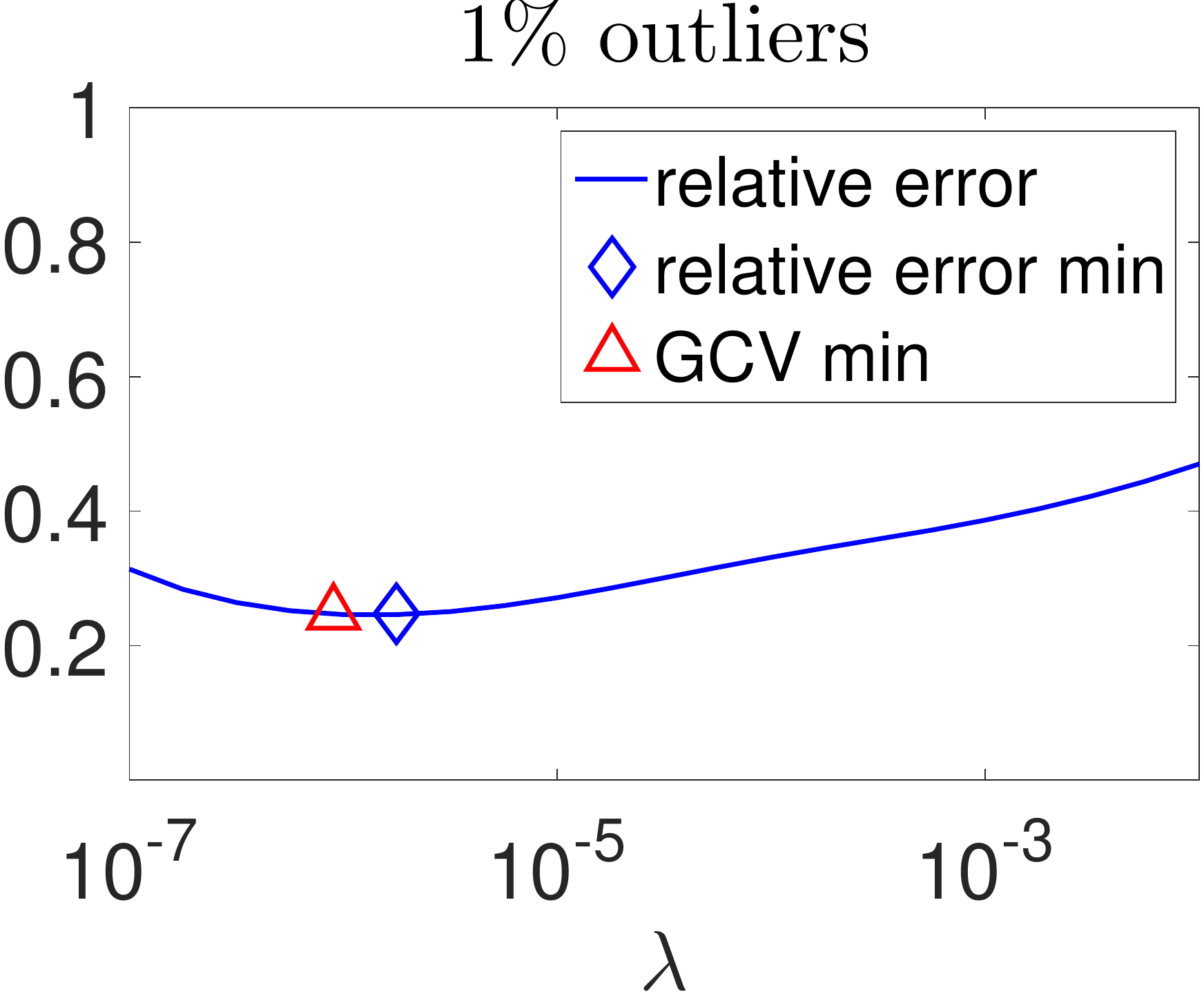} 
    \quad
	\includegraphics[width=.28\textwidth]{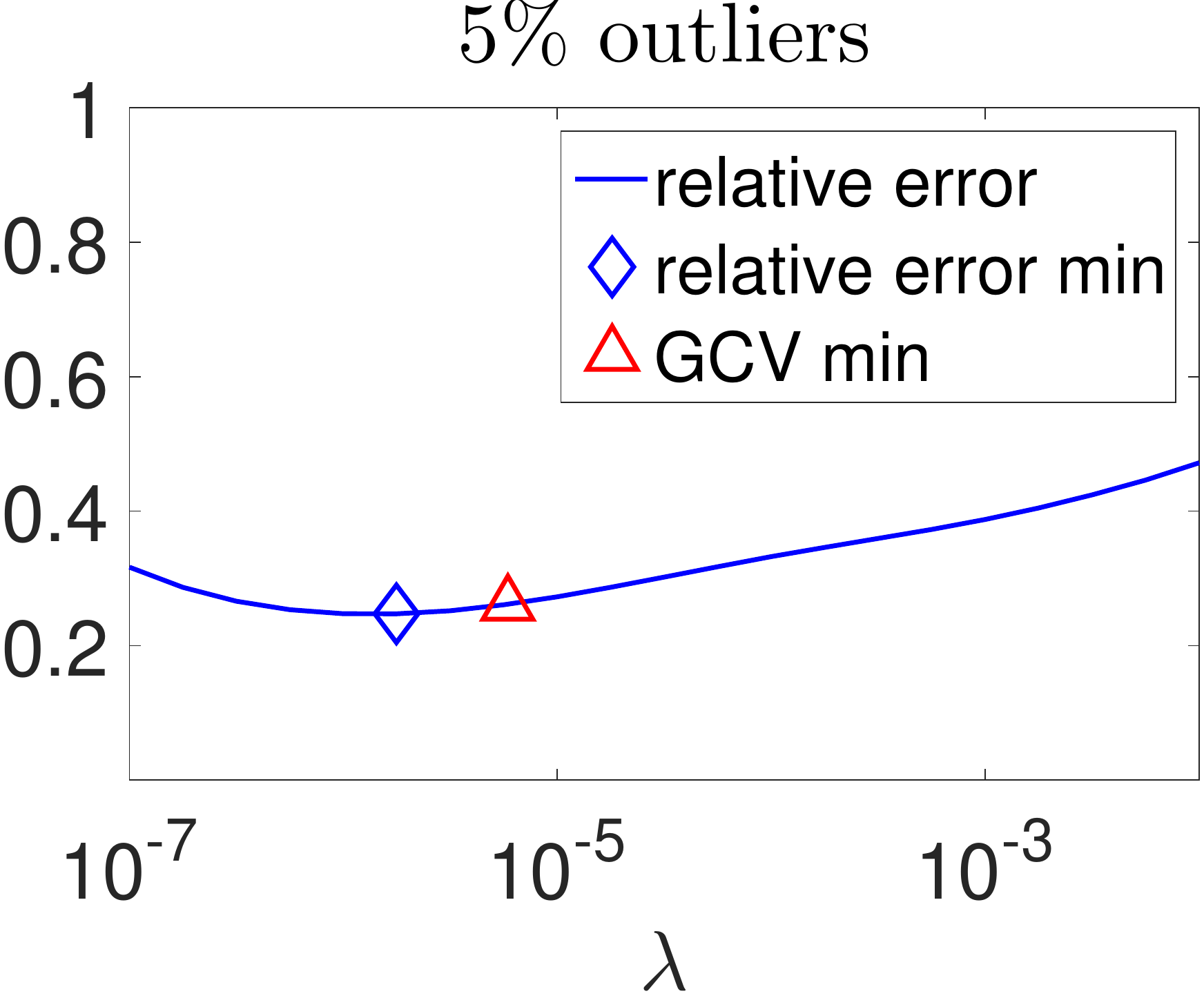}
    \quad
	\includegraphics[width=.28\textwidth]{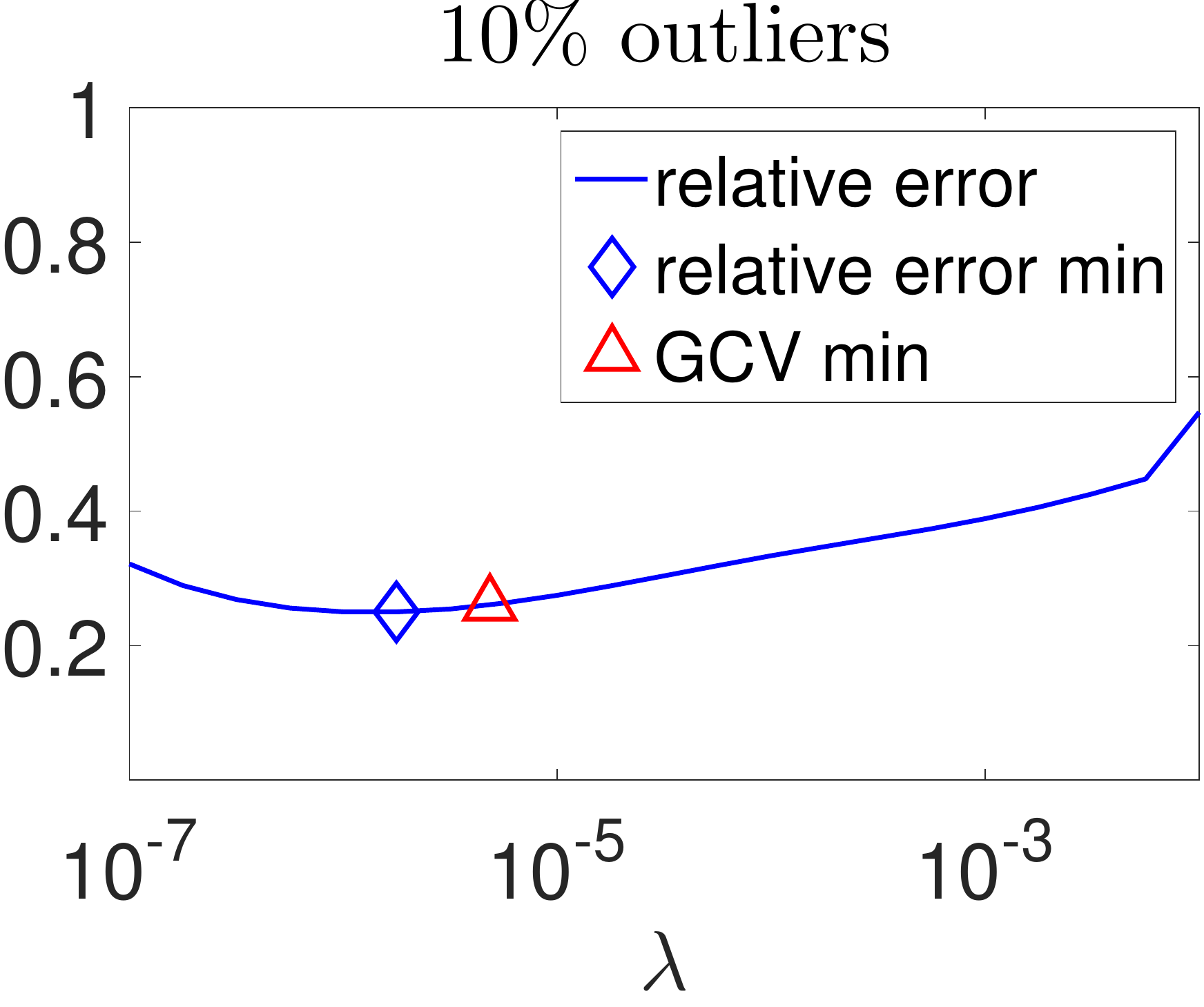} 
       	\caption{Satellite single-frame, max. intensity 2550 (SNR = 17).}\end{subfigure}
    
\vspace*{.3cm}

    \begin{subfigure}[b]{\textwidth}
    \includegraphics[width=.28\textwidth]{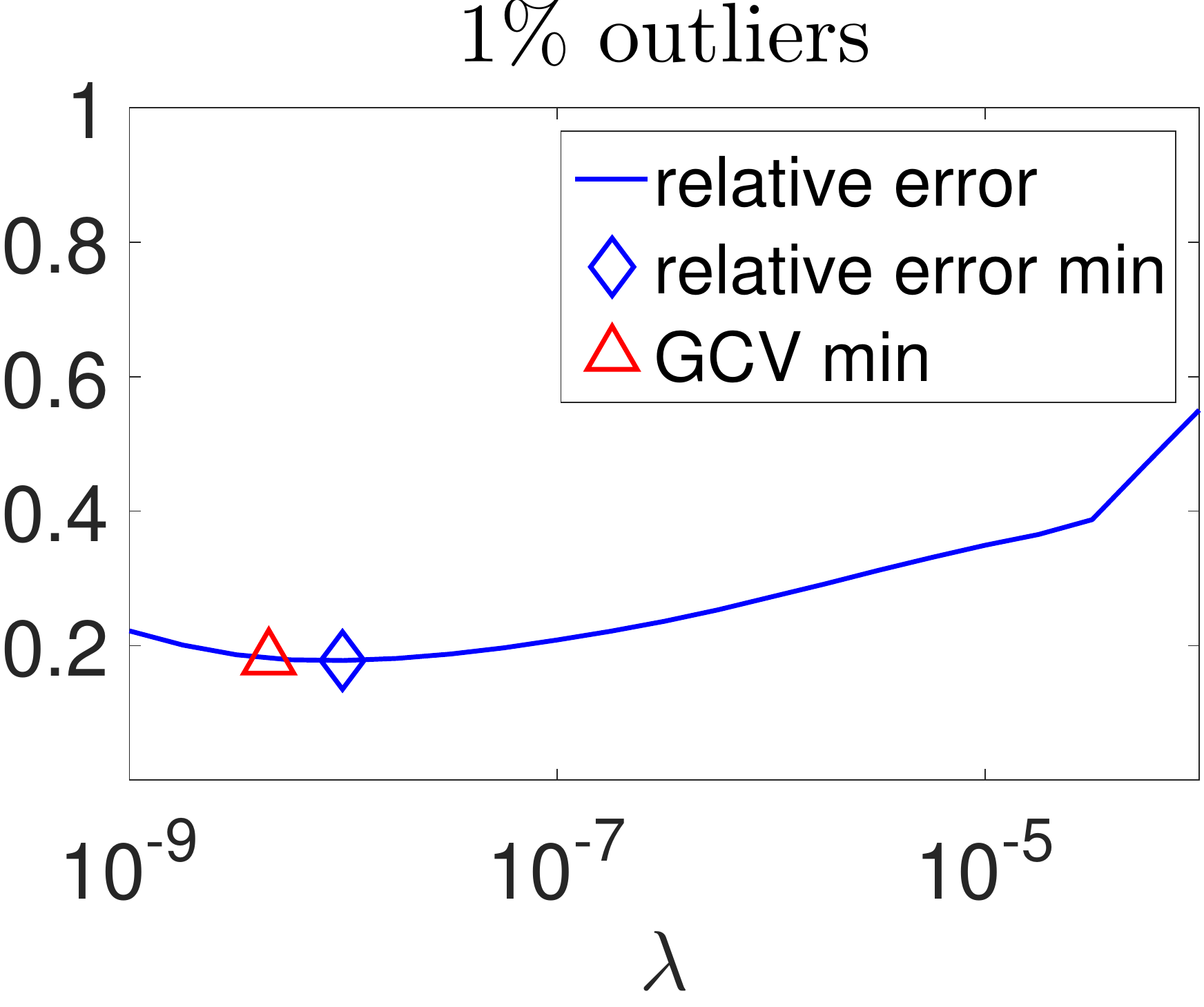} 
    \quad
	\includegraphics[width=.28\textwidth]{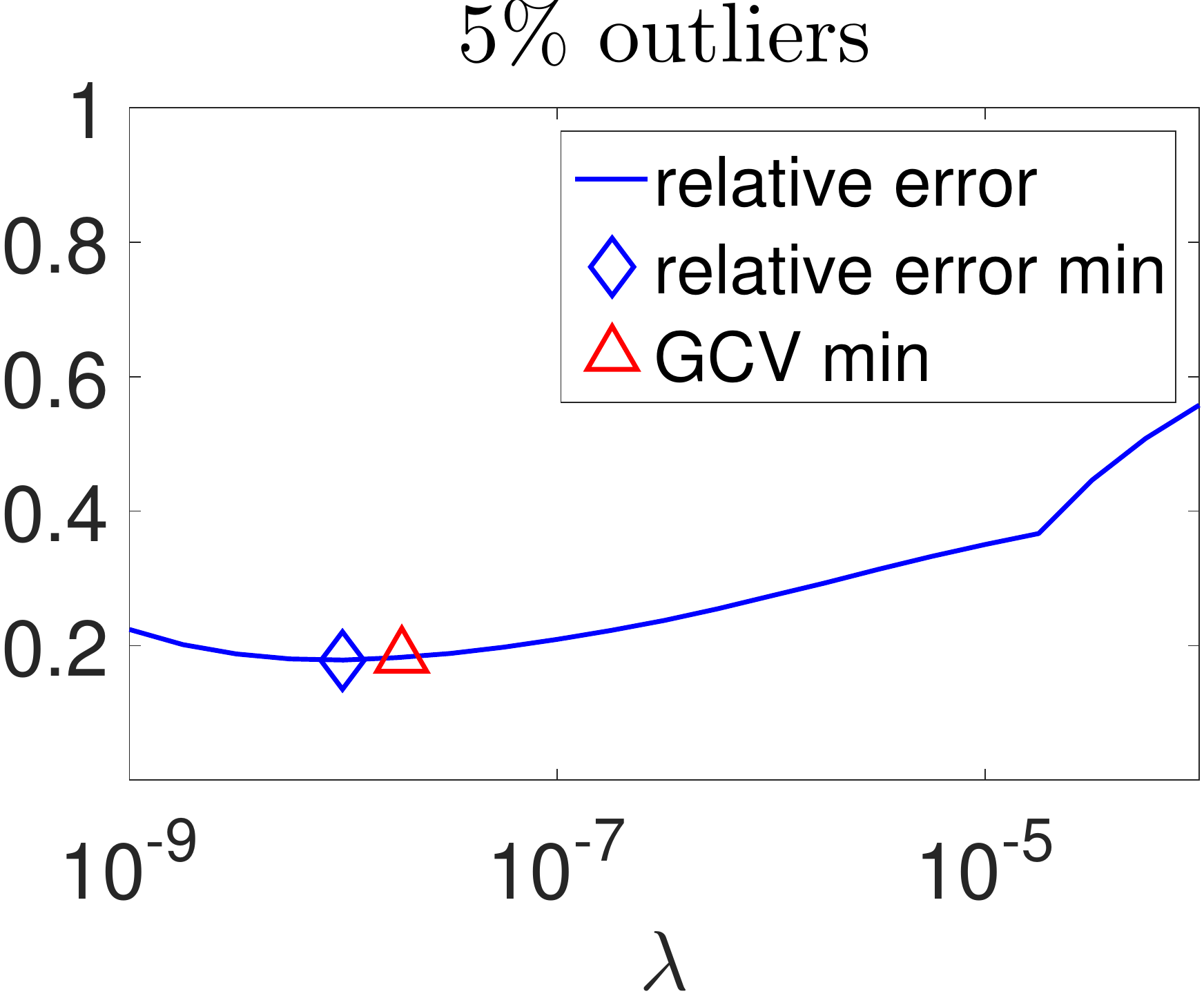}
    \quad
	\includegraphics[width=.28\textwidth]{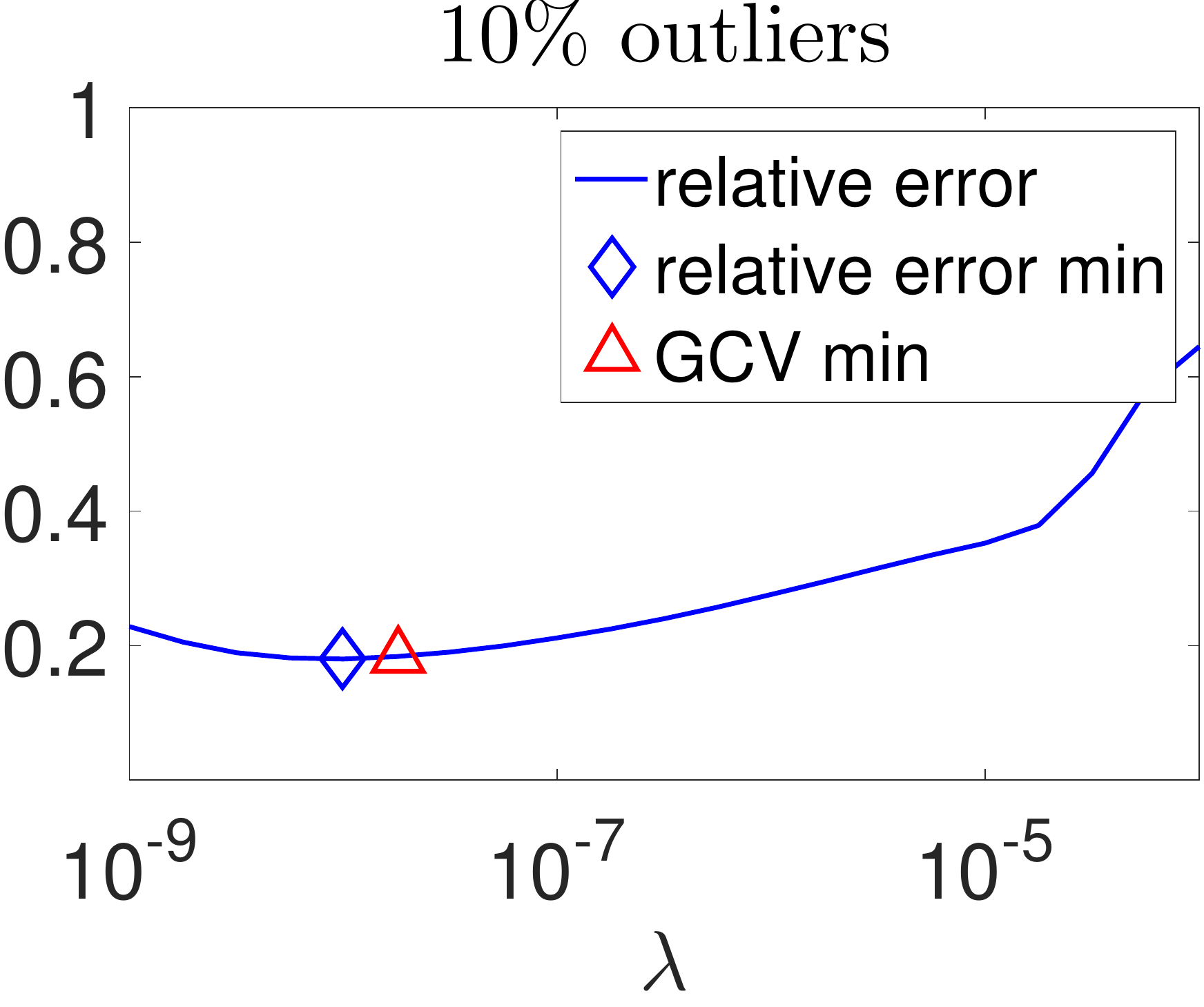} 
       	\caption{Satellite single-frame, max. intensity 25500 (SNR = 52).}\end{subfigure}
    
\vspace*{.3cm}

    \begin{subfigure}[b]{\textwidth}
    \includegraphics[width=.28\textwidth]{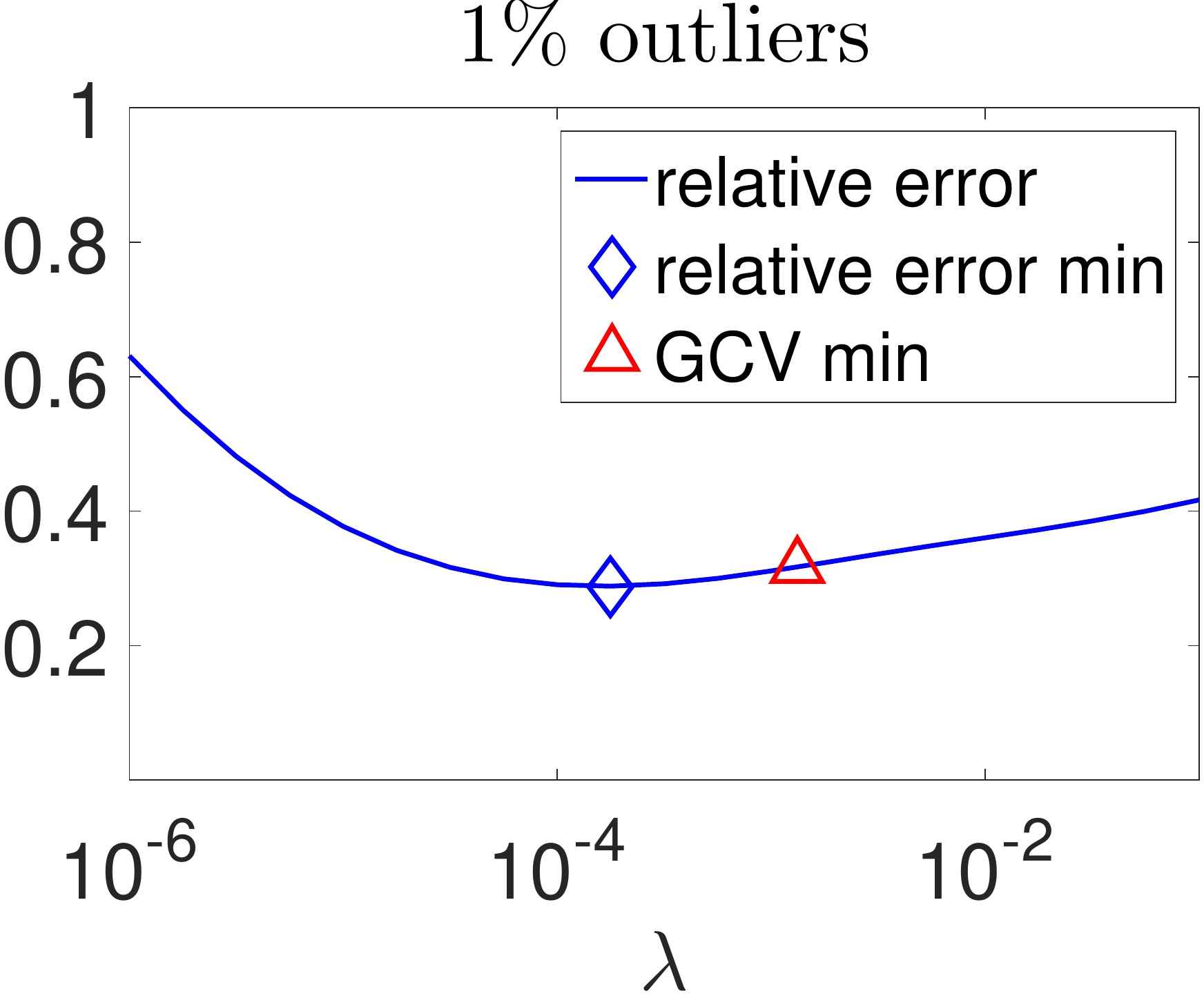} 
    \quad
	\includegraphics[width=.28\textwidth]{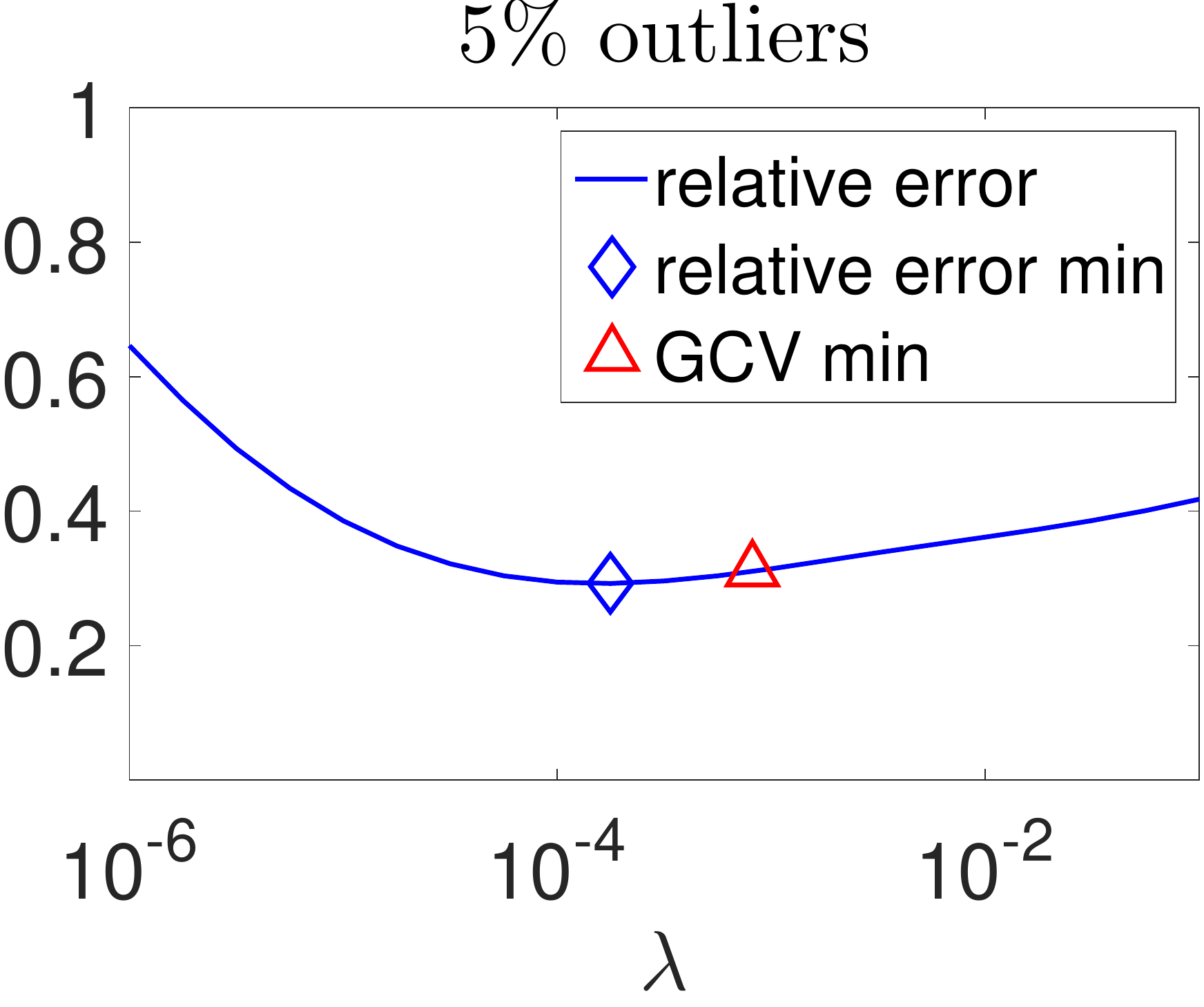}
    \quad
	\includegraphics[width=.28\textwidth]{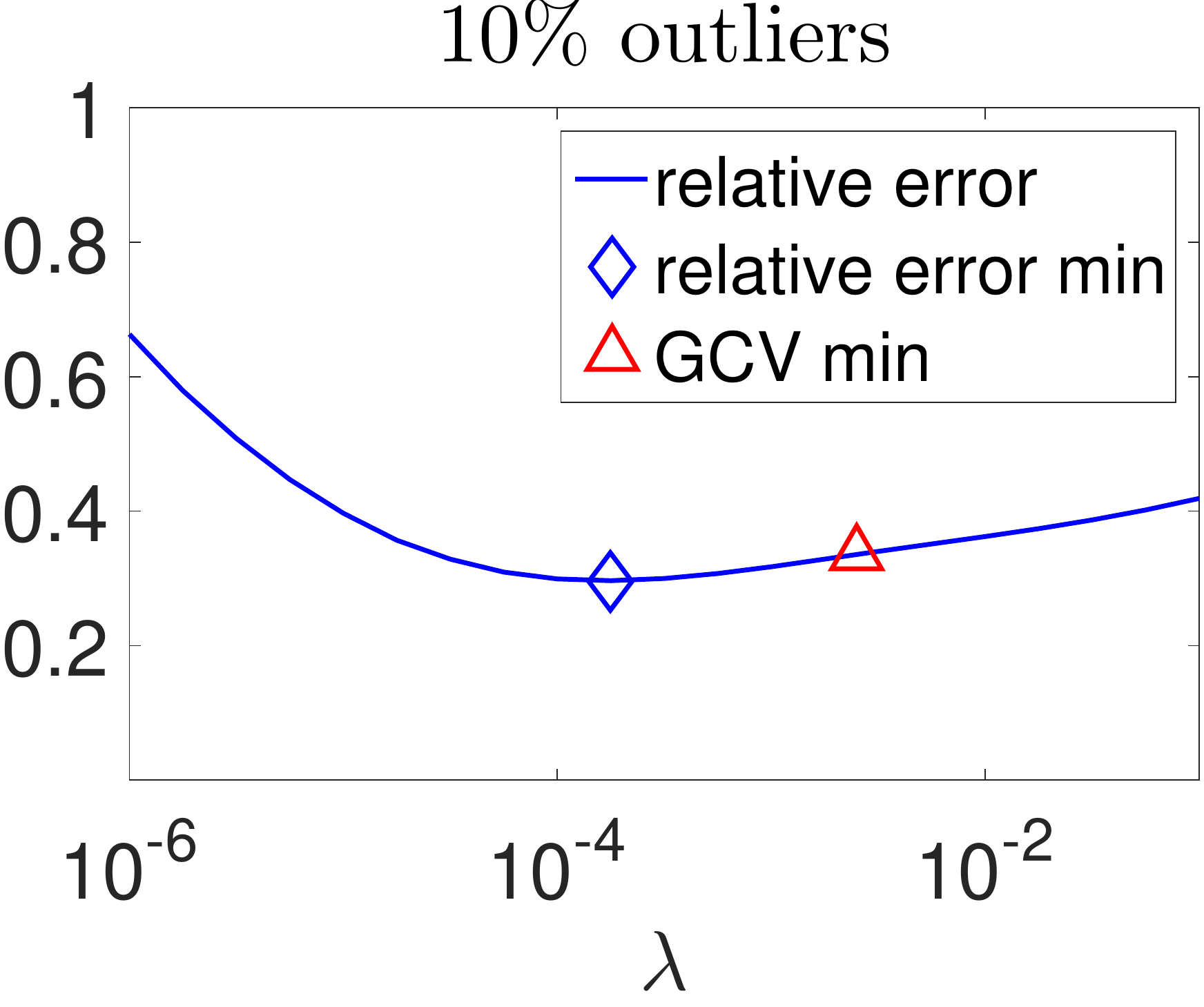} 
       	\caption{Satellite multi-frame, max. intensity 255 (SNR = 5).}\end{subfigure}
    
    	\caption{GCV for data with outliers. 
} \label{fig:GCV_satellite}
\end{figure}

For test problem Satellite, we show the semiconvergence curves including the minimum error and the error obtained using GCV in Figure~\ref{fig:GCV_satellite}. Quantitative results (averaged over 10 realizations of outliers) for both test problems are shown in Table~\ref{tab:gcv}. We observe that the proposed rule is rather stable with respect to the increasing number of outliers and generally better for the Carbon ash than for the Satellite. As expected, the method provides better approximation of the optimal regularization parameter for smaller noise levels  (larger $Ax_\text{true}$), where the functional \eqref{eq:wls} approximates better the maximum likelihood functional for the mixed Poisson--Gaussian model. Occasionally, GCV provides slightly worse reconstruction for the highest percentage (10\%) of outliers.

\begin{landscape}
\begin{table}
\centering
\caption{Relative errors of the reconstruction: optimal $\lambda$ vs. $\lambda$ obtained by minimizing the GCV functional \eqref{eq:gcv_fun}. 
}\label{tab:gcv}
\begin{subtable}[h]{1.5\textwidth}
\subcaption{Satellite}
\centering
{\scriptsize\begin{tabular}{llcllclcclc} \toprule 
\multicolumn{2}{r}{}& \multicolumn{4}{r}{error min} & \multicolumn{4}{l}{($\lambda$ optimal)} \\ 
\multicolumn{2}{r}{} & \multicolumn{4}{r}{error GCV} & \multicolumn{4}{l}{($\lambda$ GCV)} \\ 
\multicolumn{1}{r}{} & \multicolumn{9}{c}{\% outliers}\\ 
problem & max. int. & \multicolumn{2}{c}{1\%} & & \multicolumn{2}{c}{5\%} & & \multicolumn{2}{c}{10\%} \\ \midrule 
\multirow{2}{*}{single-frame} & \multirow{2}{*}{255}&  $3.39\times 10^{-1}$ & ($\lambda_\text{opt} = 2.1\times 10^{-4}$) & & $3.43\times 10^{-1}$ & ($\lambda_\text{opt} = 2.5\times 10^{-4}$) & & $3.49\times 10^{-1}$ & ($\lambda_\text{opt} = 3.2\times 10^{-4}$)\\ 
 & &  $3.60\times 10^{-1}$ & ($\lambda_\text{GCV} = 2.1\times 10^{-3}$) & & $3.68\times 10^{-1}$ & ($\lambda_\text{GCV} = 3.1\times 10^{-3}$) & & $3.80\times 10^{-1}$ & ($\lambda_\text{GCV} = 5.6\times 10^{-3}$)\\ \addlinespace[.2cm] 
\multirow{2}{*}{single-frame} & \multirow{2}{*}{2550}&  $2.46\times 10^{-1}$ & ($\lambda_\text{opt} = 1.5\times 10^{-6}$) & & $2.48\times 10^{-1}$ & ($\lambda_\text{opt} = 1.5\times 10^{-6}$) & & $2.50\times 10^{-1}$ & ($\lambda_\text{opt} = 1.7\times 10^{-6}$)\\ 
 & &  $2.48\times 10^{-1}$ & ($\lambda_\text{GCV} = 1.5\times 10^{-6}$) & & $2.57\times 10^{-1}$ & ($\lambda_\text{GCV} = 4.5\times 10^{-6}$) & & $2.69\times 10^{-1}$ & ($\lambda_\text{GCV} = 8.4\times 10^{-6}$)\\ \addlinespace[.2cm] 
\multirow{2}{*}{single-frame} & \multirow{2}{*}{25500}&  $1.78\times 10^{-1}$ & ($\lambda_\text{opt} = 1.0\times 10^{-8}$) & & $1.79\times 10^{-1}$ & ($\lambda_\text{opt} = 1.0\times 10^{-8}$) & & $1.81\times 10^{-1}$ & ($\lambda_\text{opt} = 1.0\times 10^{-8}$)\\ 
 & &  $1.79\times 10^{-1}$ & ($\lambda_\text{GCV} = 8.8\times 10^{-9}$) & & $1.81\times 10^{-1}$ & ($\lambda_\text{GCV} = 1.5\times 10^{-8}$) & & $2.42\times 10^{-1}$ & ($\lambda_\text{GCV} = 7.5\times 10^{-6}$)\\ \addlinespace[.2cm] 
\multirow{2}{*}{mutli-frame} & \multirow{2}{*}{255}& $2.89\times 10^{-1}$ & ($\lambda_\text{opt} = 1.8\times 10^{-4}$) & & $2.92\times 10^{-1}$ & ($\lambda_\text{opt} = 1.8\times 10^{-4}$) & & $2.97\times 10^{-1}$ & ($\lambda_\text{opt} = 1.8\times 10^{-4}$)\\ 
 & & $3.16\times 10^{-1}$ & ($\lambda_\text{GCV} = 1.6\times 10^{-3}$) & & $3.13\times 10^{-1}$ & ($\lambda_\text{GCV} = 1.2\times 10^{-3}$) & & $3.49\times 10^{-1}$ & ($\lambda_\text{GCV} = 5.8\times 10^{-3}$) \\ \addlinespace[.2cm] 
\bottomrule\end{tabular}
}
\end{subtable}

\begin{subtable}[h]{1.5\textwidth}
\subcaption{Carbon ash}
\centering
{\scriptsize\begin{tabular}{llcllclcclc} \toprule 
\multicolumn{2}{r}{}& \multicolumn{4}{r}{error min} & \multicolumn{4}{l}{($\lambda$ optimal)} \\ 
\multicolumn{2}{r}{} & \multicolumn{4}{r}{error GCV} & \multicolumn{4}{l}{($\lambda$ GCV)} \\ 
\multicolumn{1}{r}{} & \multicolumn{9}{c}{\% outliers}\\ 
problem & max. int. & \multicolumn{2}{c}{1\%} & & \multicolumn{2}{c}{5\%} & & \multicolumn{2}{c}{10\%} \\ \midrule 
\multirow{2}{*}{single-frame} & \multirow{2}{*}{255}& $3.08\times 10^{-1}$ & ($\lambda_\text{opt} = 1.8\times 10^{-3}$) & & $3.09\times 10^{-1}$ & ($\lambda_\text{opt} = 1.8\times 10^{-3}$) & & $3.10\times 10^{-1}$ & ($\lambda_\text{opt} = 1.7\times 10^{-3}$)\\ 
& & $3.12\times 10^{-1}$ & ($\lambda_\text{GCV} = 8.9\times 10^{-4}$) & & $3.11\times 10^{-1}$ & ($\lambda_\text{GCV} = 1.5\times 10^{-3}$) & & $3.11\times 10^{-1}$ & ($\lambda_\text{GCV} = 2.2\times 10^{-3}$) \\ \addlinespace[.2cm] 
\multirow{2}{*}{single-frame} & \multirow{2}{*}{2550}& $2.91\times 10^{-1}$ & ($\lambda_\text{opt} = 1.9\times 10^{-5}$) & & $2.92\times 10^{-1}$ & ($\lambda_\text{opt} = 2.1\times 10^{-5}$) & & $2.92\times 10^{-1}$ & ($\lambda_\text{opt} = 1.9\times 10^{-5}$)\\ 
& & $2.97\times 10^{-1}$ & ($\lambda_\text{GCV} = 9.7\times 10^{-6}$) & & $2.95\times 10^{-1}$ & ($\lambda_\text{GCV} = 1.2\times 10^{-5}$) & & $2.93\times 10^{-1}$ & ($\lambda_\text{GCV} = 2.2\times 10^{-5}$)\\ \addlinespace[.2cm] 
\multirow{2}{*}{single-frame} & \multirow{2}{*}{25500}& $2.77\times 10^{-1}$ & ($\lambda_\text{opt} = 3.0\times 10^{-7}$) & & $2.78\times 10^{-1}$ & ($\lambda_\text{opt} = 2.7\times 10^{-7}$) & & $2.79\times 10^{-1}$ & ($\lambda_\text{opt} = 2.6\times 10^{-7}$)\\ 
& & $3.00\times 10^{-1}$ & ($\lambda_\text{GCV} = 5.2\times 10^{-8}$) & & $2.84\times 10^{-1}$ & ($\lambda_\text{GCV} = 9.4\times 10^{-8}$) & & $2.82\times 10^{-1}$ & ($\lambda_\text{GCV} = 1.3\times 10^{-7}$)\\ \addlinespace[.2cm] 
\multirow{2}{*}{multi-frame} & \multirow{2}{*}{255}& $3.03\times 10^{-1}$ & ($\lambda_\text{opt} = 1.8\times 10^{-3}$) & & $3.04\times 10^{-1}$ & ($\lambda_\text{opt} = 1.8\times 10^{-3}$) & & $3.04\times 10^{-1}$ & ($\lambda_\text{opt} = 1.8\times 10^{-3}$)\\ 
 & & $3.14\times 10^{-1}$ & ($\lambda_\text{GCV} = 6.8\times 10^{-4}$) & & $3.07\times 10^{-1}$ & ($\lambda_\text{GCV} = 9.8\times 10^{-4}$) & & $3.05\times 10^{-1}$ & ($\lambda_\text{GCV} = 1.9\times 10^{-3}$)\\ \addlinespace[.2cm] 
\bottomrule\end{tabular}
}
\end{subtable}
\end{table} 
\end{landscape}

\subsection{Linear subproblems}

\begin{figure}[!ht]
    \centering
        \begin{subfigure}[b]{0.85\textwidth}
        	\includegraphics[width=.47\textwidth]{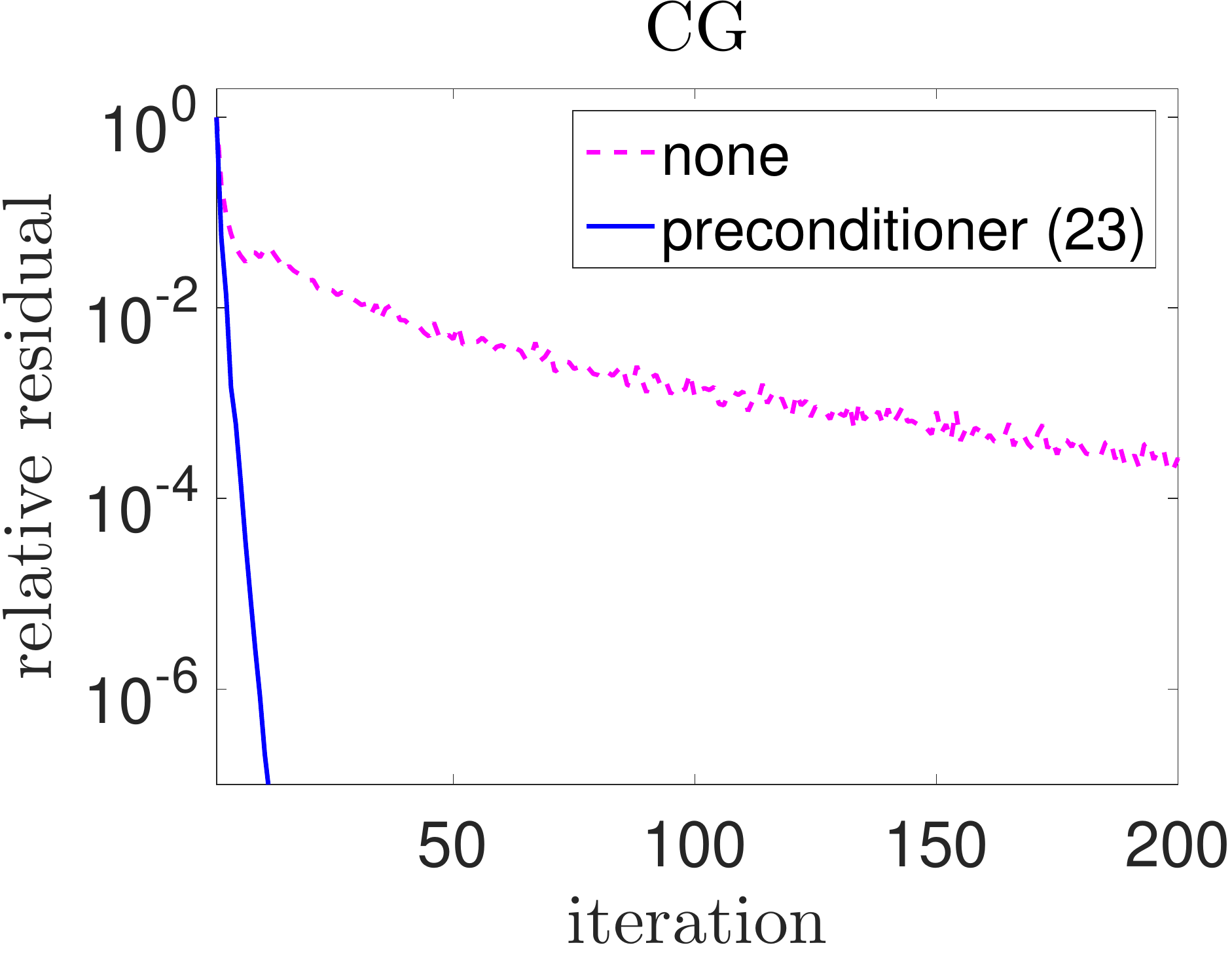} 
        	\hspace*{.1cm}
        	\includegraphics[width=.47\textwidth]{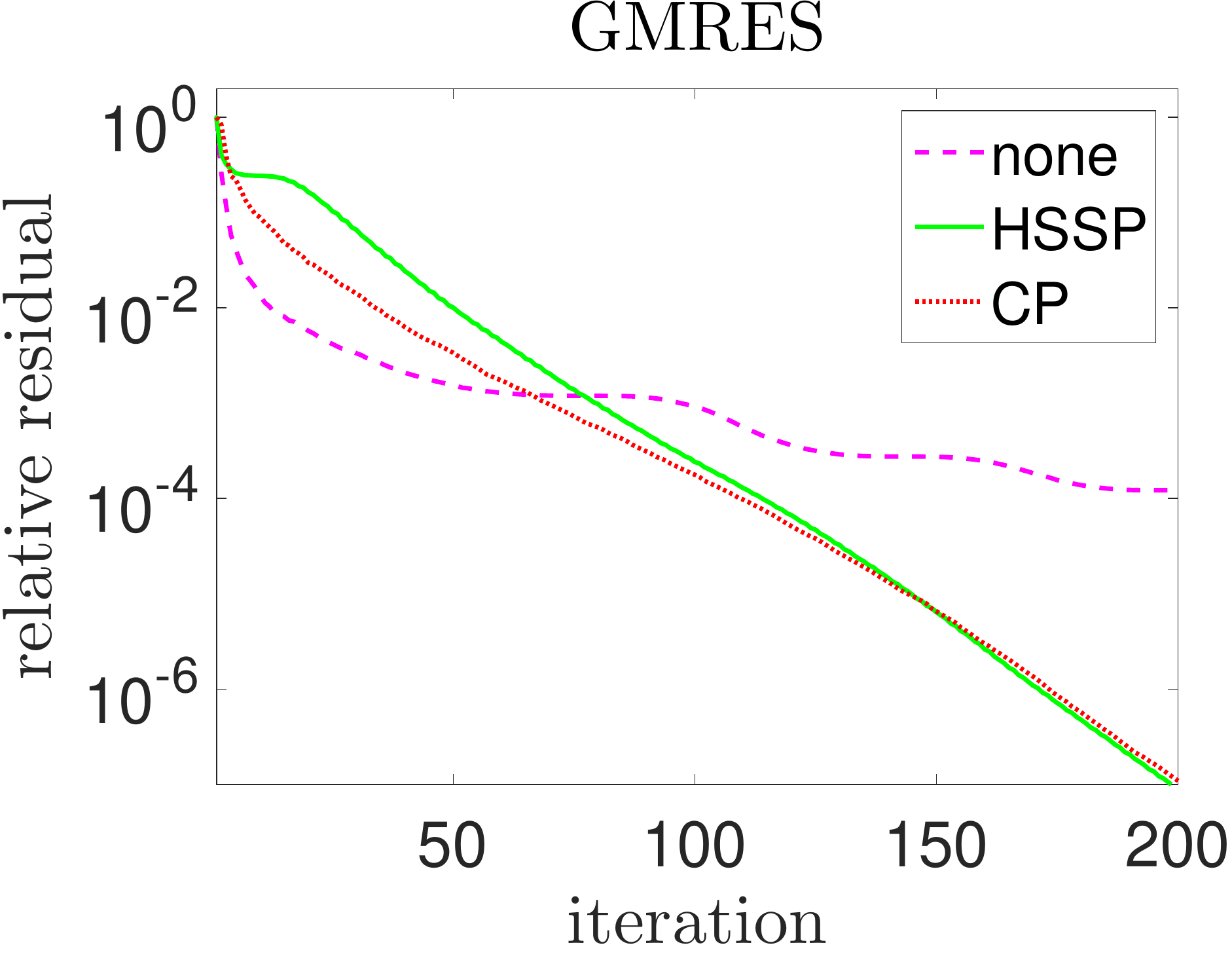} 
        	\caption{Satellite single-frame, random diagonal}
    	\end{subfigure}
    	\caption{Preconditoner defined in \eqref{eq:prec}, constraint preconditioner (CP), and Hermitian and skew-Hermitian splitting preconditioner (HSSP) performance for $(A^TDA + \lambda L^TL)s = -A^Tb$, where $A$ and $b$ are adopted from the test problem Satellite, and $D$ is a diagonal with random entries uniformly distributed in $(0,1)$. 
}\label{fig:prec1}
\end{figure}

\begin{figure}[!ht]
    \centering
        \begin{subfigure}[b]{0.85\textwidth}
        	\includegraphics[width=.47\textwidth]{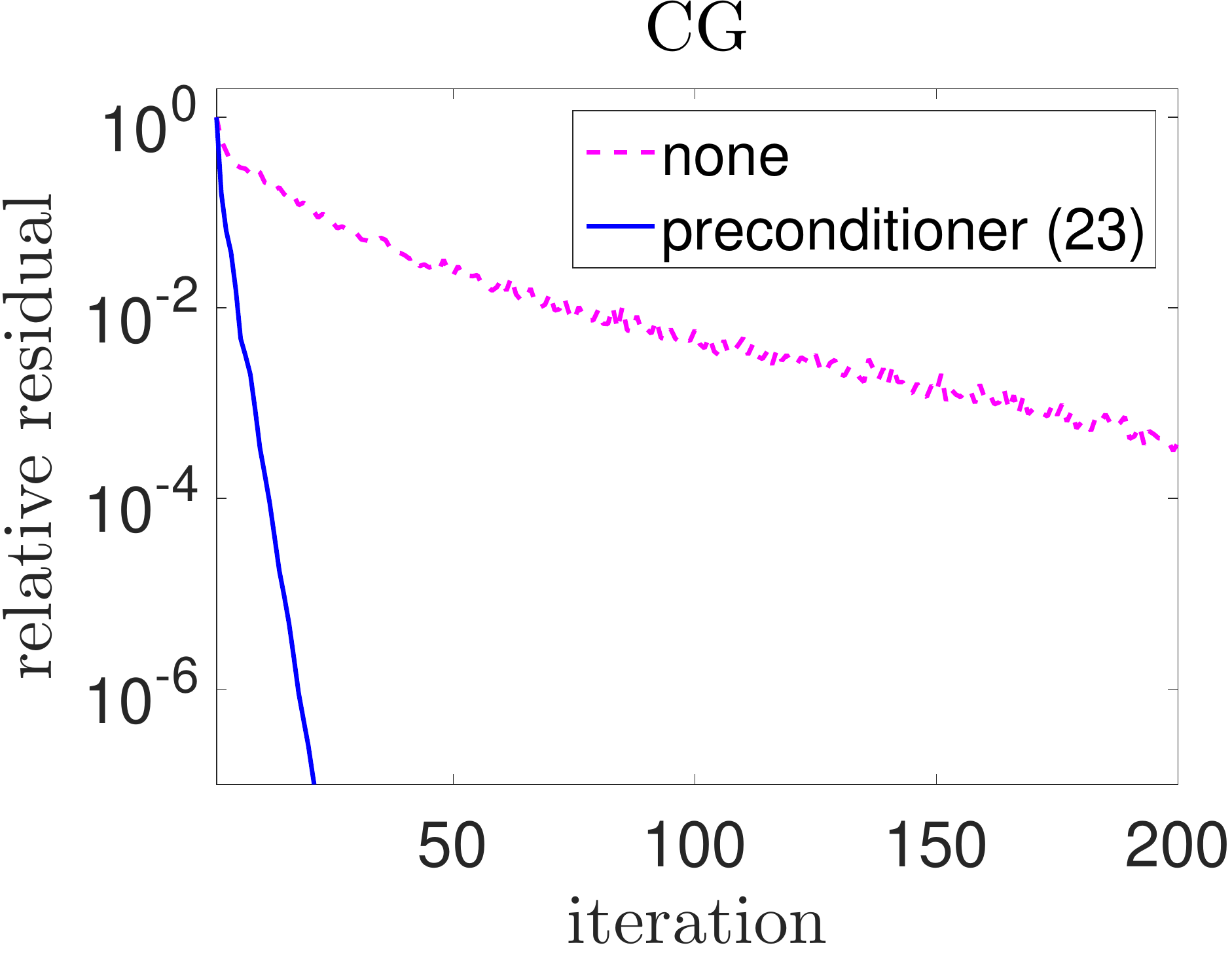} 
        	\hspace*{.1cm}
        	\includegraphics[width=.47\textwidth]{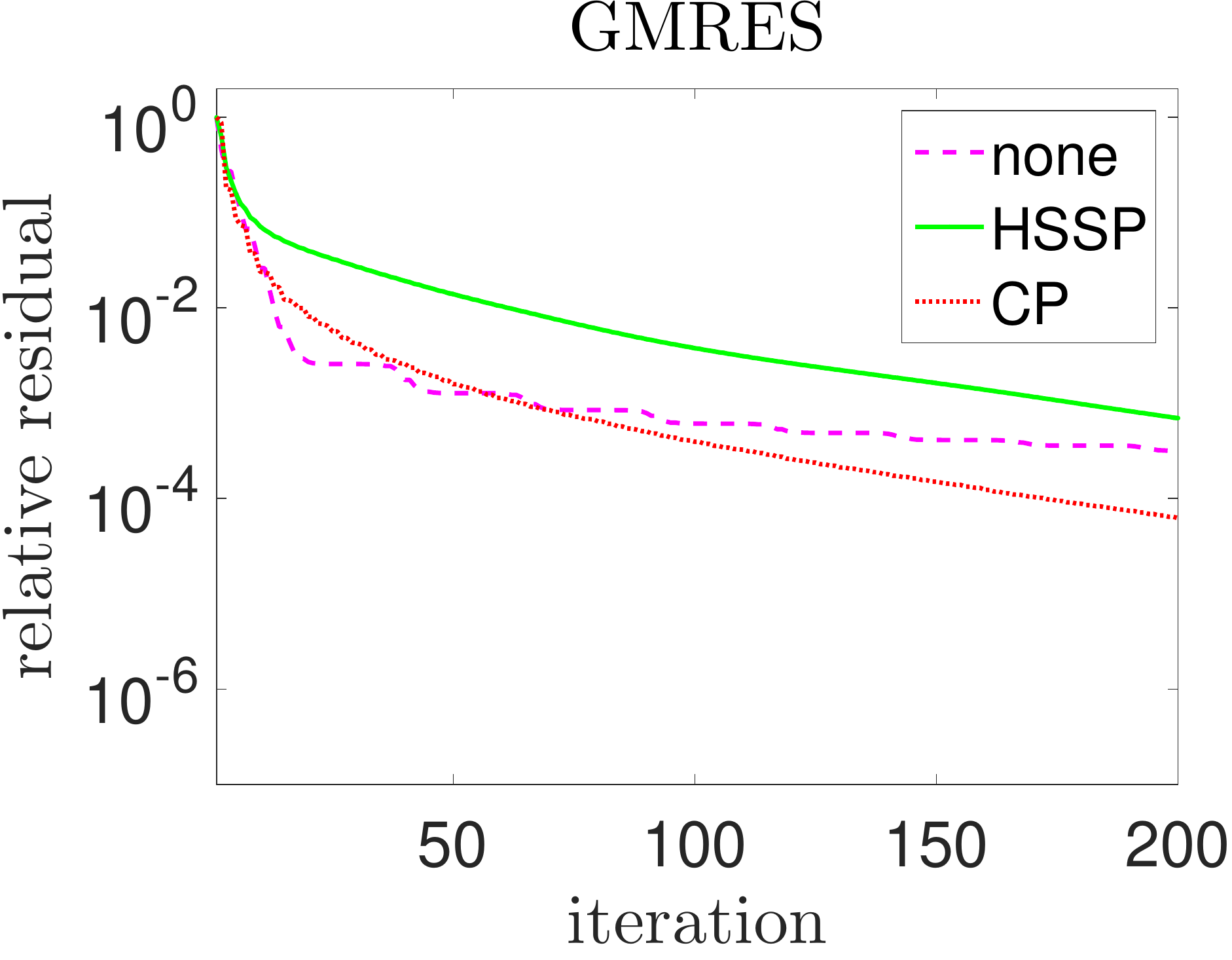} 
        	\caption{Satellite single-frame, Newton it = 3}
    	\end{subfigure}
    	\caption{Preconditoner defined in \eqref{eq:prec}, constraint preconditioner (CP), and Hermitian/skew-Hermitian preconditioner (HSSP) performance for $(A^TD^{(k)}A + \lambda L^TL)s = -\left(A^Tz^{(k)} + \lambda L^TLx^{(k)}\right)$. 
}\label{fig:prec2} 
\end{figure}

As mentioned earlier, various types of preconditioners have been developed to speed up convergence of iterative methods applied to systems of type \eqref{eq:lin_system} or its saddle-point counterpart 
\[
\begin{pmatrix}
  D^{-1} & A  \\
  A^T & \lambda L^TL\\
 \end{pmatrix} = 
 \begin{pmatrix}
  -z \\
  \lambda L^TLx\\
 \end{pmatrix}.\label{eq:saddle_point}
\]
The Hermitian and skew-Hermitian (HSS) preconditioner, as well as constraint precondtioner belong to the best known preconditioners for this type of linear system. Both were incorporated in GMRES and tested on deblurring problems with random diagonal scaling $D$ in \cite{Benzi2006Preconditioned}. Using random $D$, they indeed accelerate convergence also in our case, as shown in Figure \ref{fig:prec1}. However, our preconditioner \eqref{eq:prec} provides a much better speedup. Moreover, for real computations, e.g., when the matrix $D$ is actually generated during the Projected Newton computation, the HSS and constraint preconditioners did not perform well, and even slowed down the convergence, see Figure~\ref{fig:prec2}. This is fortunately not the case for 
our proposed preconditioner. In this experiment, we did not assume projection on the non-negative half-plane and since in \eqref{eq:saddle_point}, we need to evaluate $D^{-1}$, if some component $D_{ii}=0$, we replaced it by $2\sqrt{\epsilon_\text{mach}}$, see also \cite{Mastronardi2008Fast}.
We also did not incorporate any outliers for these initial experiments with the preconditioners; 
these results are intended to show that 
our proposed preconditioning for these problems
often performs much better than the well-known standard preconditioners.
In fact, we see that the behavior of the constraint and HSS preconditioner depends heavily on the actual setting of the problem. In the remainder of this section we will therefore focus on the preconditioner given in \eqref{eq:prec}.

\begin{figure}[!ht]
    \centering
    \begin{subfigure}[b]{\textwidth}
    \includegraphics[width=.28\textwidth]{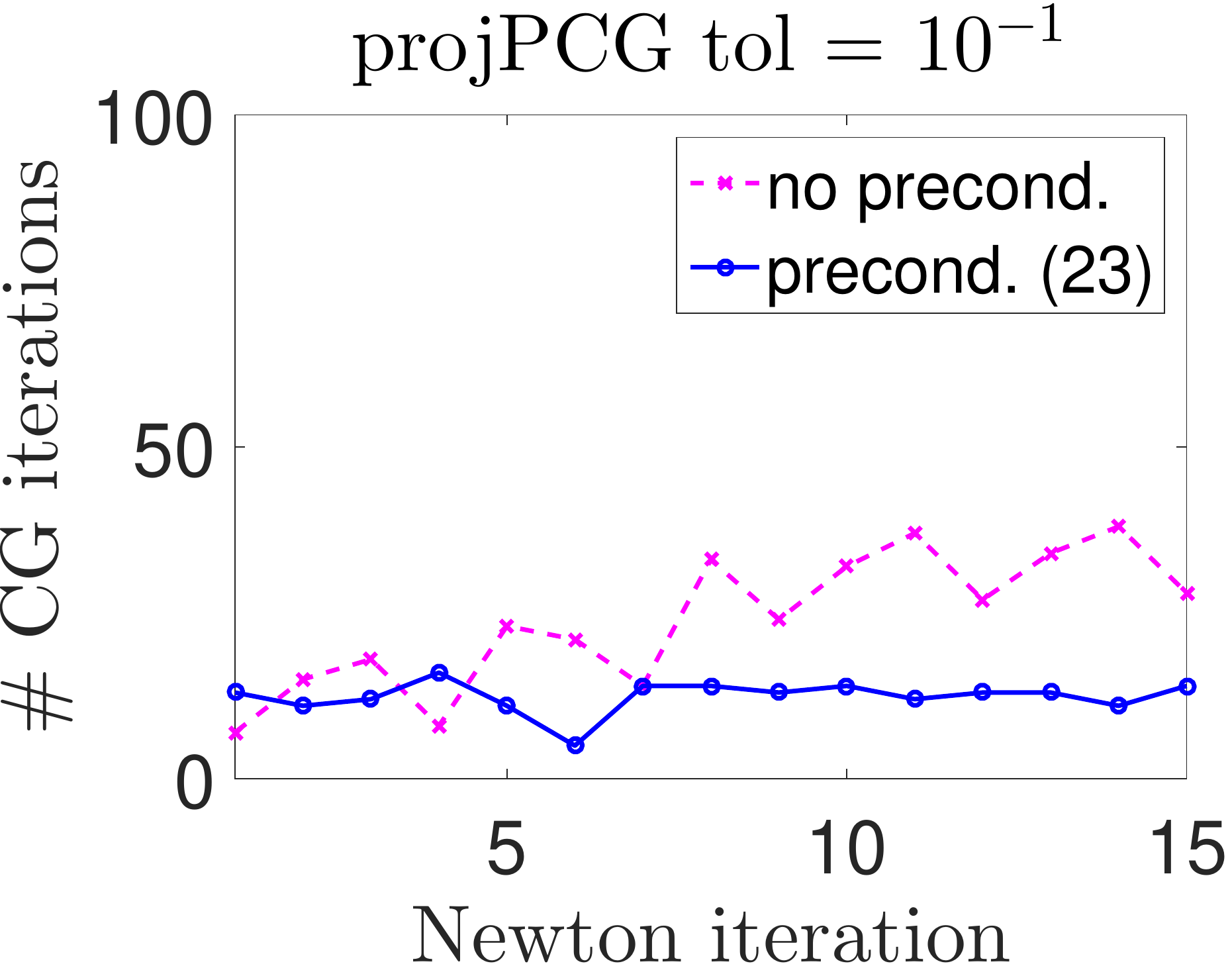} 
    \quad
	\includegraphics[width=.28\textwidth]{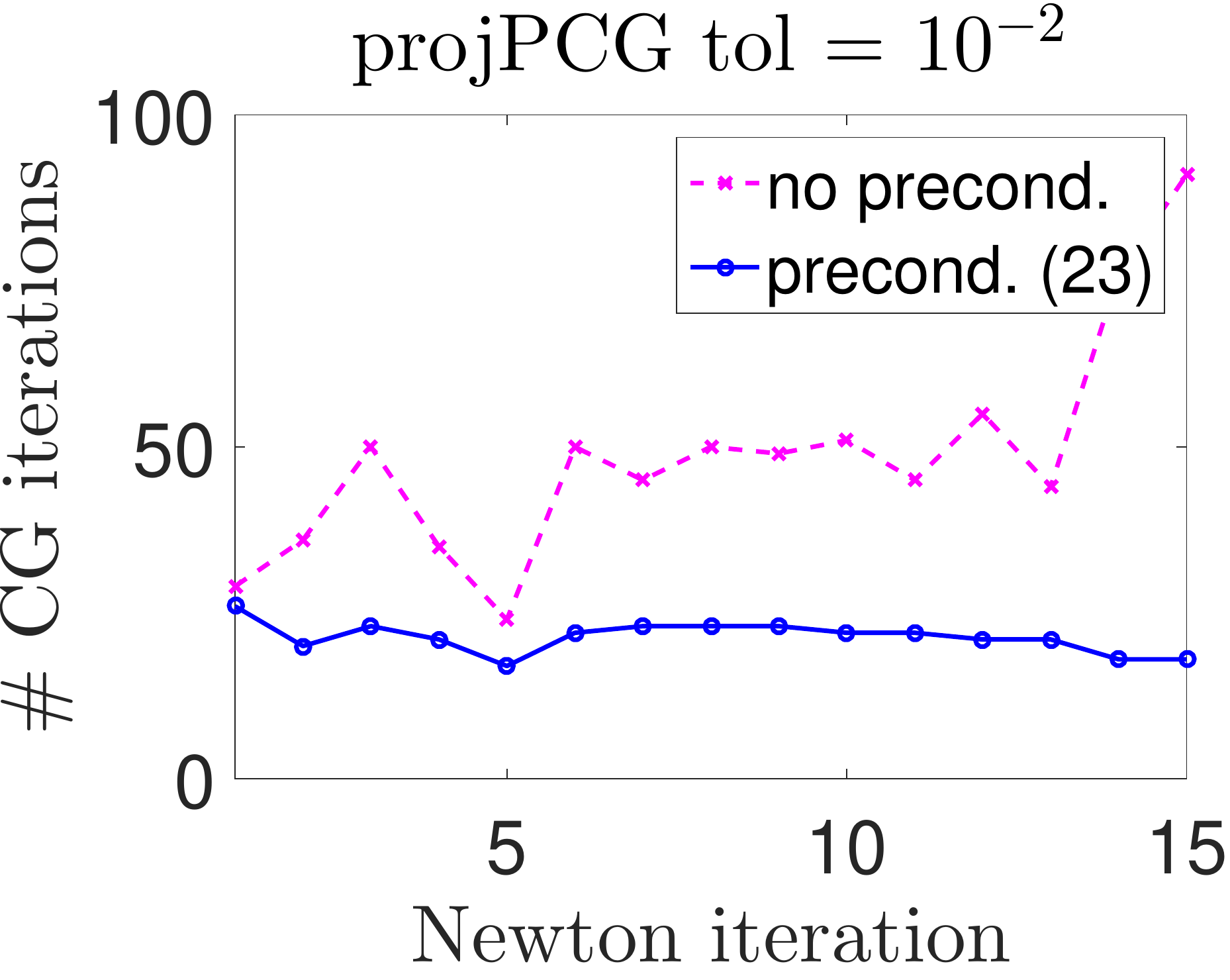} 
    \quad
	\includegraphics[width=.28\textwidth]{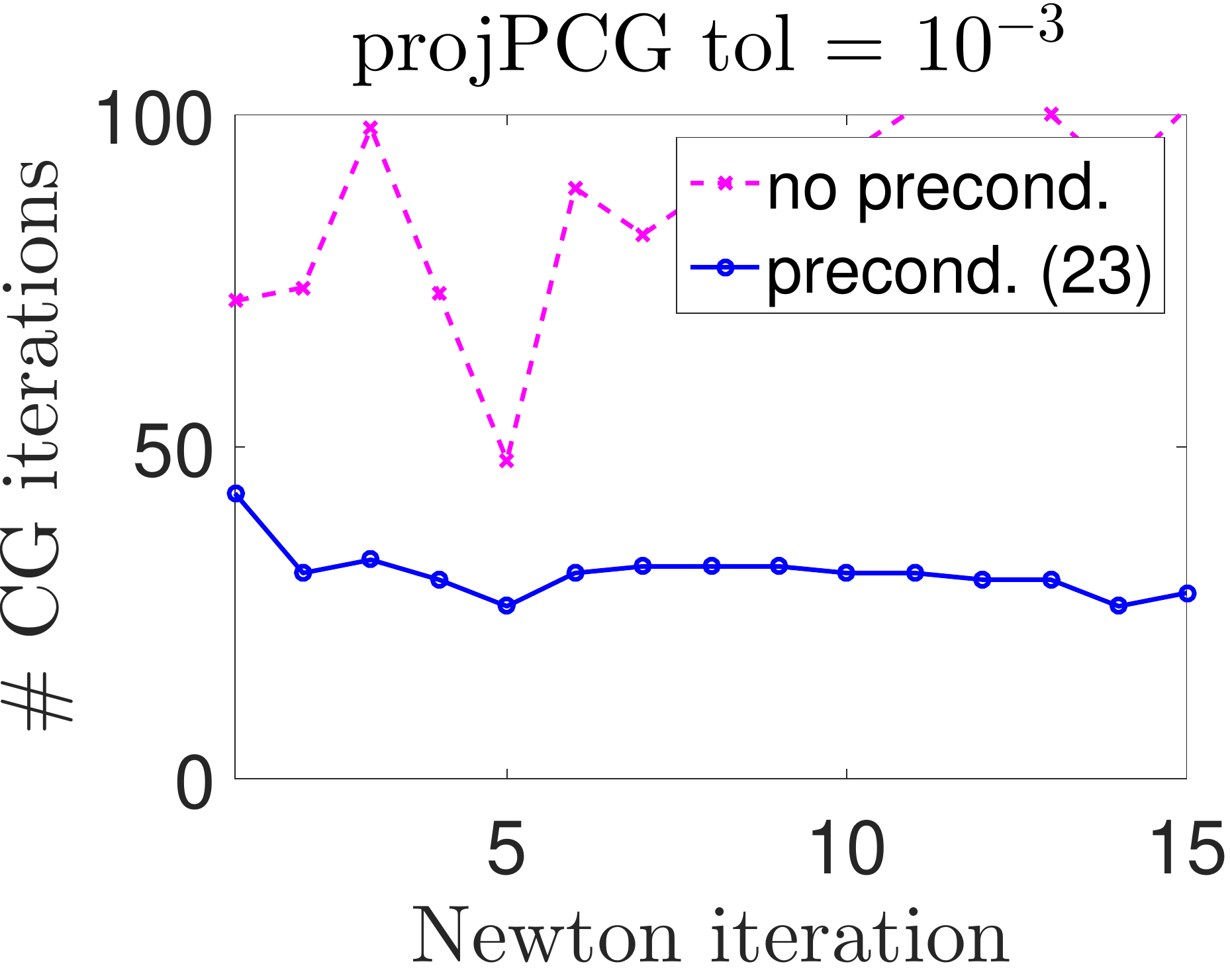} 
       	\caption{Satellite single-frame}\end{subfigure}

    \vspace*{.3cm}   	
    
   \begin{subfigure}[b]{\textwidth}
    \includegraphics[width=.28\textwidth]{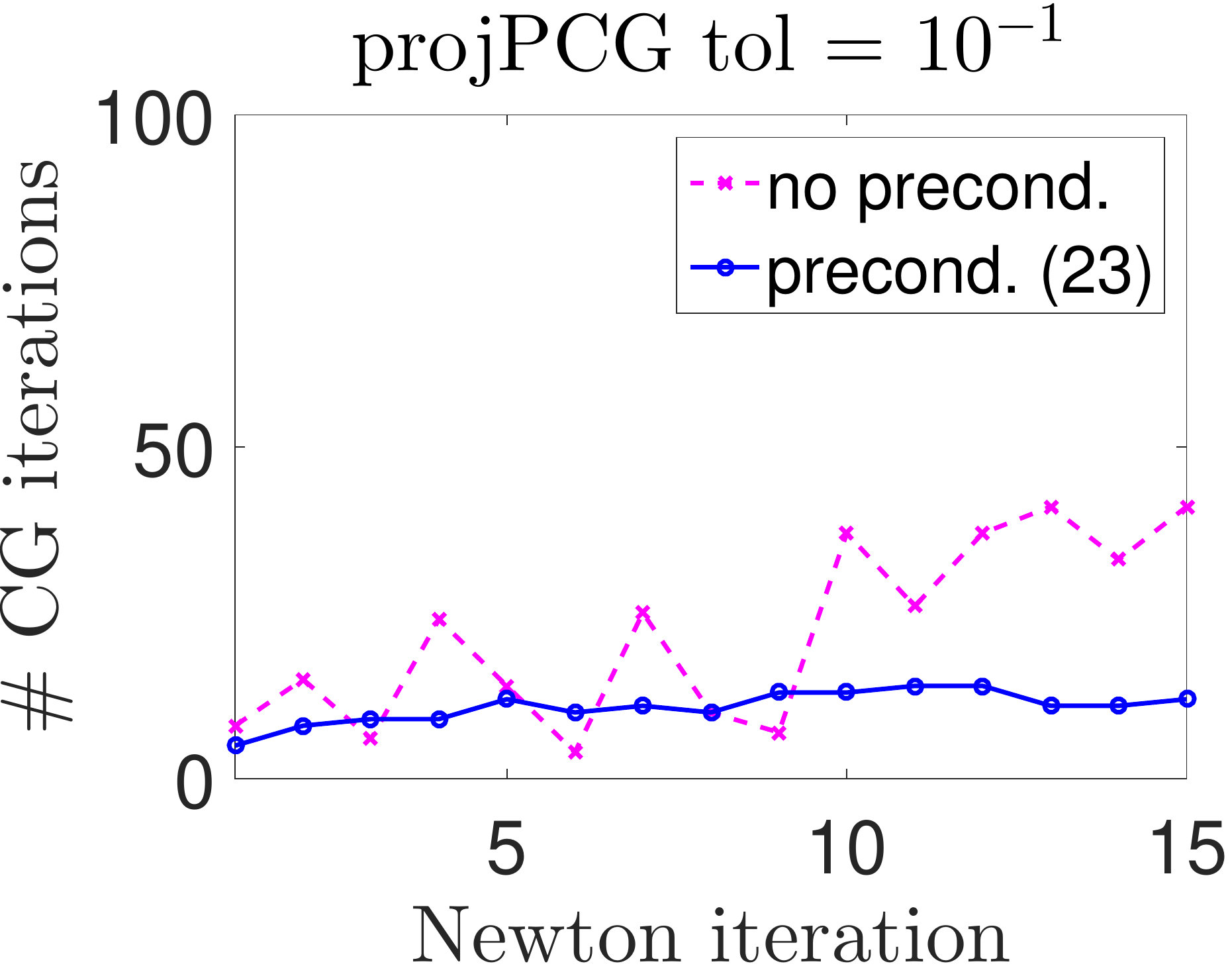} 
    \quad
	\includegraphics[width=.28\textwidth]{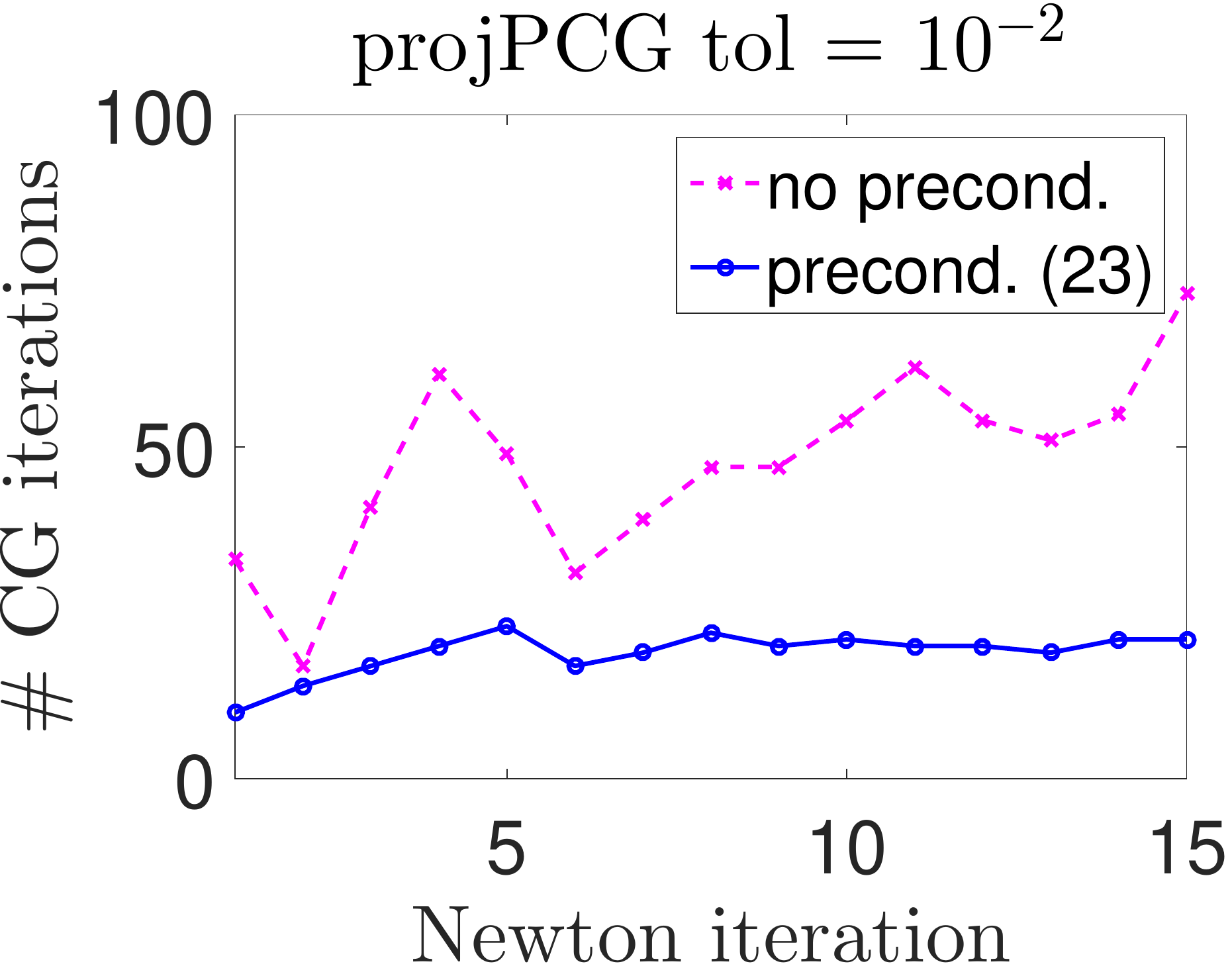}
    \quad
	\includegraphics[width=.28\textwidth]{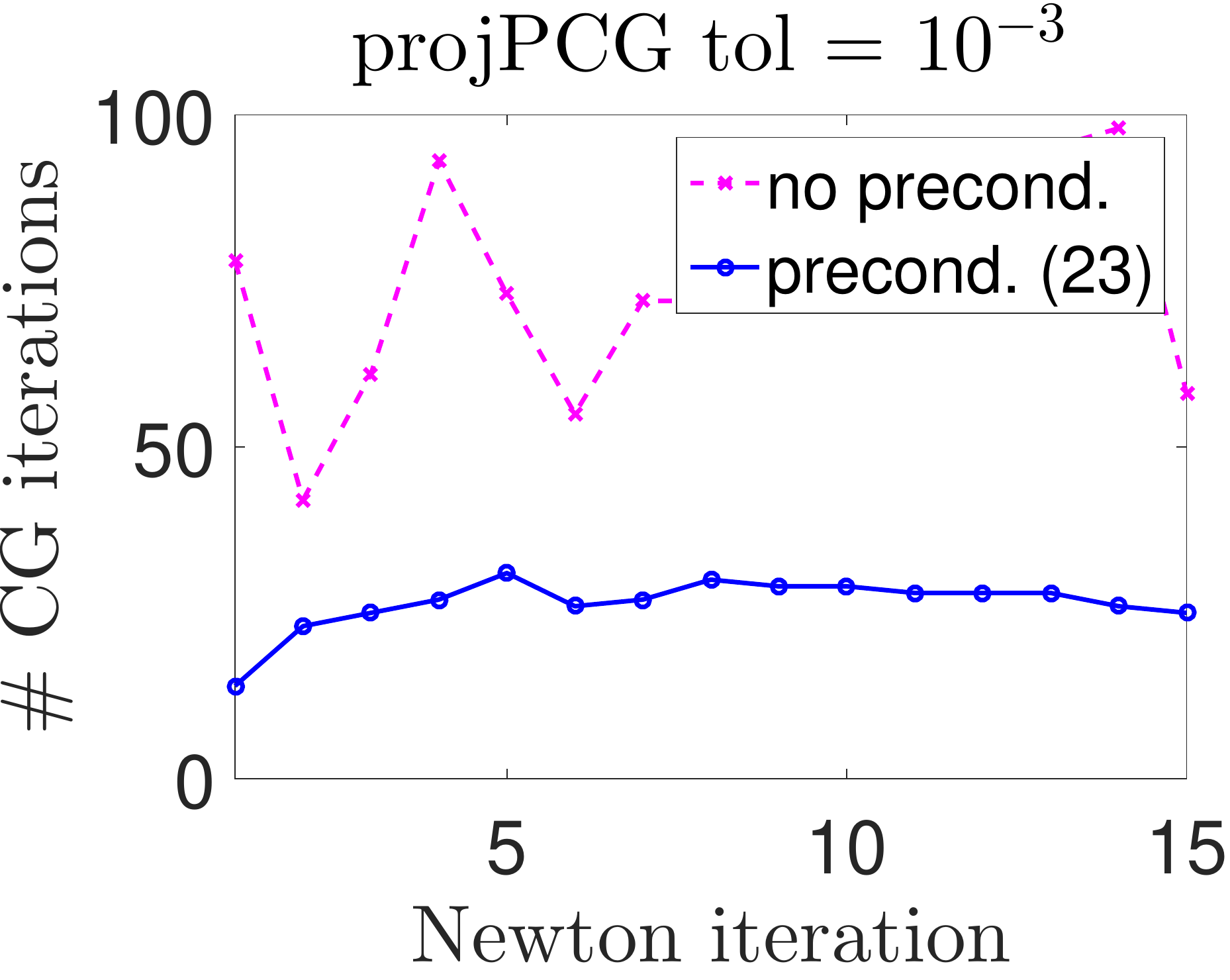} 
       	\caption{Satellite multi-frame}\end{subfigure}

    \vspace*{.3cm}   	
    
   \begin{subfigure}[b]{\textwidth}
    \includegraphics[width=.28\textwidth]{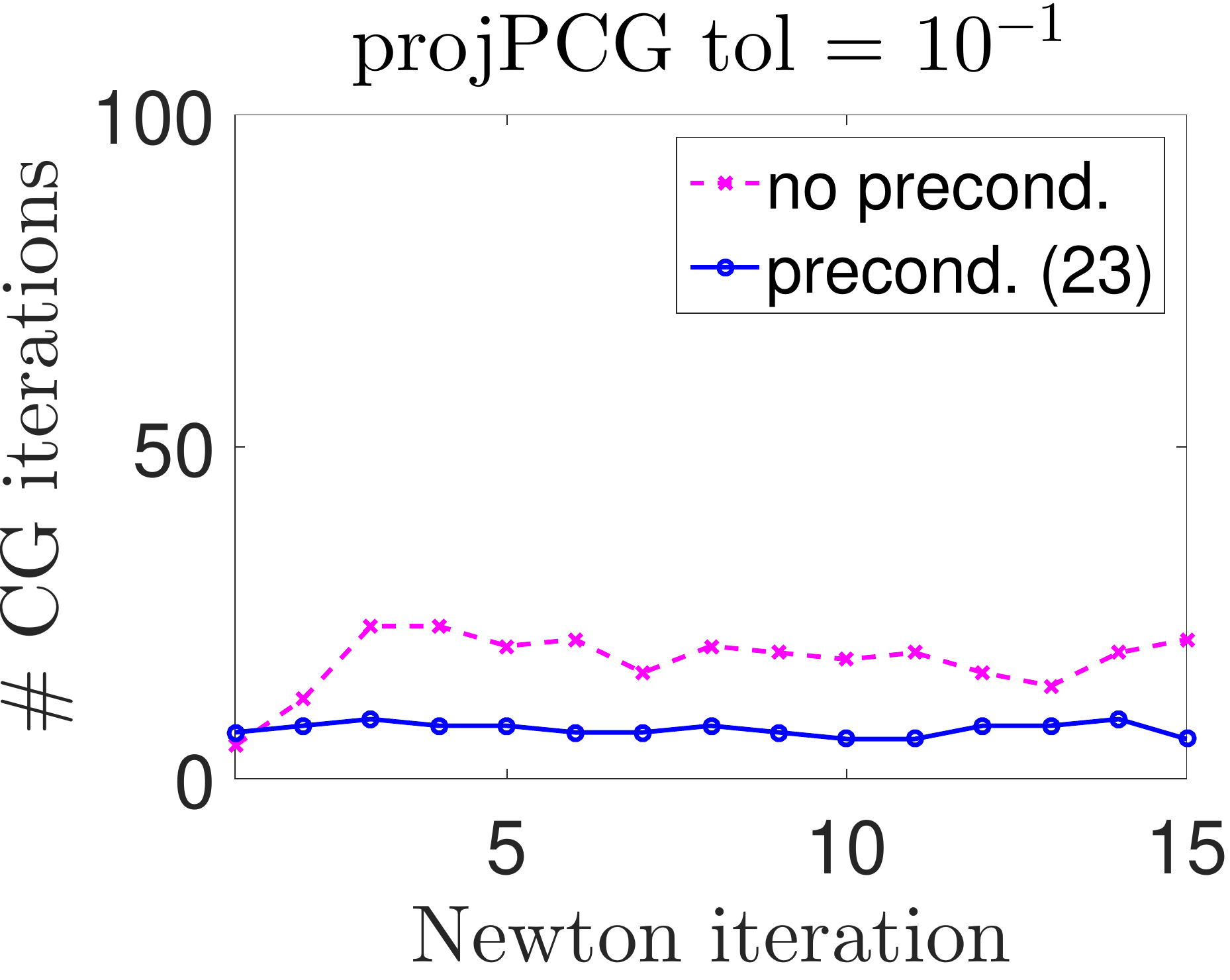} 
    \quad
	\includegraphics[width=.28\textwidth]{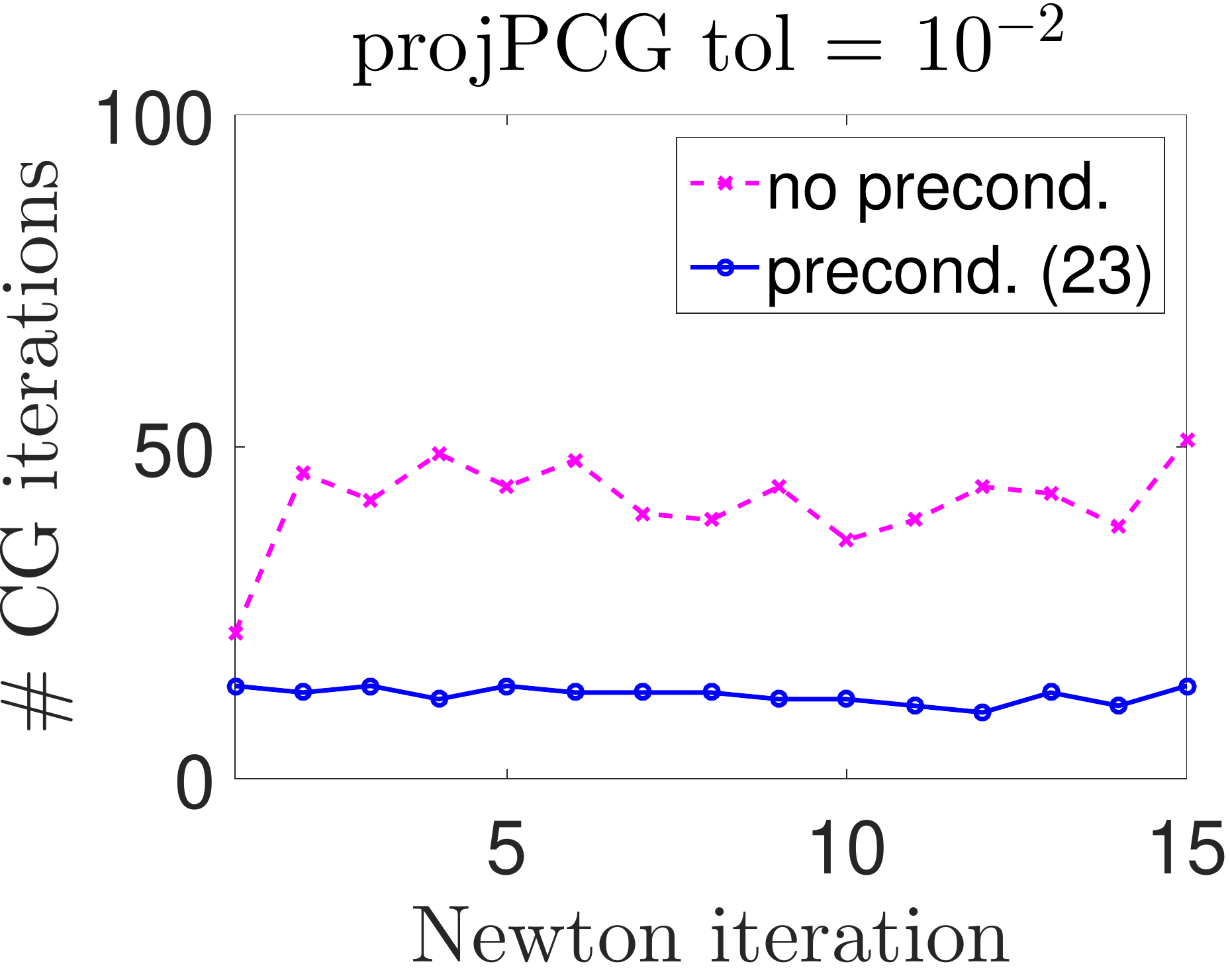} 
    \quad
	\includegraphics[width=.28\textwidth]{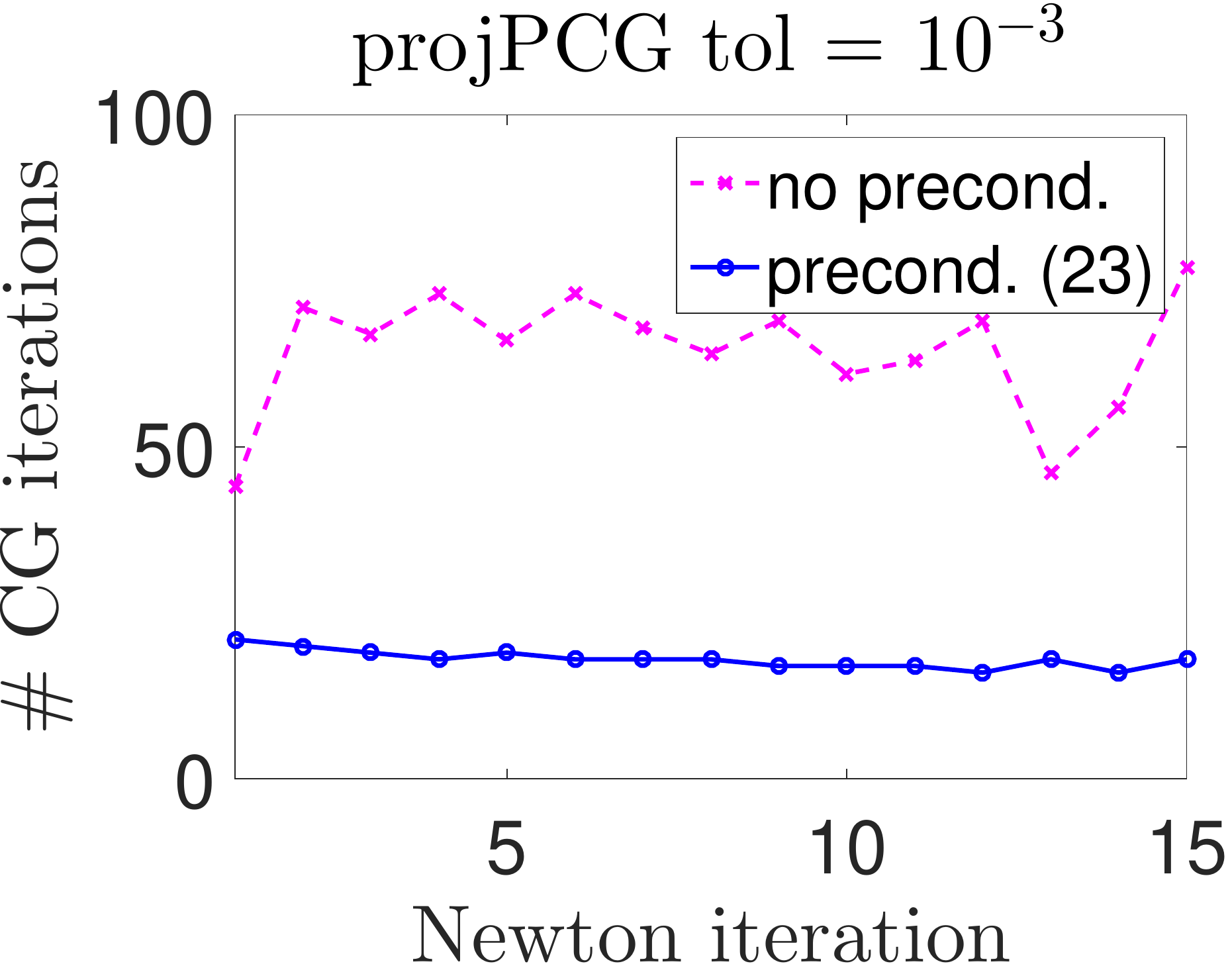} 
       	\caption{Carbon Ash single-frame}\end{subfigure}
    
    \vspace*{.3cm}   	
    
    \begin{subfigure}[b]{\textwidth}
    \includegraphics[width=.28\textwidth]{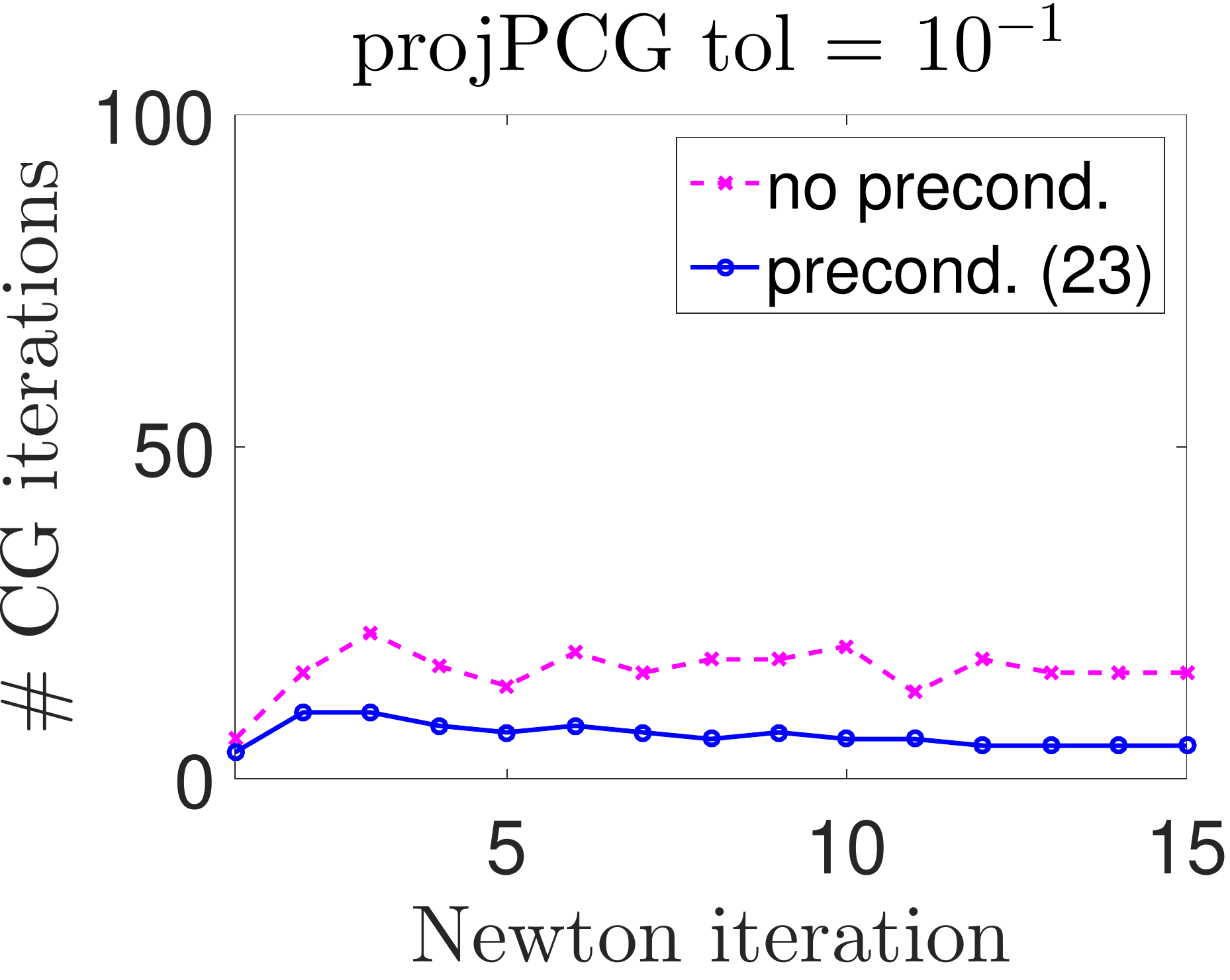} 
    \quad
	\includegraphics[width=.28\textwidth]{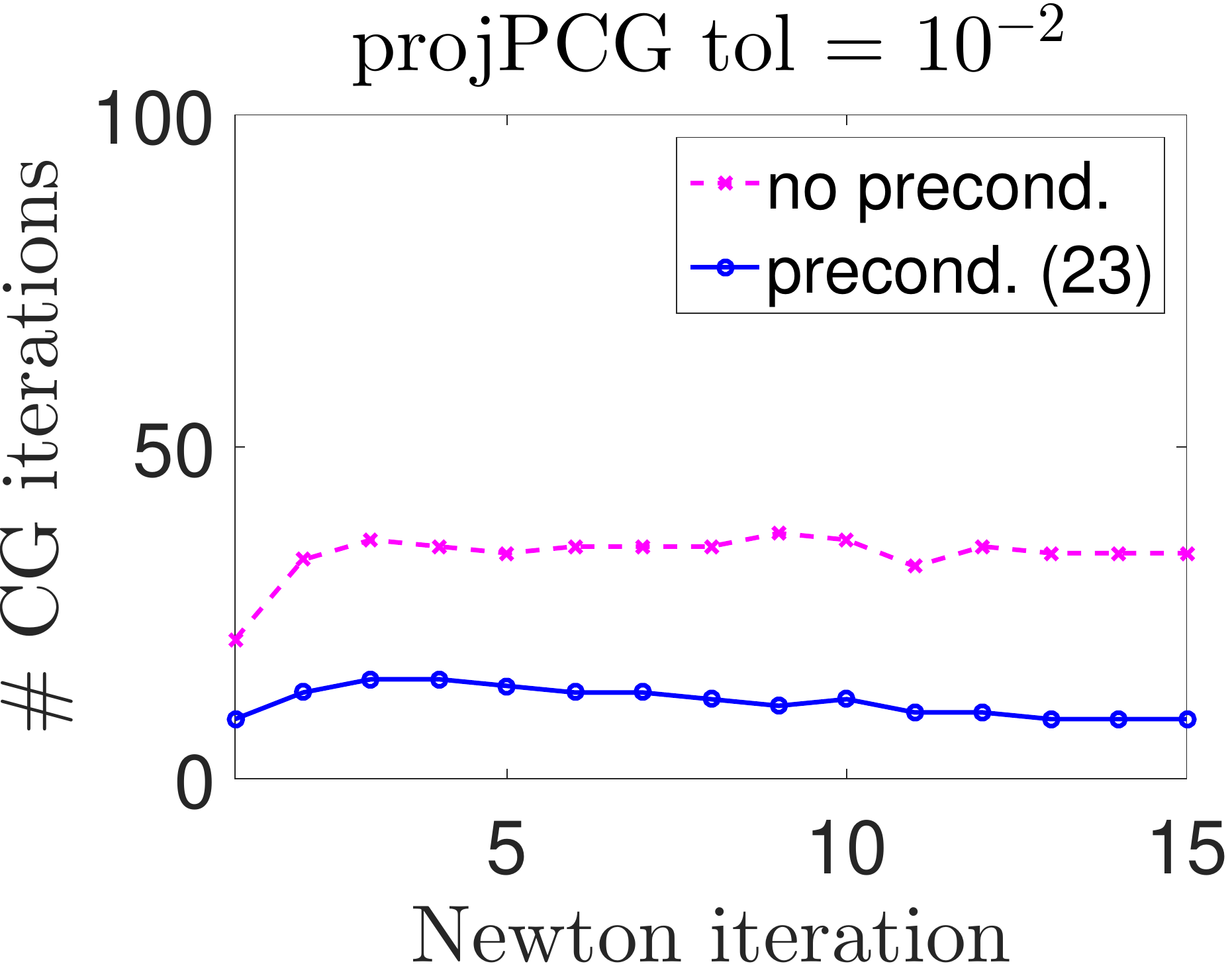} 
    \quad
	\includegraphics[width=.28\textwidth]{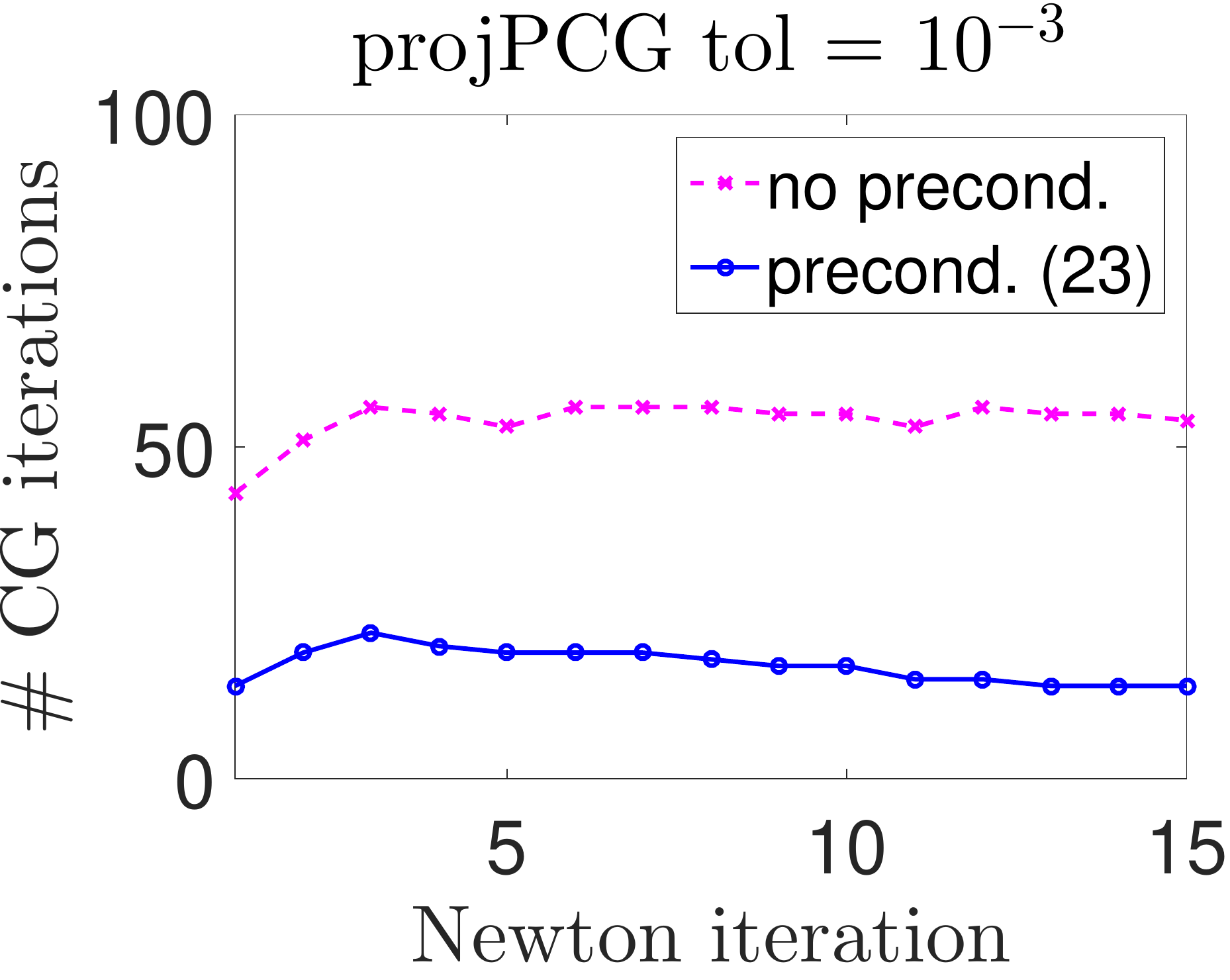} 
       	\caption{Carbon Ash multi-frame}\end{subfigure}
    
    
    	\caption{The effect of preconditioning by preconditioner defined in \eqref{eq:prec}: number of projPCG iterations performed in each Newton iteration to achieve the desired tolerance. 5 \% outliers.
} \label{fig:prec_levels}
\end{figure}

In Figure \ref{fig:prec_levels}, we investigate the overall speedup of the convergence by plotting the number of projected PCG steps needed in each Newton iteration to reach the desired tolerance on the relative size of the projected gradient. Even for the most generous tolerance $10^{-1}$, preconditioner \eqref{eq:prec} significantly reduces the number of projPCG iterations. Note that in this experiment, the linear subproblems solved in each Newton iteration are generally not identical, since the subproblems are not solved exactly and therefore the approximations $x^{(k)}$ are not the same. We set the outer tolerance to $0$ in order to perform always at least $15$ Newton iterations. 
The choice of projPCG tolerance is a difficult question, but from the average number of Newton iterations/projPCG iterations/fast Fourier transforms shown in Table \ref{tab:numit}, we observe that raising the tolerance does not considerably increase the number of Newton steps we need to perform here. Therefore larger tolerance, here $10^{-1}$, leads to a smaller total number of projPCG iterations. This is independent of the percentage of outliers. For each setting, the number of projPCG iterations is significantly smaller for the preconditioned version. This is not always the case for the total count of the fast Fourier transforms, since we need to perform 6 \texttt{fft2}/\texttt{ifft2} in each iteration vs. 4 for the unpreconditioned iterations; see Table \ref{tab:op_counts}. For large scale problems, however, the computational complexity of fast Fourier transform, which is $\mathcal{O}(n\log n)$ is comparable to other operations performed in projPCG, such as the inner products, whose complexity is $\mathcal{O}(n)$, and therefore the number of projPCG iterations seems to be the more important indicator of efficiency of the preconditioner.  Recall here that
$n$ is the number of pixels in the image, so if we have a $256 \times 256$ array of pixels, then $n = 65535$.
\begin{table}[!ht]
\centering
\caption{Average number of Newton iterations, projPCG iterations, and (inverse) 2D Fourier transforms for projPCG with and without preconditioning, and two tolerances on the relative size of the projPCG residual. Results are averaged over 10 independent realization of noise and outliers.  
}\label{tab:numit}
\begin{subtable}[h]{\textwidth}
\subcaption{projPCG tol = $10^{-1}$}
\centering
{\small\begin{tabular}{lcccc} \toprule 
 \multicolumn{5}{c}{average count: Newton/CG/\texttt{fft2}s } \\ 
& & \multicolumn{3}{c}{\% outliers} \\ 
\cmidrule(r){3-5}
\multicolumn{1}{c}{problem} & precond & \multicolumn{1}{c}{0\%} & \multicolumn{1}{c}{2\%} & \multicolumn{1}{c}{10\%} \\ \midrule
\multirow{2}{*}{Satellite single-frame}& no& 14/290/1383& 14/274/1329& 14/283/1374\\ 
 & yes& 15/161/1280& 16/172/1362& 14/158/1252\\ 
\addlinespace[.2cm]
\multirow{2}{*}{Satellite multi-frame}& no& 12/250/2398& 12/216/2104& 13/241/2364\\ 
& yes& 12/107/1545& 11/107/1535& 12/103/1507\\ 
\addlinespace[.2cm]
\multirow{2}{*}{Carbon ash single-frame}& no& 11/190/939& 10/179/891& 11/184/915\\ 
& yes& 10/71/641& 11/72/654& 13/82/753\\ 
\addlinespace[.2cm]
\multirow{2}{*}{Carbon ash multi-frame}& no& 14/221/2200& 14/219/2179& 16/254/2542\\ 
& yes& 16/88/1510& 15/85/1419& 17/99/1654\\ 
\bottomrule
\end{tabular}}
\end{subtable}

\vspace*{.5cm}

\begin{subtable}[h]{1\textwidth}
\subcaption{projPCG tol = $10^{-2}$}
\centering
{\small\begin{tabular}{lcccc} \toprule 
 \multicolumn{5}{c}{average count: Newton/CG/\texttt{fft2}s } \\ 	
& & \multicolumn{3}{c}{\% outliers} \\ 
\cmidrule(r){3-5}
\multicolumn{1}{c}{problem} & precond & \multicolumn{1}{c}{0\%} & \multicolumn{1}{c}{2\%} & \multicolumn{1}{c}{10\%} \\ \midrule
\multirow{2}{*}{Satellite single-frame}& no& 13/536/2359& 13/539/2373& 14/641/2819\\ 
& yes& 14/284/2001& 15/302/2117& 14/296/2091\\ 
\addlinespace[.2cm] 
\multirow{2}{*}{Satellite multi-frame}& no& 11/457/4082& 11/460/4130& 12/499/4511\\ 
& yes& 11/201/2536& 11/197/2499& 12/213/2747\\ 
\addlinespace[.2cm]
\multirow{2}{*}{Carbon ash single-frame}& no& 11/393/1754& 11/432/1912& 13/526/2315\\ 
& yes& 10/121/934& 11/129/1004& 13/155/1201\\ 
\addlinespace[.2cm]
\multirow{2}{*}{Carbon ash multi-frame}& no& 13/426/3813& 14/461/4127& 16/498/4467\\ 
& yes& 13/144/1954& 14/150/2038& 16/173/2359\\ 
\bottomrule
\end{tabular}}
\end{subtable}

%
\end{table} 

\section{Conclusion}
We have presented an efficient approach to compute approximate
solutions of a linear inverse problem that is contaminated with mixed Poisson--Gaussian noise, and 
when there are outliers in the measured data. 
We investigated the convexity properties of various robust regression functions and
found that the Talwar function was the best option.  We proposed a preconditioner,
and illustrated that it was more effective than other
standard preconditioning approaches on the types of problems studied in this paper.  Moreover, we showed that a variant of the GCV method
can perform well in estimating regularization parameters in robust regression.
A detailed discussion of computational
costs, and extensive numerical experiments illustrate the approach proposed in this paper is effective and efficient on
image deblurring problems.

\section{Acknowledgment}
The authors would like to thank Lars Ruthotto for pointing them to the Projected PCG algorithm and providing them with a Matlab code. The first author would like to thank Emory University for the hospitality offered in academic year 2014-2015, when part of this work was completed.

\FloatBarrier

\end{document}